\documentclass[11pt]{amsart}

\usepackage[sort,compress]{cite}
\usepackage[footnotesize,
labelfont=bf,textfont=it,up,
belowskip=-8pt,aboveskip=5pt]{caption}
\usepackage{amsaddr,amssymb,amsmath,amsthm}
\usepackage[inline]{enumitem}
\usepackage{framed,multirow}
\usepackage{geometry}
\usepackage{nicefrac}
\usepackage{latexsym}
\usepackage{url}
\usepackage{xcolor}
\usepackage{booktabs}
\usepackage{graphicx}
\usepackage{multirow}
\usepackage{siunitx}
\usepackage{algorithm,algpseudocode}

\definecolor{p1color}{RGB} {200,  16,  46} 
\definecolor{p2color}{RGB} {  0, 179, 136} 
\definecolor{p3color}{RGB} {246, 190,   0} 
\definecolor{p4color}{RGB} {136, 139, 141} 
\definecolor{p5color}{RGB} {255, 249, 217} 

\usepackage[colorlinks=true,
citecolor=p1color,
filecolor=black,
linkcolor=p1color,
urlcolor=p1color]{hyperref}

\usepackage[T1]{fontenc}
\newtheorem{remark}{Remark}

\newcounter{runid}
\newcommand\resetrid{\setcounter{runid}{0}}
\newcommand\runid{\addtocounter{runid}{1}\therunid}

\newcommand\fnum[1]{\num[scientific-notation=fixed,round-precision=1,round-mode=places]{#1}}
\newcommand\snum[1]{\num[round-precision=1,round-mode=places,tight-spacing=true,retain-zero-exponent=true,retain-explicit-plus=true]{#1}}
\newcommand\inum[1]{\num[scientific-notation=fixed,round-precision=0,round-mode=places]{#1}}

\DeclareMathOperator*{\minopt}{minimize}
\DeclareMathOperator*{\argmin}{argmin}
\DeclareMathOperator{\prox}{prox}

\newcommand{\defeq}{\ensuremath{\mathrel{\mathop:}=}}

\newcommand{\kin}{\ensuremath{\text{kin}}}
\newcommand{\diff}{\ensuremath{\text{diff}(\mathbb{R}^3)}}
\newcommand{\reg}{\ensuremath{\text{reg}}}
\newcommand{\dist}{\ensuremath{\text{dist}}}

\newcommand{\app}{\ensuremath{\!\cdot\!}}
\renewcommand{\d}{\ensuremath{\text{d}}}
\newcommand{\id}{\ensuremath{\text{id}_{\mathbb{R}^3}}}

\newcommand\secref[1]{\S\ref{#1}}
\newcommand\figref[1]{Fig.~\ref{#1}}
\newcommand\tabref[1]{Tab.~\ref{#1}}
\newcommand\ialgref[1]{Alg.~\ref{#1}}
\newcommand\idef[1]{\textbf{#1}}

\title[An Operator-Splitting Approach for Diffeomorphic Shape Matching]{An Operator-Splitting Approach for Variational Optimal Control Formulations for Diffeomorphic Shape Matching}

\author{Andreas Mang, Jiwen He \and Robert Azencott}
\email{andreas@math.uh.edu (AM), jiwenhe@math.uh.edu (JH), razencott@uh.edu (RA)}
\address{Department of Mathematics, University of Houston, 3551 Cullen Blvd, Houston, TX 77204-3008, USA}

\begin{document}

\maketitle

\begin{abstract}
We present formulations and numerical algorithms for solving diffeomorphic shape matching problems. We formulate shape matching as a variational problem governed by a dynamical system that models the flow of diffeomorphism $f_t \in \diff$. We overview our contributions in this area, and present an improved, matrix-free implementation of an operator splitting strategy for diffeomorphic shape matching. We showcase results for diffeomorphic shape matching of real clinical cardiac data in $\mathbb{R}^3$ to assess the performance of our methodology.
\end{abstract}

\vspace{0.2cm}
{\noindent \it Dedicated to the memory of Roland Glowinski.}

\section{Introduction}

In (bio)medical imaging applications, the automatic matching of $k$-dimensional deformable shapes across subjects or multi-temporal data of individual patients is a critical step to aid clinical diagnosis~\cite{Modersitzki:2004a,Modersitzki:2009a,Fischer:2008a,Hajnal:2001a,Sotiras:2013a,Younes:2019a,Hill:2001a,Toga:2001a}. From a mathematical point of view, this matching problem constitutes an \emph{inverse problem}~\cite{Fischer:2008a}: Given two (or more) shapes $s_i\in \mathcal{S}$, $i=0,1$, in $\mathbb{R}^d$ representing an object/anatomy of interest, defined in some shape space $\mathcal{S}$, we seek a \emph{plausible} spatial transformation $y \in \mathcal{Y}_{\text{ad}}$, $\mathcal{Y}_{\text{ad}} \subseteq \{\phi \mid \phi: \mathbb{R}^d \to \mathbb{R}^d\}$ with $d \in \{2,3\}$ that establishes a point-wise correspondence between these objects~\cite{Modersitzki:2004a,Modersitzki:2009a}. Based on the problem specific notion of plausibility of an admissible transformation $y$, various continuum mechanical models emerged to restrict the space of admissible transformations $\mathcal{Y}_{\text{ad}}$. In a variational setting, this prior knowledge on admissible maps $y$ is typically prescribed through regularization functionals or constraints. Examples include models of elasticity~\cite{Burger:2013a,Kybic:2003a,Davatzikos:1997a} or incompressibility~\cite{Haber:2004a,Modersitzki:2008a,Mansi:2011a,Mang:2015a,Mang:2016a}.

The preferred model typically depends on the application, of which there are many in medical imaging. Excellent reference works giving an overview of developments in this field are~\cite{Modersitzki:2004a,Modersitzki:2009a,Fischer:2008a,Hajnal:2001a,Sotiras:2013a,Younes:2019a,Hill:2001a}. First, one may simply be interested in relating local information observed in two views of the same object. This can be helpful in assessing, for example, changes in texture or morphology that occurred over time due to the progression of a disease or treatment~\cite{Mang:2008a,Toga:2001a,Misra:2009a}. Moreover, one can integrate data acquired by multiple sensors (different imaging modalities) into a common reference frame (between imaging sessions, the patient/organs move).  Likewise, one can compensate for organ motion associated with breathing or the heart cycle~\cite{Guigui:2022a,Makela:2002a,Gorbunova:2012a,Ehrhardt:2010a,Christensen:2007a}. Registration has also been used for the construction of atlases~\cite{Toga:2001a,Lee:2017a,Joshi:2004a,Campbell:2021a,Zhang:2013a,Serag:2012a} and atlas-based segmentation~\cite{Cabezas:2011a,Vemuri:2003a,Rohlfing:2004a}. Aside from these examples of integrating data into a common frame of reference, one can use these tools to study changes in the (morphological) appearance of organs. This observation has lead to the notion of \emph{computational anatomy}~\cite{Ashburner:2009a,Thompson:2002a,Guigui:2022a,Grenander:1998a,Miller:2004a,Miller:2002a,Younes:2009a,Miller:2015a}. The results presented in this work fall into this category. Here, one uses the computed maps $y$ (typically, diffeomorphisms) to derive features that might be useful to classify or characterize patients (e.g., diseased versus healthy)~\cite{Avants:2008a,Risser:2010a,Wang:2007a,Fox:1997a,Hua:2008a,Joy:2023a,Misra:2009a,Davatzikos:2001a}. These features include geodesic distances, strain values, or local volume change associated with the computed map $y$ (we specify this mathematically more rigorously below). One can either compare patients to one another~\cite{Dabirian:2022a}, or, likewise, individual patients to a standard reference template~\cite{Lee:2017a,Joshi:2004a,Zhang:2013a} and then---based on the derived features---try to classify them. Similarly, one can monitor a patient over time and use the computed maps between time frames to derive features that can be used to assess the progression or existence of a pathology. We consider both applications in the present work.

We limit $\mathcal{Y}_{\text{ad}}$ to infinite dimensional groups of diffeomorphisms $\mathcal{G}\subseteq \diff$ in $\mathbb{R}^3$. This approach provides a rich mathematical framework for studying and quantifying shape variability in $\mathbb{R}^3$ through the lens of geodesics in the group of diffeomorphisms $\mathcal{G}$. In the context of biomedical imaging, this area of research is---as mentioned above---typically referred to as computational anatomy~\cite{Grenander:1998a,Miller:2004a,Miller:2002a,Younes:2009a,Miller:2015a}. The underlying mathematical framework for diffeomorphic shape matching is referred to as \emph{large deformation diffeomorphic metric mapping} (\idef{LDDMM})~\cite{Trouve:1998a,Trouve:1995a,Glaunes:2004a,Glaunes:2008a,Younes:2019a,Beg:2005a,Joshi:2000a} and builds upon the seminal works~\cite{Arnold:1966a,Arnold:1976a,Ebin:1970a}. In~\cite{Arnold:1966a,Arnold:1976a,Ebin:1970a} it was shown that in an incompressible fluid, if $f_t(x)$, $f_t(x) \defeq f(t,x)$, is the position at time $t$ of a particle originating from initial location $x$, then the map $t \to f_t(x)$ defines a geodesic in $\mathcal{G}$ for the metric defined by the kinetic energy of the fluid (see \eqref{e:kinv}). This observation has been used in numerous works to design variational optimal control formulations that allow us to match two shapes $s_0$ and $s_1$ by computing the ``shortest'' geodesic in $\mathcal{G}$ that connects them, i.e., maps $s_0$ to $s_1$.

\subsection{Contributions}

We provide an exposition of our work in the area of diffeomorphic shape matching, including the design of numerical algorithms~\cite{Azencott:2008a,Azencott:2010a,Freeman:2014a,Jajoo:2011a,Qin:2013a,Zhang:2019a,Zhang:2021a} and their application to clinical problems~\cite{Zekry:2012a,Zekry:2016a,Zekry:2018a,ElTallawi:2021a,ElTallawi:2019a,Dabirian:2022a}. In addition, we present an improved numerical implementation of the operator-splitting algorithm originally proposed in~\cite{Zhang:2021a,Zhang:2019a}. The work in~\cite{Zhang:2021a,Zhang:2019a} uses direct methods~\cite{Davis:2006a,Duff:2017a} to solve the optimality systems, i.e., matrices that appear in our optimality systems are formed and stored explicitly. Our new implementation is matrix-free. That is, we either implement matrix-vector products instead of forming matrices and/or exploit the block-structure of the matrix operators to reduce the memory footprint. Moreover, we consider a Schur complement (reduced-space) method~\cite{Nocedal:2006a} and use iterative algorithms to solve the individual subproblems that appear in our optimality system. Overall, we obtain a speedup of approximately 2$\times$ compared to our prior work~\cite{Zhang:2021a,Zhang:2019a}. In addition, the matrix-free implementation reduces the memory footprint of our method, which allows us to solve larger problems compared to our prior work.

\subsection{Related Work}

We present numerical algorithms for solving optimal control problems for diffeomorphic shape matching and their application to cardiac imaging. We refer the reader to~\cite{Modersitzki:2004a,Sotiras:2013a,Modersitzki:2009a,Younes:2019a} for a general introduction into numerical methods, applications, and formulations for image registration and diffeomorphic shape matching. There exists a plurality of mathematical and numerical frameworks to study shape variability; examples include~\cite{Bauer:2014a,Tycowicz:2018a,Ludke:2022a,Heimann:2009a,Ambellan:2019a,Iglesias:2018a,Kilian:2007a,Torres:2004a,Bone:2020a,Junior:2018a,Rustamov:2013a,Srivastava:2010a,Fletcher:2004a}. The work presented in this manuscript builds upon the rich mathematical framework developed in~\cite{Bauer:2014a,Arguillere:2015a,Dupuis:1998a,Trouve:1998a,Arnold:1966a,Arnold:1976a,Ebin:1970a,Trouve:1995a,Glaunes:2004a,Glaunes:2008a,Younes:2019a,Beg:2005a}. In this framework, we define the shape space $\mathcal{S}$ as orbits of given template shapes under the action of the group of diffeomorphisms $\mathcal{G}$. Our past work for the matching of smooth surfaces has been described in~\cite{Azencott:2008a,Azencott:2010a,Freeman:2014a,Jajoo:2011a,Qin:2013a,Zhang:2021a,Zhang:2021a}. Likewise to~\cite{Trouve:1998a,Trouve:1995a,Glaunes:2004a,Glaunes:2008a,Younes:2019a,Beg:2005a}, we consider a variational approach for solving for the diffeomorphism $y \in \mathcal{G}$ that maps one shape $s_0 \in \mathcal{S}$ to another shape $s_1 \in \mathcal{S}$ such that $y \app s_0 \approx s_1$. These shapes $s_i$, $i=0,1$, are typically represented as curves, lines, or surfaces in $\mathbb{R}^d$ with $d\in\{2,3\}$~\cite{Kaltenmark:2017a,Azencott:2010a}.

Our work in this area was initiated by a short exposition of the mathematical problem formulation~\cite{Azencott:2008a}. This joint paper of Azencott and Glowinski represents the foundation for a series of papers of our group, spanning more than one decade~\cite{Azencott:2008a,Azencott:2010a,Freeman:2014a,Jajoo:2011a,Qin:2013a,Zhang:2021a,Zekry:2018a,ElTallawi:2021a,ElTallawi:2019a,Dabirian:2022a}. In~\cite{Azencott:2010a} Azencott, Glowinski, He, and co-authors build up on~\cite{Azencott:2008a} and develop a first, efficient gradient descent algorithm for the considered optimal control formulation~\cite{Azencott:2008a}. They applied their method to datasets of mitral valves, initiating a long-term collaboration between the Department of Mathematics at the University of Houston and the Houston Methodist Hospital. This collaboration has resulted in several joint publications over the years~\cite{Zekry:2012a,Zekry:2016a,Zekry:2018a,ElTallawi:2021a,ElTallawi:2019a,Dabirian:2022a,Zhang:2021a}. Improvements to the gradient descent (adjoint-based) methods and associated numerical strategies originally proposed in~\cite{Azencott:2010a} led to several doctoral research projects~\cite{Freeman:2014a,Jajoo:2011a,Qin:2013a,Zhang:2019a}. This work culminated in the development of an operator-splitting strategy to solve the associated control problem~\cite{Zhang:2021a}---a numerical technique that is equivalent to the \emph{Douglas--Rachford} splitting method~\cite{ODonghue:2013a,Bauschke:2011a,Glowinski:2016b,Glowinski:2016a,Glowinski:2016c,Bukac:2016a}. This approach can be traced back to the seminal works of Glowinski and Marrocco~\cite{Glowinski:1975a} and Gabay and Mercier~\cite{Gabay:1976a}. The present work describes our current numerical implementation of this approach, and showcases its application to real world problems in biomedical imaging.

Related work on adjoint-based, variational approaches for diffeomorphic registration of images by our group is described in~\cite{Mang:2015a,Mang:2016a,Mang:2016b,Mang:2017a,Mang:2017b,Mang:2019a,Brunn:2020a,Brunn:2021a,Brunn:2021b,Himthani:2022a}. As opposed to the work considered here, these approaches directly operate on image intensities rather than shape representations derived from these data. Other work on numerical algorithms for the solution of variational shape matching problems (for images as well as surface representations) of other groups are, e.g., described in~\cite{Beg:2005a,Cao:2005a,Hsieh:2022a,Polzin:2016a,Polzin:2020a,Vialard:2012a,Niethammer:2009a,Hart:2009a,Arguillere:2016a}. The majority of existing works consider an \emph{optimize-then-discretize} approach to solve the variational problem~\cite{Beg:2005a,Arguillere:2016a,Mang:2015a,Mang:2016a,Mang:2016b,Vialard:2012a,Zhang:2015b,Mang:2017a,Ashburner:2011a,Hart:2009a,Miller:2006a}. Aside from our work~\cite{Azencott:2010a,Mang:2017a,Zhang:2021a}, \emph{discretize-then-optimize} approaches for LDDMM can be found in~\cite{Polzin:2020a,Polzin:2018a}. The problem formulation as well as our discretization naturally leads to a reproducing kernel Hilbert space ({\bf RKHS}) structure~\cite{Aronszajn:1950a}. Exploiting an RKHS representation of smooth functions in the context of diffeomorphic shape matching (and associated problems) is a common strategy considered in numerous works~\cite{Jain:2013a,Miller:2015a,Richardson:2016a,Hsieh:2021a,Arguillere:2016a,Glaunes:2004a}. As for the numerical solution of the variational optimization problem, many existing works consider (first order) gradient descent-type optimization algorithms; see, e.g.,~\cite{Beg:2005a,Cao:2005a,Hart:2009a,Vialard:2012a,Hsieh:2020a}. Higher order (quasi-)Newton algorithms have, e.g., been developed in~\cite{Ashburner:2011a,Mang:2015a,Mang:2016a,Mang:2019a,Polzin:2016a}. Related work on operator-splitting algorithms for LDDMM (and related problems) can be found in~\cite{Thorley:2021a,Lee:2016a}. Other recent works that do not explicitly derive optimality conditions based on variational principles but rely on automatic differentiation can be found in~\cite{Hsieh:2021a,Franccois:2021a,Bone:2018a,Hartman:2023a,Bone:2020a}. Lastly, we note that the success of machine learning in various scientific disciplines has led to several recent works that attempt to solve the inverse problem of diffeomorphic registration based on machine learning techniques~\cite{Shen:2021a,Bone:2020a,Tian:2022a,Amor:2021a,Krebs:2019a,Sun:2022a,Yang:2017a,Wu:2023a,Bharati:2022a,Huang:2021a,Wu:2022a}.

\subsection{Outline}

We present our mathematical formulation and our numerical approaches in~\secref{s:methods}. The parameterization of diffeomorphic matching in terms of a time dependent velocity field $v$ is presented in~\secref{s:flows}. The associated ODE constitutes the state equation for our nonlinear control problem. The variational problem formulation is described in~\secref{s:varopt}. We present the numerical discretization in~\secref{s:discretization}. We describe our approach for numerical optimization in~\secref{s:optimization}. We list quantitative measures of shape variability derived from diffeomorphic shape matching in~\secref{s:shapevariability}. We report numerical experiments in~\secref{s:experiments} and conclude with~\secref{s:conclusions}.

\section{Mathematical Treatment and Numerical Methods}\label{s:methods}

In the following sections, we describe the overall mathematical formulation and different numerical strategies to approach the problem. For easier access to the material, we provide an overview of our notation in \tabref{t:notation}.

\begin{figure}
\centering
\includegraphics[totalheight=3cm]{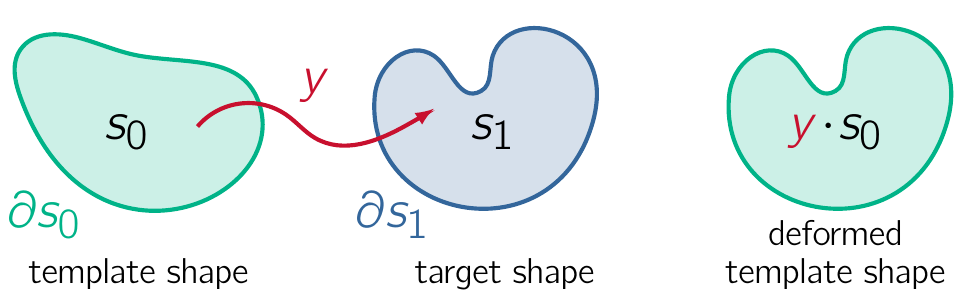}
\caption{Illustration of the diffeomorphic shape matching problem. We display two shapes $s_0, s_1 \in \mathcal{S}$ with smooth boundaries $\partial s_0$, $\partial s_1$. We seek a diffeomorphic map $y$ that, if applied to $s_0$, makes the deformed template shape $y \app s_0$ similar to the target shape $s_1$, i.e., $y \app s_0 = y(s_0) = s_1$.\label{f:diffeomorphic-registration}}
\end{figure}

Diffeomorphic shape matching can be formulated as follows: Given two $k$-dimensional shapes $s_0 \in \mathcal{S}$ and $s_1 \in \mathcal{S}$ in $\mathbb{R}^3$ with $k \in \{1,2,3\}$ we seek an $\mathbb{R}^3$-diffeomorphism $y \in \diff$ such that $y \app s_0 = y(s_0) = s_1$. This problem is ill-posed; it requires regularization for a numerical treatment~\cite{Engl:1996a}. In the following, we describe the mathematical formulation and different numerical strategies to solve the variational optimization problem presented below. One can either match biomedical imaging data by diffeomorphic image registration operating directly on image intensities (see, e.g.,~\cite{Vialard:2012a,Ashburner:2011a,Hart:2009a,Beg:2005a,Mang:2016b,Mang:2019a,Mang:2015a,Mang:2017a}), or one can match shape representations extracted from these data (see \figref{f:diffeomorphic-registration}). Here, we focus on the latter approach. In the present work, these shapes are represented as ``surfaces.'' The term ``surface'' denotes a compact smooth 2-dimensional manifold with boundary, smoothly embedded in $\mathbb{R}^3$. The boundary of such a surface is a smooth 1-dimensional manifold smoothly embedded in $\mathbb{R}^3$, and---in all the cases considered here---is a connected finite union of smooth curves. We note that the general framework presented here applies (and has been applied) to curves, lines, or boundaries of solids embedded in $\mathbb{R}^3$. We parametrize the considered smooth surfaces with boundaries as finite sets of points in $\mathbb{R}^3$ (see \secref{s:discretization}).

\begin{table}
\caption{Notation and symbols.\label{t:notation}}
\centering
\begin{small}
\begin{tabular}{ll}
\toprule
\bf Variable/Acronym & \bf Meaning \\
\midrule
$\mathcal{S}$ & shape space \\
$\mathcal{H}$ & Hilbert space of $\mathbb{R}^3$ vector fields \\
$\mathcal{H}_{\sigma}^d$ & RKHS of $\mathbb{R}^d$ vector fields parametrized by $\sigma > 0$ (we drop $d$ for $d=1$)\\
$\text{diff}(\mathbb{R}^3)$ & space of diffeomorphic mappings from $\mathbb{R}^3$ to $\mathbb{R}^3$  \\
$\mathcal{V}$ & space of $L^2$-integrable (in time) smooth $\mathbb{R}^3$ vector fields\\
$\mathcal{G}$ & group of diffeomorphisms; $\mathcal{G} \subseteq \text{diff}(\mathbb{R}^3)$\\
$\id$ & identity map in $\mathbb{R}^3$; $\id(x) = x$ for any $x \in \mathbb{R}^3$\\
$s_0 \in \mathcal{S}$ & template shape (shape to be deformed) \\
$s_1 \in \mathcal{S}$ & target/fixed shape \\
$x_0 = (x_0^i)_{i=1}^{m_0} \in \mathbb{R}^{3m_0}$ & vector of mesh points to represent $s_0$ \\
$x_1 = (x_1^i)_{i=1}^{m_1}\in \mathbb{R}^{3m_1}$ & vector of mesh points to represent $s_1$ \\
$m_i \in \mathbb{N}$ & number of points for the parameterization of $s_i$, $i=0,1$ \\
$n \in \mathbb{N}$ & number of cells to discretize time interval $[0,1]$ \\
$v : [0,1] \times \mathbb{R}^3 \to \mathbb{R}^3$ & time-dependent velocity field; $v_t\defeq v(t,\cdot\,)$ \\
$f : [0,1] \times \mathbb{R}^3 \to \mathbb{R}^3$ & diffeomorphic flow; $f_t\defeq f(t,\cdot\,)$ \\
$y \defeq f_1$ $y : \mathbb{R}^3 \to \mathbb{R}^3$ & diffeomorphic map at final time $t=1$ \\
$t \in [0,1]$ & pseudo-time variable\\
$a \in \mathbb{R}^{3m_0(p+1)}$ & discrete control variable (control vector)\\
$x \in \mathbb{R}^{3m_0(p+1)}$ & discrete state variable (state vector)\\
$u \in \mathbb{R}^{3m_0(p+1)}$ & discrete dual variable corresponding to the state vector $x$  \\
$w \in \mathbb{R}^{3m_0(p+1)}$ & discrete dual variable corresponding to the control vector $w$ \\
$\operatorname{dist}:\mathcal{S} \times \mathcal{S} \to \mathbb{R}$ & shape similarity measure \\
$\operatorname{kin}:\mathcal{V} \to \mathbb{R}$ & kinetic energy/regularization functional\\
$I_n = \operatorname{diag}(1,\ldots,1) \in \mathbb{R}^{n,n}$ & $n \times n$ identity matrix \\
$e_n = (1,\ldots,1) \in \mathbb{R}^{n}$ & 1-vector of size $n \times 1$ \\
$\otimes$ & Kronecker product (see \eqref{e:defkronprod})\\
$\odot$ & Hadamard product (see \eqref{e:defhadaprod}) \\
\midrule
ADMM & alternating direction method of multipliers\\
KKT & Karush--Kuhn--Tucker (conditions)\\
LDDMM & large deformation diffeomorphic metric mapping\\
RKHS & reproducing kernel Hilbert space \\
PCG & preconditioned conjugate gradient method \\
\bottomrule
\end{tabular}
\end{small}
\end{table}

\subsection{Flows of Diffeomorphisms}\label{s:flows}

Let $k \defeq (k_1,k_2,k_3) \in \mathbb{N}^3$ denote a multi-index and $\partial^k$ denote the differential operator of order $|k| = \sum_{i=1}^3 k_i$. Moreover, let $q \in \mathbb{N}$, $1 \leq p \leq \infty$, and $u : \mathbb{R}^3 \supseteq \Omega \to \mathbb{R}$. Let
\[
C^q(\Omega)
 \defeq
\left\{u : \Omega \to \mathbb{R} : \partial^k u \text{ is continous for }|k| \leq q \right\}
\]

\noindent denote the space of $q$-times continuously differentiable functions on $\Omega \subseteq \mathbb{R}^3$. Moreover, let
\[
W^{q,p}(\Omega)
\defeq
\left\{
u \in L^p(\Omega) : \partial^k u \in L^p(\Omega) \text{ for } 0 \leq |k| \leq q
\right\}
\]

\noindent denote the Sobolev space with norm
\[
\|u \|_{q,p}
\defeq
\begin{cases}
\sum_{0 \leq |k| \leq q} \left(\|\partial^k u\|^p_p\right)^{1/p}  &\text{if } 1 \leq p < \infty,\\
\max_{0 \leq |k| \leq q} \|\partial^k u \|_{\infty} &\text{if } p = \infty,
\end{cases}
\]

\noindent where $\|\,\cdot\,\|_{\infty}$ denotes the standard supremum norm. Using these definitions, we denote by $C^q_0(\Omega)^3 \subset C^q(\Omega)^3$ with $q \in \mathbb{N}$ the completion of the space of vector fields of class $C^q$ who along with their derivatives of order less than or equal to $q$ converge to zero at infinity. The space $C^q_0(\Omega)^3$ is a Banach space for the norm $\|u \|_{q,\infty}$. Similarly, we define the Sobolev space $W^{q,p}_0(\Omega)^3$ as a space that consists of elements with compact support on $\Omega \subseteq \mathbb{R}^3$.

We introduce a pseudo-time variable $t \in [0,1]$, a suitable Hilbert space $\mathcal{H}$ of smooth vector fields in $\mathbb{R}^3$ (e.g., a RKHS; see \secref{s:rkhs}), and parameterize diffeomorphisms using smooth vector fields $v \in \mathcal{V}$, $\mathcal{V} \defeq L^2([0,1],\mathcal{H})$, $t \mapsto v_t \defeq v(t,\cdot\,)$, $v_t \in \mathcal{H}$. With this, we model the flow of $\mathbb{R}^3$-diffeomorphisms $f_t \defeq f(t, \cdot\,)$ as the solution of the ordinary differential equation (the flow equation)
\begin{equation}\label{e:flow}
\begin{aligned}
\partial_t f_t &= v_t(f_t) && \text{in}\; (0,1],\\
f_0 & = \id,
\end{aligned}
\end{equation}

\noindent where $\id : \mathbb{R}^3 \to \mathbb{R}^3$, $\id(x) = x$, is the identity transformation in $\mathbb{R}^3$ and the vector field $v_t$ tends to zero as $x \to \infty$. More precisely, we assume that $v_t$ is an element of $C^q_0(\Omega)^3 \supseteq \mathcal{H}$, as defined above, for any $t \in [0,1]$. This assumption, along with suitable regularity requirements in time guarantees that~\eqref{e:flow} admits a unique solution, and that the flow $f_t$ is a diffeomorphism of $\Omega$~\cite{Younes:2019a}. We refer to~\cite{Dupuis:1998a,Glaunes:2008a,Azencott:2010a,Younes:2019a} for a more rigorous discussion. We have already specified above that $v$ is $L^2$-integrable in time. However, we note that $L^1$-integrability is another option; see, e.g., \cite{Younes:2019a}. The smoothness of $v_t \in \mathcal{H}$ is typically determined using a Sobolev norm $\|\, \cdot\, \|_{q,p}$ as defined above, where $q$ defines the differentiability class and $p$ the integrability order. This choice is motivated by the fact that the above result is true for vector fields $v$ belonging to Banach or Hilbert spaces that are continuously embedded in $C^1_0(\Omega)^3$. A common and natural choice for $p$ is again $p=2$. The choice of $q$ depends on the dimension $d$ of the ambient space $\Omega \subseteq \mathbb{R}^d$, where $d\in\{2,3\}$ in most applications. Based on the Sobolev embedding theorem (see, e.g., \cite{Ziemer:2012a}), we can conclude that for $p=2$ and $q > \nicefrac{d}{2} + 1$ the embedding $W^{q,2}_0(\Omega) \hookrightarrow C^1(\bar{\Omega})$ is compact. Since this embedding holds for all components of $v_t$, we have that $v_t \in \mathcal{H} = W^{q,2}(\Omega)^3$ with $q > \nicefrac{5}{2}$ for $d=3$ is an admissible space that yields a diffeomorphic flow $f_t$ of smoothness class $1 \leq r < q - \nicefrac{3}{2}$, $r \in \mathbb{N}$.

The set of all admissible flows $y \defeq f_1$ at time $t=1$ is a subgroup \[\mathcal{G} \defeq \{ f_1 : \int_0^1 \|v_t\|_{\mathcal{H}} \,\d t < \infty\}\] of $C^r$-diffeomorphisms in $\mathbb{R}^3$. This subgroup can be equipped with a right-invariant metric defined as the minimal path length of all geodesics joining two elements $\phi,\psi \in \mathcal{G}$. We refer, e.g., to \cite{Trouve:1995a,Miller:2002a,Younes:2020a,Bauer:2014a} for a precise definition and mathematically rigorous discussion. The geodesic distance between $\id$ and a mapping $\phi$ corresponds to the square root of the kinetic energy
\begin{equation}\label{e:kinv}
\kin(v)
\defeq
\| v \|_{L^2([0,1],\mathcal{H})}^2
=
\int_0^1 \| v_t \|_{\mathcal{H}}^2 \,\d t
\end{equation}

\noindent subject to the constraint that $\phi$ is equal to the solution of~\eqref{e:flow} at time $t=1$ for the energy minimizing velocity $v$. We denote this geodesic distance by $\rho_{\mathcal{G}}(\id,\phi)$. We have,
\[
\rho_{\mathcal{G}}(\id,\phi)^2 \defeq \inf_{v\in \mathcal{V}}
\left\{
\int_0^1 \| v_t \|^2_{\mathcal{H}} \,\d t
: v \in \mathcal{V}, \,\, \phi = f_1,\,\,
\partial_t f_t = v_t(f_t),\,\,
f_0 = \id
\right\}.
\]

\noindent The geodesic distance between two diffeomorphic maps $\phi$ and $\psi$ is given by $\rho_{\mathcal{G}}(\id,\phi\circ{\psi}^{-1})^2$. Similarly, we can measure the distance between two shapes $s_0,s_1 \in \mathcal{S}$ in terms of the kinetic energy associated with the energy minimizing $v_t$ that gives rise to the diffeomorphic flow $f_t$ that maps the shape $s_0$ to the shape $s_1$. That is, we measure the distance between two shapes $s_0$ and $s_1$ as the smallest-cost deformation that maps $s_0$ to $s_1$. Using the geodesic distance $\rho_{\mathcal{G}}$ introduced above we obtain
\begin{equation}\label{e:shapemetric}
\rho_{\mathcal{S}}(s_0,s_1) \defeq \inf_{\phi \in \mathcal{G}}
\left\{\,
\rho_{\mathcal{G}}(\id,\phi)
: s_1 = \phi \app s_0\right\}.
\end{equation}

\noindent Consequently, to compare two shapes $s_0$ and $s_1$, we have to find the energy minimizing velocity $v$ for $\kin$ in \eqref{e:kinv} that parameterizes the map that maps $s_0$ to $s_1$. This is exactly the inverse problem we seek to solve in the present work. More precisely, we seek $y\in \mathcal{G}$ such that $s_1 = y \app s_0$, where $y \defeq f_1$ is the end point at $t=1$ of the geodesic $f_t$, $t \in [0,1]$ satisfying~\eqref{e:flow}.

\subsection{Variational Problem Formulation}\label{s:varopt}

We consider diffeomorphic shape matching as a variational, optimal control problem. We relax the hard constraint $y \app s_0 = s_1$ by introducing a distance functional $\dist : \mathcal{S} \times \mathcal{S} \to \mathbb{R}$ that quantifies the proximity between $y \app s_0$ and $s_1$. We model the diffeomorphism $y$ to be generated by some time dependent flow of smooth $\mathbb{R}^3$-vector fields $v_t$, as prescribed by the dynamical system in~\eqref{e:flow}. We arrive at the nonlinear optimal control problem~\cite{Beg:2005a,Dupuis:1998a,Azencott:2010a}
\begin{equation}\label{e:varoptprob}
\begin{aligned}
\minopt_{(f_t) \in \mathcal{F}\!,\;v\in\mathcal{V}} \;\; & \text{dist}(f_1\app s_0,s_1) + \reg(v)\\
\begin{aligned}
\text{subject to} \;\;\\{}\\
\end{aligned}
&
\begin{aligned}
\partial_t f_t & = v_t(f_t) && \text{in } (0,1],\\
f_0 & = \id.
\end{aligned}
\end{aligned}
\end{equation}

We specify the distance measure in~\eqref{e:varoptprob} in \secref{s:discdist}. The regularization functional $\reg : \mathcal{V} \to \mathbb{R}$ corresponds to the kinetic energy in \eqref{e:kinv}; more details can be found in \secref{s:disckinv}.

\begin{remark}
We note that the formulation above has been extended to the matching of finite sequences of shapes~\cite{Azencott:2008a,Zhang:2021a,Niethammer:2009a}. For simplicity, we omit this case here.
\end{remark}

\subsection{Numerical Discretization}\label{s:discretization}

We consider a \emph{discretize-then-optimize} approach. Other works using a discretize-then-optimize approach for LDDMM are described in~\cite{Polzin:2020a,Polzin:2018a,Mang:2017a}. For a discussion of tradeoffs between \emph{discretize-then-optimize} and \emph{optimize-then-discretize} approaches we refer to~\cite{Gunzburger:2003a}.

In our work, we model the surfaces $s_0$ (template shape) and $s_1$ (target shape) by a mesh of points $x_0 = (x_0^i)_{i=1}^{m_0} \in \mathbb{R}^{3m_0}$ and $x_1 = (x_1^i)_{i=1}^{m_1} \in \mathbb{R}^{3m_1}$ in $\mathbb{R}^3$, respectively. We subdivide the time horizon $[0,1]$ into $n$ equispaced cells of size $h = \nicefrac{1}{(n+1)}$ on a nodal grid. The associated discrete time points are denoted by $t^j \defeq (j-1)h$, $j = 1, \ldots, n+1$. We denote the state and control vectors of our problem by $x(t^j) = ( x^i(t^j) )_{i=1}^{m_0} \in \mathbb{R}^{3m_0}$ and $a(t^j) = (a^i(t^j) )_{i=1}^{m_0} \in \mathbb{R}^{3m_0}$ with $j = 1,\ldots,n+1$, respectively. The control $a^i(t^j) \in \mathbb{R}^3$ are coefficients associated with the $i$th state vector $x^i(t^j)\in\mathbb{R}^3$ at time $t^j$; they parameterize the continuous velocity field $v$ (see \secref{s:rkhs} for details). We rearrange the vectors $x(t^j)$ and $a(t^j)$ for $j=1,\ldots,n+1$ in lexicographical orderings
\begin{equation}\label{e:lexord}
x \defeq (x(t^1), \ldots, x(t^{n+1})) \in \mathbb{R}^{3m_0(n+1)}
\quad\text{and}
\quad a \defeq (a(t^1), \ldots, a(t^{n+1})) \in \mathbb{R}^{3m_0(n+1)},
\end{equation}

\noindent respectively, to be able to present the discretized version of the optimal control problem in~\eqref{e:varoptprob} in a compact format.

\begin{remark}
We note that our implementation allows for matching shape representations discretized with a varying number $m_i \in \mathbb{N}$ of points $(x_i^j)_{j=1}^{m_i}$, $i \in \{0,1\}$, despite the fact that the numerical experiments included in this study are performed for meshes $x_i$ with an identical number of discretization points $m_i$.
\end{remark}

\subsubsection{Reproducing Kernel Hilbert Space Representation}\label{s:rkhs}

RKHS are a key element of diffeomorphic shape matching~\cite{Jain:2013a,Miller:2015a,Richardson:2016a,Hsieh:2021a,Azencott:2010a,Arguillere:2016a,Glaunes:2004a}. As we have stated above, we assume that $v_t$ are smooth vector fields on $\mathbb{R}^3$. To guarantee this we have introduced an adequate Hilbert space $\mathcal{H}$, as outlined in \secref{s:flows}. Assuming that the Hilbert space $\mathcal{H}$ of $\mathbb{R}^3$ vector fields is continuously embedded into $C^q_0(\Omega)^3$, $\Omega \subseteq \mathbb{R}^3$, for some $q \geq 1$, there exists a constant $c > 0$ such that
\[
\|v_t\|_{q,\infty} \leq c \|v_t\|_{\mathcal{H}}
\]

\noindent for some $v_t \in \mathcal{H}$, where $\|\,\cdot\,\|_{\mathcal{H}}$ denotes the Hilbert norm on $\mathcal{H}$ with associated inner product $\langle\,\cdot\,, \,\cdot\,\rangle_{\mathcal{H}}$. From this it follows that $\mathcal{H}$ is an RKHS. Since $\mathcal{H}$ is a Cartesian product space for vector fields in $\mathbb{R}^3$, the construction of the associated RKHS is slightly more involved than in the scalar case. In particular, we have to assume that the kernel $\ker : \mathbb{R}^3 \times \mathbb{R}^3 \to \mathbb{R}^{3,3}$ is matrix valued. Based on the Riesz representation theorem we know that there exists a matrix valued kernel function such that for all $z,a\in\mathbb{R}^3$ the vector field $\ker(\,\cdot,z)a$ belongs to $\mathcal{H}$ and $\langle\ker(\,\cdot,z)a,v_t \rangle_{\mathcal{H}} = a^\mathsf{T} v_t(z)$ for $v_t \in \mathcal{H}$. Using this construction, it is possible to reduce the minimization in~\eqref{e:varoptprob} from an infinite-dimensional, nonparametric variational problem to a parametric one. This can be accomplished using an RKHS argument similar to the so called kernel trick in kernel methods~\cite{Scholkopf:2002a,Aronszajn:1950a,Meinguet:1979a}. Since the distance between the deformed shape $f_1 \app s_0$ and $s_1$ only depends on $f$ at time $t=1$ and we parameterize the template shape $s_0$ as a set of a few select points $\{x^i_0\}_{i=1}^{m_0}$ (see \secref{s:discretization}), it suffices to compute trajectories for these select points $x^i$, $i=1,\ldots,m_0$, that satisfy the flow equation~\eqref{e:flow}, i.e., $\partial_t x^i = v_t(x^i)$. It follows that for an optimal velocity $v$ we have that the squared norm $\|v_t\|_{\mathcal{H}}$ is minimal over all $\|u\|_{\mathcal{H}}^2$ with $u(x^i) = v_t(x^i)$ for all $t\in[0,1]$. Therefore, we have that $v_t$ is of the form~\cite{Azencott:2010a,Scholkopf:2002a}
\[
v_t(z) = \sum_{i=1}^{m_0} \ker(x^i(t),z) \, a^i(t) \quad \text{for all}\; z \in \mathbb{R}^3
\]

\noindent with $\ker :\mathbb{R}^3 \times \mathbb{R}^3 \to \mathbb{R}^{3,3}$ as defined above. With this, we reduced the search for $v \in L^2([0,1],\mathcal{H})$ to a set of time dependent, unknown coefficients $a(t) = (a^i(t))_{i=1}^{m_0}$ in $\mathbb{R}^3$.

Next, we specify more concretely the implementation details. We simplify the RKHS construction by modeling the kernel as a scalar function $\ker_{\sigma_v} : \mathbb{R}^3 \times \mathbb{R}^3 \to \mathbb{R}$ parameterized by $\sigma_v > 0$. The matrix representation introduced above is, equivalently, given by $\ker(\cdot, \cdot) = \ker_{\sigma_v}(\cdot, \cdot) I_3$, $I_3 \defeq \operatorname{diag}(1,1,1) \in \mathbb{R}^{3,3}$. With this, we obtain
\begin{equation}\label{e:velrkhs}
v_t(z) = \sum_{i=1}^{m_0} \ker_{\sigma_v}(x^i(t),z) \, a^i(t) \quad \text{for all}\; z \in \mathbb{R}^3,
\end{equation}

\noindent where the coefficients $a(t) = (a^i(t))_{i=1}^{m_0} \in \mathbb{R}^{3m_0}$, $m_0 \in \mathbb{N}$, are associated with the state vector $x(t) = (x^i(t))_{i=1}^{m_0} \in \mathbb{R}^{3m_0}$. Another key feature of this construction is that we can specify the kernel $\ker$ instead of defining the space $\mathcal{H}$~\cite{Younes:2019a,Aronszajn:1950a}.

\noindent The choice of the kernel $\ker_{\sigma_v}$ is crucial, since it uniquely determines the associated RKHS. We consider a Gaussian kernel defined as
\begin{equation}\label{e:gaussiankernel}
\ker_{\sigma_v}(u,u') = \frac{1}{\left(\sqrt{2\pi}\sigma_v\right)^3} \exp\left(-\frac{1}{2} \|u - u'\|^2/\sigma_v^2\right)
\end{equation}

\noindent for any $u,u'\in\mathbb{R}^3$. This kernel is symmetric, positive definite and yields a metric that is invariant under rotation and translation. (We note that other kernels have been considered~\cite{Jain:2013a}.) The Hilbert space $\mathcal{H}_{\sigma_v}^3$ of smooth vector fields in $\mathbb{R}^3$ is defined as
\begin{equation}\label{e:rkhs}
\mathcal{H}_{\sigma_v}^3
\textstyle
\defeq
\overline{
\left\{
w : w(\,\cdot\,) = \sum_{i=1}^{m_0} \ker_{\sigma_v}(x^i,\,\cdot\,)\, a^i,
\,\, a^i \in \mathbb{R}^3,
\,\, x^i \in \mathbb{R}^3,
\,\,m \in \mathbb{N}
\right\}}.
\end{equation}

\noindent We note that the Sobolev embedding hypothesis stated earlier in this work is satisfied if we model the space $\mathcal{H}$ using the self-reporducing Hilbert space $\mathcal{H}_{\sigma_v}^3$ for any $q > \nicefrac{5}{2}$.

\subsubsection{Discretized Kinetic Energy}\label{s:disckinv}

The norm on $\mathcal{H}_{\sigma_v}^3$ in~\eqref{e:rkhs} is defined by the square root of
\begin{equation}\label{e:rkhs-norm}
\| \sum_{l=1}^{m_0} \ker_{\sigma_v}(\,\cdot\,,x^l)\,a^l\|_{\mathcal{H}_{\sigma_v}^3}^2
= \sum_{k=1}^{m_0}\sum_{l=1}^{m_0} \ker_{\sigma_v}(x^k,x^l)\,
\langle a^k, a^l \rangle_{\mathbb{R}^3},
\end{equation}

\noindent where $\langle \,\cdot\,, \,\cdot\, \rangle_{\mathbb{R}^3}$ denotes the standard inner product in $\mathbb{R}^3$, and $a^l \in \mathbb{R}^3$ are coefficient vectors associated with the mesh points $x^l \in \mathbb{R}^3$. We use a left Riemann summation~\cite{Quarteroni:2010a} to discretize the integration in time that appears in \eqref{e:kinv}. With this, the discrete version of the kinetic energy in~\eqref{e:kinv} is given by
\begin{equation}\label{e:kinvdisc}
\kin^h (a) \defeq h \sum_{j=1}^{n} \text{hilb}(t^{j}),
\end{equation}

\noindent where $h = \nicefrac{1}{(n+1)}$ denotes the size of the cells of the temporal mesh $[t^1, \ldots, t^{n+1}] \in \mathbb{R}^{n+1}$ and
\[
\text{hilb}(t^j) \defeq \sum_{k=1}^{m_0}\sum_{l=1}^{m_0} \ker_{\sigma_v}(x^k(t^j),x^l(t^j))\, \langle a^i(t^j), a^l(t^j)\rangle_{\mathbb{R}^3}.
\]

\begin{figure}
\centering
\includegraphics[width=0.98\textwidth]{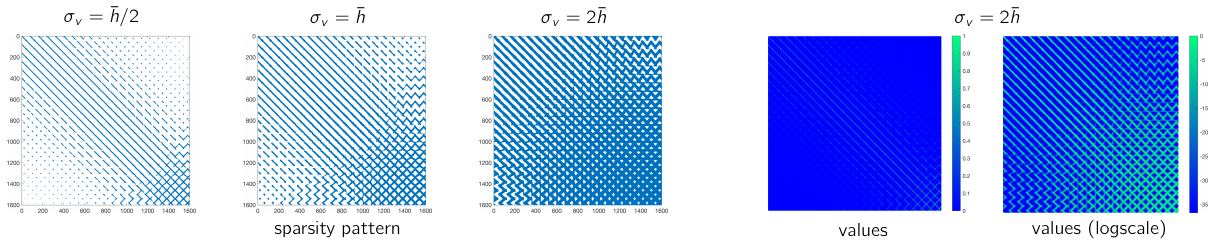}
\caption{Left block: Spy plot for the $m_0 \times m_0$ kernel matrix $K$ in~\eqref{e:kernelmat-k} for different choices of the bandwidth $\sigma_v$ as a function of the mean edge length $\bar{h}$. The plots have been generated for an exemplary mitral valve shape. The matrices are of size $1600 \times 1600$. The mean edge length $\bar{h}$ is $\fnum{1.2127}$. We set $\sigma_v = s \bar{h}$ with $s \in\{1/4,1/2,1,2,4\}$. We show the sparsity pattern for $s = 1/2$ (left), $s=1$ (middle) and $s=2$ (right). The number of nonzero entries (more precisely, entries with a value below $\num{1e-16}$) is $42\,966$  (1.68\%; $s=1/4$; not shown), $157\,988$ (6.17\%; $s=1/2$), and $523\,740$ (percentage: 20.46\%; $s=1$), $1\,470\,372$  (57.44\%; $s=2$), and $2\,559\,418$ (99.97\%; $s=4$; not shown), respectively. Right block: Visualization of the entries of the kernel matrix $K$ for $\sigma_v = 2 \bar{h}$ (left: values for the entries $(K_{ij})_{ij}$; right: logarithmic scale)\label{f:spy-plot-kernel-mat}}
\end{figure}

\noindent Under the assumption that the controls $a$ are ordered according to~\eqref{e:lexord}, we can rewrite~\eqref{e:kinvdisc} as
\begin{equation}\label{e:kinvdisclex}
\kin^h (a) \defeq h a^\mathsf{T} B a,
\end{equation}

\noindent where $B$ is a block-diagonal matrix of the form
\begin{equation}
\label{e:kernelmat}
B \defeq B[x]
= \operatorname{diag}\left(B^1[x],\ldots,B^{n+1}[x]\right)
\in \mathbb{R}^{3m_0(n+1),3m_0(n+1)},\quad B^j[x] = I_3 \otimes K[x(t^j)] \in \mathbb{R}^{3m_0,3m_0},
\end{equation}

\noindent where $\otimes$ denotes the Kronecker product, i.e.,
\begin{equation}\label{e:defkronprod}
A \otimes B
=
\begin{pmatrix}
a_{11} B & \hdots & a_{1n}B \\
\vdots   & \ddots & \vdots  \\
a_{m1} B & \hdots &  a_{mn}B
\end{pmatrix}
\in \mathbb{R}^{mp,nq}
\end{equation}

\noindent for any $A\in\mathbb{R}^{n,m}$, $B \in \mathbb{R}^{p,q}$, $I_3 = \operatorname{diag}(1,1,1)\in\mathbb{R}^{3,3}$, and $K[x(t^j)]$ are Gram matrices of the form
\begin{equation}\label{e:kernelmat-k}
K[x(t^j)]
= \begin{pmatrix}K_{lk}(x(t^j))\end{pmatrix}_{l,k=1}^{m_0,m_0} \in \mathbb{R}^{m_0,m_0},
\quad K_{lk}(x(t^j)) \defeq \ker_{\sigma_v}(x^l(t^j) , x^k(t^j)) \in \mathbb{R},
\end{equation}

\noindent with $x^l(t^j) \in \mathbb{R}^3$ and $x^k(t^j) \in \mathbb{R}^3$ denoting the $l$th and $k$th entry of the state vector $x$ at time point $t^j$ and $\ker_{\sigma_v} : \mathbb{R}^3 \times \mathbb{R}^3 \to \mathbb{R}$ as in~\eqref{e:gaussiankernel}. We illustrate the matrix $B[x]$ for a particular shape $x$ for different choices of the bandwidth $\sigma_v$ in \figref{f:spy-plot-kernel-mat}.

\begin{remark}
We note that the kernel matrices $K$ in~\eqref{e:kernelmat-k} that appear in the regularization operator $B$ in \eqref{e:kernelmat} are functions of the state vector $x(t^j)$. This significantly complicates the optimization problem since we need to differentiate not only with respect to the control variable $a$ but also the state variable $x$ (and the dual variables) to derive the optimality conditions. To avoid this complication, we fix the kernel matrices that enter the regularizer in \eqref{e:kinvdisclex} by replacing $K[x(t^j)]$ by $K[x(t^0)] = K[x_0]$ for all $j = 1, \ldots, n$. This approximation significantly simplifies the optimization problem: The regularization operator $B$ does not change during optimization since it depends on the fixed (static) shape representation $x_0$ and not on the state variable $x(t^j)$ that is updated during the course of the optimization. We end up with a quadratic term that does not need to be differentiated with respect to the state vector $x$. We have seen that this approximation works well in practice, even if the shapes to be registered, i.e., the shapes parameterized by $x_0$ and $x_1$, differ significantly. (If the shapes are different, we expect the state vector $x$ to change significantly. Consequently, the kernel entries in $B$ are expected to change as we move from $x_0$ to $x_1$.) We note that the norm in~\eqref{e:kinvdisclex} still depends on the dynamics of our problem, since the coefficients $a$ in~\eqref{e:kinvdisclex} change in time.
\end{remark}

\subsubsection{Discrete Distance Measures between Shapes}\label{s:discdist}

Before we discuss our framework, we note that several geometric distance measures have been proposed, based on different shape representations. They include representations as measures~\cite{Azencott:2010a,Glaunes:2004a} (see also below), currents~\cite{Vaillant:2005a}, or varifolds~\cite{Kaltenmark:2017a,Antonsanti:2021a,Hsieh:2022b,Charon:2013a}. We refer to these papers for more details.

Let $\mathcal{S}$ be the set of all compact, smooth 3D-surfaces properly embedded in $\mathbb{R}^3$ with piecewise smooth boundaries. Let $\mathcal{B}$ denote the space of bounded Borel measures on $\mathbb{R}^3$ endowed with the Hilbert norm $\|\,\cdot \,\|_{\mathcal{H}_{\sigma}}$ associated with the scalar inner product
\begin{equation}\label{e:scalarprodhilb}
\langle \mu,\tilde{\mu} \rangle_{\mathcal{H}_{\sigma}}
\defeq \int_{\mathbb{R}^3} \int_{\mathbb{R}^3} \ker_{\sigma}(x,\tilde{x})\,\d \mu(x) \,\d \tilde{\mu}(\tilde{x})
\end{equation}

\noindent for any $\mu,\tilde{\mu} \in \mathcal{B}$, where $\ker_{\sigma}: \mathbb{R}^3 \times \mathbb{R}^3 \to \mathbb{R}_{>0}$ denotes any smooth, symmetric, translation invariant and bounded positive definite kernel parameterized by $\sigma > 0$. (We are going to consider Gaussian kernels as defined in~\eqref{e:gaussiankernel}; hence the dependence on the parameter $\sigma$.) We can identify each submanifold $s$ regularly embedded in $\mathbb{R}^3$ with Borel measure $\mu \in \mathcal{B}$ induced on $s$ by the Lebesgue measure of $\mathbb{R}^3$. The corresponding squared distance between two bounded Borel subsets $s$ and $\tilde{s}$ with measures $\mu$ and $\tilde{\mu}$, respectively, is then defined by
\begin{equation}\label{e:boreldist-generic}
\dist_{\mathcal{H}_{\sigma}}(s,\tilde{s}) = \|\mu - \tilde{\mu}\|^2_{\mathcal{H}_{\sigma}}.
\end{equation}

As we have outlined in \secref{s:discretization}, we identify these surfaces as a mesh of points $x_i \in \mathbb{R}^{3m_i}$, $i=0,1$, respectively. In this setting, it is natural to represent the shapes $s_i \in \mathcal{S}$ as bounded positive Radon measures on $\mathbb{R}^3$, namely as weighted sums of Dirac measures $\delta(x^k_i)$, $k=1,\ldots,m_i$, $i=0,1$, associated with the point sets $x_i \in \mathbb{R}^m$~\cite{Glaunes:2004a,Azencott:2010a,Azencott:2008a}. In particular,
\begin{equation}\label{e:weighteddirac}
\mu_0 = \sum_{k=1}^{m_0} p_k\, \delta(x_0^k),
\quad
\mu_1 = \sum_{k'=1}^{m_1} q_{k'}\, \delta(x_1^{k'}),
\end{equation}

\noindent with weights $p_k, q_{k'} \in \mathbb{R}$. With this, we arrive at the discrete representation of the distance functional in~\eqref{e:varoptprob} given by
\begin{equation}\label{e:boreldist}
\text{dist}^h(x_0,x_1)
= \frac{\alpha}{2} \|\mu_0 - \mu_1\|_{\mathcal{H}_{\sigma_s}}^2.
\end{equation}

\noindent As stated above, the set of bounded Borel measures is endowed with the Hilbert scalar product in \eqref{e:scalarprodhilb}. The associated squared Hilbert norm $\|\mu_0 - \mu_1\|_{\mathcal{H}_{\sigma_s}}^2$ in~\eqref{e:boreldist} represents the Borel distance between the shapes $s_0$ and $s_1$, associated with the Gaussian kernels $\ker_{\sigma_s}$ defined in~\eqref{e:gaussiankernel} with an adequate choice for the parameter $\sigma_s > 0$. With the definition of the inner product in~\eqref{e:scalarprodhilb} and the distance in \eqref{e:boreldist-generic} we obtain
\[
\|\mu_0 - \mu_1\|_{\mathcal{H}_{\sigma_s}}^2
 =
\langle\mu_0 - \mu_1,\mu_0 - \mu_1\rangle_{\mathcal{H}_{\sigma_s}}
=
\langle\mu_0, \mu_0\rangle_{\mathcal{H}_{\sigma_s}}
-2 \langle\mu_0, \mu_1\rangle_{\mathcal{H}_{\sigma_s}}
+ \langle\mu_1, \mu_1\rangle_{\mathcal{H}_{\sigma_s}},
\]

\noindent where the inner products $\langle\,\cdot\,, \,\cdot\,\rangle_{\mathcal{H}_{\sigma_s}}$ are, based on the representation of $\mu_i$, $i=0,1$, in~\eqref{e:weighteddirac}, given by
\begin{subequations}
\label{e:scalarprodhilb:shapes}
\begin{align}
\langle\mu_0, \mu_0\rangle_{\mathcal{H}_{\sigma_s}}
&= \int_{\mathbb{R}^3}\int_{\mathbb{R}^3} \ker_{\sigma_s}(x,x') \,\d \mu_0(x) \,\d \mu_1(x')
= \sum_{k = 1}^{m_0} \sum_{k' = 1}^{m_0} p_k p_{k'} \ker_{\sigma_s}(x_0^k, x_0^{k'}),
\\
\langle\mu_0, \mu_1\rangle_{\mathcal{H}_{\sigma_s}}
&= \int_{\mathbb{R}^3}\int_{\mathbb{R}^3} \ker_{\sigma_s}(x,x') \,\d \mu_0(x) \,\d \mu_1(x')
= \sum_{k = 1}^{m_0} \sum_{k' = 1}^{m_1} p_k q_{k'} \ker_{\sigma_s}(x_0^k, x_1^{k'}),
\\
\langle\mu_1, \mu_1\rangle_{\mathcal{H}_{\sigma_s}}
&= \int_{\mathbb{R}^3}\int_{\mathbb{R}^3} \ker_{\sigma_s}(x,x') \,\d \mu_1(x) \,\d \mu_1(x')
= \sum_{k = 1}^{m_1} \sum_{k' = 1}^{m_1} q_k q_{k'} \ker_{\sigma_s}(x_1^k, x_1^{k'}),
\end{align}
\end{subequations}

\noindent respectively.

The parameter $\alpha > 0$ in~\eqref{e:boreldist} is a weight that controls the contribution of the distance in~\eqref{e:boreldist} versus the contribution of the kinetic energy in~\eqref{e:kinvdisc}. We discuss policies for selecting $\alpha$ in \secref{s:parameterchoices}.

\begin{remark}
We note that $\sigma_v \not= \sigma_s$, in general. That is, $\sigma_v$ is chosen to control the smoothness of the velocity $v$ whereas $\sigma_s$ controls the contribution of all points $x_j^l \in \mathbb{R}^3$ to the distance for a particular point $x_i^k \in \mathbb{R}^{3}$ with $i,j\in\{0,1\}$ for all $k,l=1,\ldots,m$. We specify heuristic choices for the parameters $\sigma_s$ and $\sigma_v$ in \secref{s:parameterchoices}.
\end{remark}

Likewise to the kinetic energy in \eqref{e:kinvdisclex}, we consider a matrix-vector representation to compactly represent the distance in~\eqref{e:boreldist}. We have
\begin{equation}\label{e:discretedist}
\dist^h(x_0,x_1) = \frac{\alpha}{2}\left(\phi(x_0, x_0) -2 \phi(x_0, x_1) + \phi(x_1, x_1)\right),
\end{equation}

\noindent with
\[
\phi(x_i,x_j)
=
\frac{1}{m_im_j} \sum_{k = 1}^{\,\,\,m_i\phantom{j}} \sum_{l = 1}^{m_j}\ker_{\sigma_s}(x_i^k, x_j^l)
= \frac{1}{m_im_j} e_{m_im_j}^\mathsf{T} \exp\left(-\frac{1}{2\sigma_s}S(d(x_i,x_j) \odot d(x_i,x_j))\right),
\]

\noindent where $S \defeq e_3^\mathsf{T} \otimes I_{m_im_j}$, $S \in \mathbb{R}^{m_im_j,3m_im_j}$, $e_3 = (1,1,1) \in \mathbb{R}^3$, $I_{m_im_j} \defeq \operatorname{diag}(1,\ldots,1) \in \mathbb{R}^{m_im_j,m_im_j}$, $\exp : \mathbb{R}^{m_im_j} \to \mathbb{R}^{m_im_j}$ is the pointwise exponential, $\otimes$ is the Kronecker product (see \eqref{e:defkronprod}), $\odot$ is the Hadamard product defined as the elementwise/entrywise product of two matrices/vectors, i.e.,
\begin{equation}\label{e:defhadaprod}
(A\odot B)_{kl} = (A)_{kl}(B)_{kl}, \quad k, l = 1,\ldots,n,
\end{equation}

\noindent for any $A,B \in \mathbb{R}^{n,n}$, and all-to-all distance vector $d(x_i,x_j) \defeq P_1x_i - P_2x_j \in \mathbb{R}^{3m_im_j}$ with permutation matrices $P_1 \defeq I_3 \otimes I_{m_i} \otimes e_{m_j}$, $P_1\in \mathbb{R}^{3m_im_j,3m_i}$, $P_2 \defeq I_3 \otimes e_{m_i} \otimes I_{m_j}$, $P_2 \in \mathbb{R}^{3m_im_j,3m_j}$, respectively. The all-to-all distance vector $d$ introduced above is a vector whose entries represent the difference between any point $\{x_i^l\}_{l = 1}^{m_i}$ of $x_i$ to any point $\{x_j^k\}_{k = 1}^{m_j}$ of $x_j$ with $i,j \in \{0,1\}$.

Comparing the expression for $\phi$ in~\eqref{e:discretedist} to~\eqref{e:scalarprodhilb:shapes} reveals that we select uniform weights $p_k = 1/m_0$ and $q_{k'} = 1/m_1$ for all $k = 1,\ldots,m_0$, $k'=1,\ldots,m_1$, respectively, for the weighted sums in~\eqref{e:weighteddirac}. We use uniform weights for simplicity. Based on our past experience and the results reported in this work, this choice yields good performance in terms of convergence and shape matching accuracy. However, we note that it is certainly possible to use weights that are adaptive and depend on the parameterization of the shape. It remains subject to future work to explore if such an adaptive scheme is beneficial for the problems considered in our work. A simple modification we have implemented, that does not quite depend on the parameterization, is to increase the contribution of the points located on the boundary of the curved surface to the matching term (the boundaries of the mitral valve leaflets in our case; see \figref{f:mitral-valve-regions} for an illustration). That is, we increase $\alpha$ in~\eqref{e:discretedist} for points located on the boundary. We observed in our past work~\cite{Zhang:2021a} that this leads to a more accurate matching close to the boundary. We omit these details for simplicity of presentation.

\begin{remark}
The permutation matrices $P_l$, $l=1,2$, are introduced to provide a compact representation of the computation of the all-to-all difference between the points of $x_i$, $x_j$, $i,j=0,1$, across all spatial dimensions, respectively. We note that the Kronecker product notation is only used for presentation purposes. In computation, tensor operations are used to avoid explicit formation of these large, sparse matrices.
\end{remark}

\subsubsection{Numerical Time Integration}\label{s:timeintegration}

We use an explicit Euler method to integrate~\eqref{e:flow} forward in time~\cite{Azencott:2008a,Azencott:2010a}. Using the RKHS representation in~\eqref{e:velrkhs} to model $v_t$ we have
\begin{equation}\label{e:expeuler}
x^l(t^{j+1})
= x^l(t^j) + h \sum_{k=1}^{m_0} \ker_{\sigma_v}(x^k(t^j),x^l(t^j)) \, a^k(t^j),
\quad l = 1,\ldots,m_0,
\end{equation}

\noindent for the $l$th entry $x^l(t^{j+1})\in\mathbb{R}^3$ of the state vector $x = (x(t^1), \ldots, x(t^{n+1}))\in \mathbb{R}^{3m_0(n+1)}$ at time $t^{j+1}$, with uniform time step size $h = 1/(n+1)$. Using the lexicographical ordering in~\eqref{e:lexord}, we can represent~\eqref{e:expeuler} as
\[
\begin{pmatrix}
G^x\;\;G^a
\end{pmatrix}
\begin{pmatrix}
x\\a
\end{pmatrix}
= q
\]

\noindent with state-control vector $(x,a) \in \mathbb{R}^{6m_0(n+1)}$. The matrices $G^x$ and $G^a \defeq G^a[x]$ are given by the lower block-bidiagonal and lower block-diagonal matrices
\begin{subequations}
\label{e:expeulermat}
\begin{equation}
G^a[x] =
-h
\begin{pmatrix}
0      & \cdots & 0          \\
B^1[x] & \ddots & \vdots     \\
\vdots & \ddots & 0          \\
0      & \cdots & B^{n+1}[x]
\end{pmatrix}
\in \mathbb{R}^{3m_0(n+1),3m_0(n+1)}
\end{equation}

\noindent and
\begin{equation}
G^x
=
\begin{pmatrix}
 I_{3m_0} & \cdots & 0       \\
-I_{3m_0} & \ddots & \vdots  \\
\vdots  & \ddots & I_{3m_0}  \\
 0      & \cdots & -I_{3m_0}
\end{pmatrix}
\in \mathbb{R}^{3m_0(n+1),3m_0(n+1)},
\end{equation}
\end{subequations}

\noindent $I_{3m_0} = \operatorname{diag}(1,\ldots,1)\in\mathbb{R}^{3m_0,3m_0}$, respectively. The right hand side $q\in\mathbb{R}^{3m_0(n+1)}$ is given by
\begin{equation}\label{e:expeulerrhs}
q = \begin{pmatrix}x_0\\0\end{pmatrix}
\end{equation}

\noindent and $B^j[x]$ is as defined in~\eqref{e:kernelmat}.

\begin{remark}
Notice that the matrix $G^x$ is a sparse matrix whereas the lower block-diagonal entries of $G^a$ are, in general, dense (as a function of $\sigma_v$; see \figref{f:spy-plot-kernel-mat}). In our implementation, the matrices $G^x$ and $G^a$ (and the associated adjoint operators) are not assembled and/or stored; we implement matrix-vector products using the block matrices, instead. However, we assemble and store the matrices $K[x(t^j)]$ which imposes memory restrictions on our current implementation.
\end{remark}

\begin{remark}
The explicit Euler scheme considered in this work is only first order accurate in time. Our methodology is modular in the sense that we can replace this scheme with any other (high-order) numerical time integration scheme (at the expense of computational complexity and runtime). In past work on a different problem formulation we observed that using higher-order time integration schemes for computing the  diffeomorphic flow can be beneficial~\cite{Mang:2017a}. In particular, we observed that a fourth order Runge–Kutta method for numerical time integration is more expressive when it comes to modelling complex diffeomorphic flows~\cite{Mang:2017a}. In other past work, we consider second order accurate time integration schemes~\cite{Mang:2017b,Mang:2015a,Mang:2016b,Mang:2019a}. We did not observe issues in the present application with our first order scheme. Exploring high-order schemes for numerical time integration remains subject to future work.
\end{remark}

\subsubsection{Discretized Optimal Control Problem}

If we put all building blocks described above together, we arrive the following discrete version of~\eqref{e:varoptprob}:
\begin{equation}\label{e:optprobdisc}
\begin{aligned}
\minopt_{a,\, x} \;\; & \dist^h(Qx,x_1) + \kin^h(a)\\
\text{subject to} \;\; & \begin{pmatrix}G^x\;\; G^a\end{pmatrix}
\begin{pmatrix}
x\\a
\end{pmatrix}
= q.
\end{aligned}
\end{equation}

\noindent Here, $(x,a) \in \mathbb{R}^{6m_0(n+1)}$ denotes the entire state-control trajectory, with $x$ and $a$ as in~\eqref{e:lexord}, discretized kinetic energy $\kin^h$ as in~\eqref{e:kinvdisc}, discretized distance $\dist^h$ as in~\eqref{e:discretedist}, matrix-operators $G^a$ and $G^x$ as in~\eqref{e:expeulermat}, and right-hand side $q$ as in \eqref{e:expeulerrhs}, respectively. Notice that we apply an observation operator $Q$ to the state variable $x$, such that we compute the proximity between the solution of the state equation (equality constraint in~\eqref{e:optprobdisc}), i.e., the deformed shape at time $t^{n+1}$ for a given control $a$, to the reference shape $x_1$. More precisely, $x(t^{n+1}) = Qx$, $Q = (0_{3m_0}\,\ldots\,0_{3m_0}\,I_{3m_0}) \in \mathbb{R}^{3m_0,3m_0(n+1)}$ with zero matrices $0_{3m_0} = \operatorname{diag}(0,\ldots,0) \in \mathbb{R}^{3m_0,3m_0}$ and identity matrix $I_{3m_0} \in \mathbb{R}^{3m_0,3m_0}$, respectively.

\subsection{Numerical Optimization}\label{s:optimization}

In what follows, we describe different methods to solve~\eqref{e:optprobdisc} numerically.

\subsubsection{Adjoint-Based Line Search Methods}

In our past work, our group has developed several adjoint based first and second order numerical optimization algorithms~\cite{Azencott:2010a,Freeman:2014a,Jajoo:2011a,Qin:2013a}. These approaches can be summarized by the iterative scheme~\cite{Nocedal:2006a,Boyd:2004a}
\begin{equation}\label{e:linesearchmethod}
a^{(k+1)} = a^{(k)} - t^{(k)} s^{(k)}, \quad k = 1,2,\ldots,
\end{equation}

\noindent where $k \in \mathbb{N}$ denotes the iteration index, $a^{(k)} \in \mathbb{R}^{3m_0}$ is the iterate, and $s^{(k)} \defeq - B^{(k)} g^{(k)}$ represents the search direction with $B^{(k)} \in \mathbb{R}^{3m_0,3m_0}$, $m_0 \in \mathbb{N}$, reduced gradient $g^{a,(k)}\in \mathbb{R}^{3m_0}$ (with respect to the control variable $a$), and line search parameter $t^{(k)} > 0$. We note that we consider a reduced space method in~\eqref{e:linesearchmethod}. That is, we only iterate on the reduced space of the control variable $a$. In \emph{full space} or \emph{all-at-once methods}, all unknowns of the problem (i.e., the state variable(s), the adjoint (or dual) variable(s), and the control variable(s)) are updated simultaneously. This affects the structure and size of the vectors and matrices involved in the computations associated with \eqref{e:linesearchmethod}, respectively. We refer, e.g., to~\cite{Biros:2005a,Biros:2005b,Biegler:2003a,Biegler:1995a,Herzog:2010a,Akccelik:2006a} for more details. The choice of $B^{(k)}$ determines the type of algorithm~\cite{Nocedal:2006a,Boyd:2004a}. If we set $B^{(k)}$ to $I_{3m_0} \defeq \operatorname{diag}(1,\ldots,1)\in\mathbb{R}^{3m_0,3m_0}$, we obtain a gradient descent algorithm. If we select $B^{(k)} = (H^{(k)})^{-1}$, where $H^{(k)} \in \mathbb{R}^{3m_0,3m_0}$ is the Hessian matrix, we arrive at a Newton algorithm. Other options for the matrix $B^{(k)}$ include approximations to the Hessian, resulting in different flavors of quasi-Newton methods. These methods vary in their rate of convergence, ranging from linear (gradient descent), to super-linear (quasi-Newton methods), to quadratic (Newton method)~\cite{Nocedal:2006a,Boyd:2004a}. In the context of LDDMM (and related formulations), quasi--Newton algorithms have been considered in~\cite{Ashburner:2011a,Mang:2015a,Mang:2016a,Mang:2019a} (specifically, adjoint-based Gauss--Newton) and \cite{Hsieh:2021a,Polzin:2016a} (BFGS; adjoint-based or automatic differentiation), respectively.

We omit additional details for the sake of brevity and refer to our past work~\cite{Azencott:2010a,Freeman:2014a,Jajoo:2011a,Qin:2013a}. We provide bits and pieces of the building blocks for these algorithms in the following sections, including derivative information of the variational optimization problem~\eqref{e:optprobdisc}. We note that we demonstrated in~\cite{Zhang:2021a} that our (original) implementation of the splitting approach discussed below outperforms the Newton algorithm developed in our prior work~\cite{Qin:2013a} in terms of runtime. Our improved implementation presented in this work is about $2\times$ faster and requires less memory, allowing for the solution of bigger problems than those considered in~\cite{Zhang:2021a} due to a matrix-free implementation.

\subsubsection{Operator Splitting}

In~\cite{Zhang:2021a}, we propose the use of an operator splitting strategy to solve the nonlinear control problem in~\eqref{e:varoptprob}. The considered approach is typically referred to as \emph{alternating direction method of multipliers} ({\bf ADMM})~\cite{Boyd:2011a,ODonoghue:2013a,Parikh:2013a}. It is equivalent to the \emph{Douglas--Rachford splitting method}~\cite{ODonghue:2013a,Bauschke:2011a,Glowinski:2016b,Glowinski:2016a,Glowinski:2016c,Bukac:2016a}, which was initially developed to numerically solve PDEs~\cite{Douglas:1956a}. In its modern form, this algorithm was introduced by Glowinski and Marrocco~\cite{Glowinski:1975a} and Gabay and Mercier~\cite{Gabay:1976a}. Recent work of Glowinski and colleagues in the area of operator splitting includes~\cite{Glowinski:2016a,Glowinski:2017a,Glowinski:1989a,Dean:2002a,Bukavc:2014a,Deng:2019a}. We refer to~\cite{Glowinski:2016a} for a general exposition of operator splitting techniques and to~\cite{Boyd:2011a,ODonghue:2013a,Parikh:2013a,Goldfarb:2012a} for a lucid discussion of operator splitting approaches in the context of optimization.

\paragraph{Problem Formulation}

To concisely present the splitting algorithm developed in \cite{Zhang:2021a}, we introduce the set
\begin{equation}\label{e:staconpair}
\mathcal{C} \defeq
\left\{
\begin{pmatrix}
x\\a
\end{pmatrix}
\in \mathbb{R}^{6m_0(n+1)}
:
\left[G^x\;\; G^a\right]
\begin{pmatrix}
x\\a
\end{pmatrix}
= q
\right\}
\end{equation}

\noindent of state-control pairs $(x,a) \in \mathbb{R}^{6m_0(n+1)}$ that satisfy the dynamical system (constraint) in~\eqref{e:optprobdisc}. Moreover, let $\operatorname{ind}_{\mathcal{C}} : \mathbb{R}^{3m_0(n+1)} \times \mathbb{R}^{3m_0(n+1)} \to \{0,\infty\}$ denote the indicator function for the set $\mathcal{C}$ in \eqref{e:staconpair}, i.e.,
\[
\text{ind}_{\mathcal{C}}(x,a) \defeq
\begin{cases}
0 & \text{if } (x,a) \in \mathcal{C}, \\
\infty & \text{otherwise}.
\end{cases}
\]

In addition, with slight abuse of notation, we define
\[
\dist^h(x,a) \defeq \dist^h(Qx,x_1)
\quad
\text{and}
\quad
\kin^h(x,a) \defeq \kin^h(a)
\]

\noindent with discretized distance $\dist^h$ and discretized kinetic energy $\kin^h$ as in \eqref{e:optprobdisc}, respectively. Using this notation allows us to better illustrate the derivation of the splitting strategy. We can rewrite~\eqref{e:optprobdisc} to obtain
\begin{equation}\label{e:optprobdisc-notation}
\begin{aligned}
\minopt_{a,\, x} \;\; & \dist^h(x,a) + \kin^h(x,a)\\
\text{subject to} \;\; & \begin{pmatrix}G^x\;\; G^a\end{pmatrix}
\begin{pmatrix}
x\\a
\end{pmatrix}
= q
\end{aligned}
\end{equation}

\noindent with variables $(x,a) \in \mathbb{R}^{6m_0(n+1)}$. With the indicator function $\operatorname{ind}_{\mathcal{C}}$ introduced above, we can represent the equality constrained optimization problem in~\eqref{e:optprobdisc} as an unconstrained optimization problem of the form
\begin{equation}
\label{e:optprobdiscunc}
\minopt_{x,\,a} \;\;
\text{ind}_\mathcal{C}(x,a) + \text{kin}^h(x,a) + \text{dist}^h(x,a).
\end{equation}

To derive our splitting approach, we consider the \emph{consensus form}~\cite{Boyd:2011a,Goldfarb:2012a} of~\eqref{e:optprobdiscunc}. The fundamental idea of the consensus form of operator splitting is to solve a coupled, multi-objective optimization problem (that is difficult to solve in its original form) by splitting it into multiple (easier) subproblems to be solved separately. These subproblems are then tied together by a constraint that enforces convergence to a single, consensus solution at the end of optimization. For the particular case of~\eqref{e:optprobdiscunc}, we split the problem into two separate parts by introducing a second pair of state-control variables $(\tilde{x},\tilde{a}) \in \mathbb{R}^{6m_0(n+1)}$. To ensure convergence to one solution, we introduce an equality constraint that enforces the state-control pairs $(x,a)$ and $(\tilde{x},\tilde{a})$ to be equal---the \emph{consensus} or \emph{consistency constraint}. We obtain
\begin{equation}
\label{e:consensusform}
\begin{aligned}
\minopt_{x,\,a} & \;\;
\left(
\text{ind}_{\mathcal{C}}(x,a) + \text{kin}^h(x,a)
\right)
+ \text{dist}^h(\tilde{x},\tilde{a})
\\
\text{subject to} & \;\;
\begin{pmatrix}
x \\ a
\end{pmatrix}
=
\begin{pmatrix}
\tilde{x} \\ \tilde{a}
\end{pmatrix}.
\end{aligned}
\end{equation}

As indicated in~\eqref{e:consensusform}, we have split the objective in \eqref{e:optprobdisc-notation} into two separate parts, with different variables. The first part contains the indicator function for the set $\mathcal{C}$ (i.e., the model for the discretized ODE constraint) and the kinetic energy. The second part encodes the mismatch term. The equality constraint states that these terms are ``in consensus.'' Using this formulation, we can derive an iterative scheme as follows. We use the method of Lagrange multipliers to solve \eqref{e:consensusform}. We introduce the dual variables $(\nu,\omega) \in \mathbb{R}^{6m_0(n+1)}$ for the equality constraint $(x, a) -  (\tilde{x}, \tilde{a}) = 0$ in \eqref{e:consensusform}. With this, the augmented Lagrangian for the problem in~\eqref{e:consensusform} is given by
\[
\mathcal{L}_{\rho}(\phi) =
\left(
\text{ind}_{\mathcal{C}}(x,a) + \text{kin}^h(x,a)
\right)
+ \text{dist}^h(\tilde{x},\tilde{a}) +
\begin{pmatrix}
\nu \\ \omega
\end{pmatrix}^\mathsf{T}
\left(\begin{pmatrix}
x \\ a
\end{pmatrix}
-
\begin{pmatrix}
\tilde{x} \\ \tilde{a}
\end{pmatrix}
\right)
+\frac{\rho}{2}
\left\|
  \begin{pmatrix} \tilde{x} \\ \tilde{a} \end{pmatrix}
- \begin{pmatrix} x \\ a \end{pmatrix}
\right\|_2^2,
\]

\noindent where $\rho > 0$ is a parameter and $\phi \defeq (x,a,\tilde{x},\tilde{a},\nu,\omega)$. Let $(\tilde{x}^{(1)}, \tilde{a}^{(1)}) \in \mathbb{R}^{6m_0(n+1)}$ and $(\nu^{(1)}, \omega^{(1)}) \in \mathbb{R}^{6m_0(n+1)}$ denote our initial guess for the state-control pair $(\tilde{x},\tilde{a}) \in \mathbb{R}^{6m_0(n+1)}$ and the associated dual variables $(\nu,\omega) \in \mathbb{R}^{6m_0(n+1)}$ at iteration $k=1$, respectively. We can now minimize the Lagrangian $\mathcal{L}_{\rho}$ over each individual pair of primal variables, while the value for the other pair of primal variables as well as the dual variables remain fixed. We obtain
\[
\begin{aligned}
\begin{pmatrix}
x^{(k+1)} \\ a^{(k+1)}
\end{pmatrix}
&=
\argmin_{x,\, a}
\mathcal{L}_{\rho}(x,a,\tilde{x}^{(k)},\tilde{a}^{(k)},\nu^{(k)},\omega^{(k)})
\\
\begin{pmatrix}
\tilde{x}^{(k+1)} \\ \tilde{a}^{(k+1)}
\end{pmatrix}
&=
\argmin_{\tilde{x},\, \tilde{a}}
\mathcal{L}_{\rho}(x^{(k+1)},a^{(k+1)},\tilde{x},\tilde{a},\nu^{(k)},\omega^{(k)})
\\
   \begin{pmatrix} \nu^{(k+1)} \\ \omega^{(k+1)} \end{pmatrix}
&= \begin{pmatrix} \nu^{(k)} \\ \omega^{(k)} \end{pmatrix}
+  \rho\left(
    \begin{pmatrix} \tilde{x}^{(k+1)} \\ \tilde{a}^{(k+1)} \end{pmatrix}
-  \begin{pmatrix} x^{(k+1)} \\ a^{(k+1)} \end{pmatrix}
\right).
\end{aligned}
\]

\noindent We can observe that the dual variables of this problem correspond to a scaled running sum of the consensus error~\cite{Boyd:2011a}. We combine the linear and quadratic terms in the augmented Lagrangian to derive the proximal version of the scheme stated above. Dropping constant terms from the augmented Lagrangian $\mathcal{L}_{\rho}$ in each step, we obtain
\[
\begin{aligned}
\begin{pmatrix}
x^{(k+1)} \\ a^{(k+1)}
\end{pmatrix}
&=
\argmin_{x,\, a}
\left(
\text{ind}_{\mathcal{C}}(x,a) + \text{kin}^h(x,a)
+
\begin{pmatrix}
\nu^{(k)} \\ \omega^{(k)}
\end{pmatrix}^\mathsf{T}
\begin{pmatrix}
x \\ a
\end{pmatrix}
+
\frac{\rho}{2}
\left\|
  \begin{pmatrix} \tilde{x}^{(k)} \\ \tilde{a}^{(k)} \end{pmatrix}
- \begin{pmatrix} x \\ a \end{pmatrix}
\right\|_2^2
\right),
\\
\begin{pmatrix}
\tilde{x}^{(k+1)} \\ \tilde{a}^{(k+1)}
\end{pmatrix}
&=
\argmin_{\tilde{x},\, \tilde{a}}
\left(
\text{dist}^h(\tilde{x},\tilde{a})
-
\begin{pmatrix}
\nu^{(k)} \\ \omega^{(k)}
\end{pmatrix}^\mathsf{T}
\begin{pmatrix}
\tilde{x} \\ \tilde{a}
\end{pmatrix}
+\frac{\rho}{2}
\left\|
  \begin{pmatrix} \tilde{x} \\ \tilde{a} \end{pmatrix}
- \begin{pmatrix} x^{(k+1)} \\ a^{(k+1)} \end{pmatrix}
\right\|_2^2
\right),
\\
   \begin{pmatrix} \nu^{(k+1)} \\ \omega^{(k+1)} \end{pmatrix}
&= \begin{pmatrix} \nu^{(k)} \\ \omega^{(k)} \end{pmatrix}
+  \rho\left(
    \begin{pmatrix} \tilde{x}^{(k+1)} \\ \tilde{a}^{(k+1)} \end{pmatrix}
-  \begin{pmatrix} x^{(k+1)} \\ a^{(k+1)} \end{pmatrix}
\right).
\end{aligned}
\]

\noindent We rescale the dual variables to obtain $(u^{(k)}, w^{(k)}) \defeq ((1/\rho)\nu^{(k)}, (1/\rho)\omega^{(k)})$, and pull the linear terms into the quadratic terms. Then, the iterative scheme for successively updating the primal and dual variables $(x, a)$, $(\tilde{x}, \tilde{a})$, $(u, w)$, is given by
\begin{subequations}\label{e:opsplit}
\begin{align}
\begin{pmatrix}
x^{(k+1)} \\ a^{(k+1)}
\end{pmatrix}
&=
\argmin_{x,\, a}
\left(
\text{ind}_{\mathcal{C}}(x,a) + \text{kin}^h(x,a)
+ \frac{\rho}{2}
\left\|\,
\begin{pmatrix} x \\ a \end{pmatrix}
-
\begin{pmatrix} \tilde{x}^{(k)} \\ \tilde{a}^{(k)} \end{pmatrix}
-
\begin{pmatrix} u^{(k)} \\ w^{(k)} \end{pmatrix}
\,\right\|_2^2
\right),
\label{e:stepkin}
\\
\begin{pmatrix}
\tilde{x}^{(k+1)} \\ \tilde{a}^{(k+1)}
\end{pmatrix}
&=
\argmin_{\tilde{x},\, \tilde{a}}
\left(
\text{dist}^h(\tilde{x},\tilde{a})
+ \frac{\rho}{2}
\left\|\,
  \begin{pmatrix} \tilde{x} \\ \tilde{a} \end{pmatrix}
- \begin{pmatrix} x^{(k+1)} \\ a^{(k+1)} \end{pmatrix}
+ \begin{pmatrix} u^{(k)} \\ w^{(k)} \end{pmatrix}
\,\right\|_2^2
\right),
\label{e:stepdist}
\\
   \begin{pmatrix} u^{(k+1)} \\ w^{(k+1)} \end{pmatrix}
&= \begin{pmatrix} u^{(k)} \\ w^{(k)} \end{pmatrix}
+  \begin{pmatrix} \tilde{x}^{(k+1)} \\ \tilde{a}^{(k+1)} \end{pmatrix}
-  \begin{pmatrix} x^{(k+1)} \\ a^{(k+1)} \end{pmatrix}
\label{e:consens},
\end{align}
\end{subequations}

\noindent with iteration index $k = 1,2,\ldots$. We can represent \eqref{e:opsplit} more concisely in its proximal form as
\[
\begin{aligned}
\begin{pmatrix}
x^{(k+1)} \\ a^{(k+1)}
\end{pmatrix}
&=
\prox_{\gamma(\text{ind}_{\mathcal{C}} + \text{kin}^h)} \left(\tilde{x}^{(k)} + u^{(k)}, \tilde{a}^{(k)} + w^{(k)}\right),
\\
\begin{pmatrix}
\tilde{x}^{(k+1)} \\ \tilde{a}^{(k+1)}
\end{pmatrix}
&=
\prox_{\gamma\text{dist}^h}\left(x^{(k+1)} - u^{(k)}, a^{(k+1)} - w^{(k)}\right),
\\
\begin{pmatrix} u^{(k+1)} \\ w^{(k+1)} \end{pmatrix}
&= \begin{pmatrix} u^{(k)} \\ w^{(k)} \end{pmatrix}
+  \begin{pmatrix} \tilde{x}^{(k+1)} \\ \tilde{a}^{(k+1)} \end{pmatrix}
-  \begin{pmatrix} x^{(k+1)} \\ a^{(k+1)} \end{pmatrix},
\end{aligned}
\]

\noindent where $\gamma \defeq 1/\rho$. Compared to the optimization problem in~\eqref{e:optprobdiscunc}, we now have to tackle two optimization problems. We see below that they are significantly simpler than the original, coupled problem in~\eqref{e:optprobdiscunc}. The first optimization problem for $(x,a)$ in~\eqref{e:stepkin} solely involves the kinetic energy $\kin^h$ in \eqref{e:kinvdisc} (along with the indicator function for the set of state-control pairs $\mathcal{C}$ in~\eqref{e:staconpair}). We refer to this subproblem as the \emph{kinetic energy subproblem}. The second optimization problem for $(\tilde{x},\tilde{a})$ in~\eqref{e:stepdist} solely involves the kernel distance $\dist^h$ in~\eqref{e:discretedist}; it minimizes the distance between the deformed shape $x(t^{n+1}) = Qx$ (i.e., the solution of the state equation for a candidate $a$) and the reference shape $x_1$. We refer to this optimization problem as the \emph{distance subproblem}. The update in~\eqref{e:consens} involves the residual of the primal variables $(x, a)$ and $(\tilde{x}, \tilde{a})$ at iteration $k+1$. As we converge to the solution of the original problem, $(x, a)$ and $(\tilde{x}, \tilde{a})$ converges to the same vectors and, consequently, we expect the dual variables $(u, w)$ to vanish at optimality. Evaluating~\eqref{e:consens} is prosaic, but solutions to~\eqref{e:stepkin} and~\eqref{e:stepdist} are nontrivial. We discuss our approaches for each of these subproblems next.

\paragraph{Kinetic Energy Subproblem}

At each iteration of our splitting algorithm, we solve~\eqref{e:stepkin} by setting up and solving the linear system associated with the first order KKT conditions. To this end, we form the Lagrangian for the quadratic, equality constrained optimization problem associated with the implicit, unconstrained formulation in~\eqref{e:stepkin}: Given current estimates $(u^{(k)}, w^{(k)})$ and $(\tilde{x}^{(k)}, \tilde{a}^{(k)})$ find $(x,a)$ such that
\begin{equation}\label{e:kineticsubproblem}
\begin{aligned}
\minopt_{x,\, a}\quad
& h a^\mathsf{T} B a + \frac{\rho}{2}
\left\|\,
\begin{pmatrix} x \\ a \end{pmatrix}
-
\begin{pmatrix} \tilde{x}^{(k)} \\ \tilde{a}^{(k)} \end{pmatrix}
-
\begin{pmatrix} u^{(k)} \\ w^{(k)} \end{pmatrix}
\,\right\|_2^2
\\
\text{subject to} \;\; & \begin{pmatrix}G^x\;\; G^a\end{pmatrix}
\begin{pmatrix}
x\\a
\end{pmatrix}
= q.
\end{aligned}
\end{equation}

As stated above, we use the method of Lagrange multipliers to handle the equality constraint. We have
\[
\mathcal{L}(x,a, \nu) =
h a^\mathsf{T} Ba + \frac{\rho}{2}
\left\|\,
\begin{pmatrix} x \\ a \end{pmatrix}
-
\begin{pmatrix} \tilde{x}^{(k)} \\ \tilde{a}^{(k)} \end{pmatrix}
-
\begin{pmatrix} u^{(k)} \\ w^{(k)} \end{pmatrix}
\,\right\|_2^2
- \nu^\mathsf{T}(G^aa +  G^xx  - q),
\]

\noindent where $\nu \in \mathbb{R}^{\tilde{n}}$, $\tilde{n} \defeq 3m_0(n+1)$ represents the vector of Lagrangian multipliers. Differentiating the Lagrangian with respect to each block of variables, we can derive the first order KKT optimality conditions for~\eqref{e:kineticsubproblem}. These state that the optimal variables $(a^{\text{opt}}, x^{\text{opt}}, \nu^{\text{opt}})$ satisfies the block linear system
\[
\begin{pmatrix}
B + \rho I_{\tilde{n}} & 0      & - (G^a)^{\mathsf{T}} \\
0          & \rho I_{\tilde{n}} & - (G^x)^{\mathsf{T}} \\
G^a        & G^x                & 0
\end{pmatrix}
\begin{pmatrix}
a^{\text{opt}} \\
x^{\text{opt}} \\
\nu^{\text{opt}}
\end{pmatrix}
=
\begin{pmatrix}
\tilde{a}^{(k)} + u^{(k)} \\
\tilde{x}^{(k)} + w^{(k)}  \\
q
\end{pmatrix}.
\]

To solve this KKT system numerically, we perform the change of variables $(a^{\text{opt}},x^{\text{opt}}) = (a + p^a, x + p^x)$. We obtain the linear system
\begin{equation}\label{e:kkt}
\begin{pmatrix}
B + \rho I_{\tilde{n}} & 0                  & (G^a)^{\mathsf{T}} \\
0                      & \rho I_{\tilde{n}} & (G^x)^{\mathsf{T}} \\
G^a                    & G^x                & 0
\end{pmatrix}
\begin{pmatrix}
-p^a \\
-p^x \\
\nu^{\text{opt}}
\end{pmatrix}
=
\begin{pmatrix}
g^a \\
g^x \\
c
\end{pmatrix}
\end{equation}

\noindent with $g^a = B a + \rho(a - \tilde{a}^{(k)} - u^{(k)})$, $g^a \in \mathbb{R}^{\tilde{n}}$, $g^x = \rho(x - \tilde{x}^{(k)} - v^{(k)})$, $g^x \in \mathbb{R}^{\tilde{n}}$, and $c =  G^a a +  G^x x  - q$, $c \in \mathbb{R}^{\tilde{n}}$. In exact arithmetic, the solution to this linear system gives the step $(p^a, p^x)$ that sums with our current approximation to the solution $(a, x)$ to get the solution $(a^{\text{opt}}, x^{\text{opt}})$ that satisfies the first order optimality conditions. Here, the vectors $g^a$ and $g^x$ correspond to the gradient of~\eqref{e:kineticsubproblem} with respect to $a$ and $x$, respectively, and $c$ corresponds to the linear constraint. This system is (in exact arithmetic) symmetric positive semi-definite.

There are several options to solve the linear system~\eqref{e:kkt}. The matrix on the left hand side of \eqref{e:kkt} is highly structured and sparse. Additionally, for the fixed kernel matrices we use, this matrix does not change during the optimization. These traits make direct factorization approaches~\cite{Duff:2017a,Davis:2006a} feasible for smaller problems; we considered a symmetric indefinite factorization~\cite[p.~454f.]{Nocedal:2006a} (i.e., an LDLT factorization~\cite{Higham:2009a,Ashcraft:1998a}) in our prior work~\cite{Zhang:2021a}. Since we fix the block matrices, this factorization needs to be computed only once (at the beginning of the algorithm). This trick provides significant computational gains. On the downside, the factorization is computational demanding for fine discretizations and quickly results in excessive memory pressure since we need to form and store~\eqref{e:kkt} and its factorization.

Another option is the Schur--complement method---a projection based approach~\cite[p.~]{Nocedal:2006a}. This method projects the problem onto the space of Lagrangian multipliers, resulting in the linear system
\begin{equation}
\label{e:shurcomplement}
GM^{-1}G^{\mathsf{T}} \nu = GM^{-1} g - c
\end{equation}

\noindent with
\[
G \defeq
\begin{pmatrix} G^a & G^x\end{pmatrix}
\in \mathbb{R}^{\tilde{n},2\tilde{n}},
\qquad
M \defeq
\begin{pmatrix}
B + \rho I_{\tilde{n}} & 0 \\
0                      & \rho I_{\tilde{n}}
\end{pmatrix}
\in \mathbb{R}^{2\tilde{n},2\tilde{n}},
\qquad
g \defeq
\begin{pmatrix} g^a \\ g^x\end{pmatrix}
\in \mathbb{R}^{2\tilde{n}}.
\]

We solve this equation iteratively using a preconditioned conjugate gradient method ({\bf PCG}). After solving for $\nu$, the step updates are recovered by evaluating
\[
\begin{pmatrix} p^a \\  p^x \end{pmatrix} =
M^{-1} \left(G^\mathsf{T} \nu - g\right).
\]

\begin{remark}
We make several remarks about the Schur--complement method we use here. First, we reiterate that the block matrices used above do not need to be computed explicitly but can be accessed via function calls to matrix-vector multiplication with a vector. This extends to $M$ and its inverse, both of which are block diagonal. We compute matrix-vector products with the inverse of this matrix in a block-wise fashion. This requires the inversion of the $2(n + 1)$ diagonal blocks. The first $n+1$ of these blocks are of the form $h B^j + \rho I_{3m_0}$ for $j = 1,\ldots,n+1$ with $B^j$ as in~\eqref{e:kernelmat}. These matrices are symmetric positive definite, and for small problems their Cholesky factors can be used to compute and store their inverses across all operator splitting iterations. For larger problems, explicit inversion may prove infeasible. This is a limitation of the current implementation, and the inversion of these blocks for larger problems represents future work. The remaining $n+1$ blocks are identity matrices scaled by the factor $\rho$, so their inversion is trivial.

We also note that solving~\eqref{e:shurcomplement} to a low tolerance using an iterative method results in steps $(p^a, p^x)$ that do not yield the exact solution $(a^{\text{opt}}, x^{\text{opt}})$ satisfying the optimality conditions. In practice, the solution only represent an update resulting in a new iterate, $(a^{(k+1)}, x^{(k+1)}) = (a^{(k)} + p^a,x^{(k)} + p^x)$. Balancing accuracy and computational load, we set the tolerance for the PCG to \num{1e-4} with an upper bound of 100 iterations.

Lastly, we note that our (modified) PCG implementation allows us to deal with optimization problems that tackle systems with matrices that are not positive definite. That is, we have added a condition that checks for negative curvature. We refer to, e.g.,~\cite[p.~168f.]{Nocedal:2006a} for details.
\end{remark}

\begin{figure}
\centering
\includegraphics[width=0.6\textwidth]{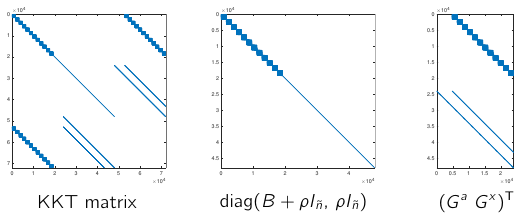}
\caption{Illustration of the KKT matrix in \eqref{e:kkt}. We show the matrix on the left and the individual blocks on the middle and right. The matrices are constructed for two representative shapes from our mitral valve databases with $m_0 = m_1 = 1\,600$.}
\end{figure}

\paragraph{Preconditioner for Reduced Space KKT System}

We have designed a block preconditioner for~\eqref{e:shurcomplement}. First, we note that
\[
G M^{-1} G^\mathsf{T}
=
G^a (B + \rho I_{\tilde{n}})^{-1} (G^a)^{\mathsf{T}}
+ \rho^{-1} G^x (G^x)^{\mathsf{T}}.
\]

\noindent Using block matrix operations, we then see that the matrix has block tridiagonal form with
\[
G M^{-1} G^\mathsf{T}
=
\frac{1}{\rho}
\begin{pmatrix}
 I_{3m_0} & -I_{3m_0} &           &           &           &           \\
-I_{3m_0} & L^1       & -I_{3m_0} &           &           &           \\
          & -I_{3m_0} & L^2       & -I_{3m_0} &           &           \\
          &           & \ddots    & \ddots    & \ddots    &           \\
          &           &           & -I_{3m_0} & L^{n-1}   & -I_{3m_0} \\
          &           &           &           & -I_{3m_0} & L^n
\end{pmatrix}
\in \mathbb{R}^{3m_0(n+1),3m_0(n+1)},
\]

\noindent where $L^j = h^2 B^j (B^j + \rho I_{3m_0})^{-1} (B^j)^{\mathsf{T}} + 2\rho I_{3m_0}$ for $j = 1, \ldots, n$. For the preconditioner, we drop the off diagonal identity blocks and invert to get
\[
P^{-1} = \operatorname{diag} \left(I_{3m_0}, L^1, \ldots, L^n \right)^{-1} \in \mathbb{R}^{3m_0(n+1), 3m_0(n+1)}
\]

\noindent as a left preconditioner. We note that $P^{-1}$ is block diagonal, so the computation of the inverse requires only the inverses of the component blocks. These can be computed offline using a Cholesky factorization and reused throughout the optimization. Furthermore, the matrix does not need to be computed explicitly and can be accessed as a function call to a matrix-vector product.%

\paragraph{Distance Subproblem}

Next, we develop a numerical method to solve subproblem~\eqref{e:stepdist}. The associated unconstrained optimization problem is as follows: Given current estimates $(x^{(k+1)}, a^{(k+1)})$, $(u^{(k)}, w^{(k)})$ find $(\tilde{x}, \tilde{a})$ such that
\begin{equation}\label{e:distsubproblem}
\minopt_{\tilde{x},\, \tilde{a}}\,\,
\text{dist}^h(\tilde{x},\tilde{a})
+ \frac{\rho}{2}
\left\|\,
  \begin{pmatrix} \tilde{x} \\ \tilde{a} \end{pmatrix}
- \begin{pmatrix} x^{(k+1)} \\ a^{(k+1)} \end{pmatrix}
+ \begin{pmatrix} u^{(k)} \\ w^{(k)} \end{pmatrix}
\,\right\|_2^2.
\end{equation}

Notice that the only variable appearing in the distance functional is $x^{n+1} = x(t^{n+1}) \in \mathbb{R}^{3m_0}$, i.e., the ``transported'' shape/state variable at the terminal time point $t^{n+1} = 1$. The remaining variables are included only in the mismatch term involving a 2-norm. Consequently, we can compute the minimizer for these variables explicitly as
\[
\tilde{a}^{(k+1)} = a^{(k+1)} - w^{(k)}
\quad \text{and} \quad
\tilde{x}^{(k+1)}(t^j) = x^{(k+1)}(t^j) - u^{(k)}(t^j) \quad \text{for } j = 1,\ldots,n.
\]

The remaining problem for $\tilde{x}$ at time $t^{n+1}$ with $z \defeq \tilde{x}^{n+1}$ is
\begin{equation}\label{e:distancesub}
\minopt_{z}\,\,
\dist^h(z,x_1)
+ \frac{\rho}{2} \| z - x^{(k+1)}(t^{n+1}) + v^{(k)}(t^{n+1}) \|^2_2
\end{equation}

We solve the unconstrained optimization problem in~\eqref{e:distancesub} at each outer iteration using a matrix-free Newton--Krylov method~\cite{Nocedal:2006a}. This requires first order (gradient) and second order (curvature) derivative information. The derivatives of the squared-norm are straightforward to obtain and therefore omitted. The first and second derivative of the distance as defined in~\eqref{e:discretedist} with respect to an arbitrary vector $z\in \mathbb{R}^{3m_0}$ is formally given by
\[
\d_{z} \dist^h(z,x_1)
= \frac{\alpha}{2}\d_z  \left(\phi(z,z) - 2 \phi(z,x_1) + \phi(x_1,x_1)\right)
= \frac{\alpha}{2}\left(\d_z  \phi(z,z) - 2 \d_z\phi(z,x_1)\right)
\]

\noindent and
\[
\d_{zz} \dist^h(z,x_1)
=\frac{\alpha}{2}\left(\d_{zz}  \phi(z,z) - 2 \d_{zz}  \phi(z,x_1)\right),
\]

\noindent respectively. Starting with the second term of the sum, the first derivative of $\phi$ with respect to the vector $z \in \mathbb{R}^{3m_0}$ is given by
\[
\d_z \phi(z,x_1)
= -\frac{1}{m_0m_1\sigma^2_s} P^{\mathsf{T}}_1 D^{\mathsf{T}}  S^{\mathsf{T}}
\exp \left(-\frac{1}{2\sigma^2_s} S(d(z,x_1) \odot d(z, x_1))\right)
\in \mathbb{R}^{3m_0},
\]

\noindent where $D \defeq \operatorname{diag}(d(z,x_1))$ denotes a diagonal matrix with the entries of the vector $d(z,x_1)$ on the main diagonal. The second derivative is given by the $3m_0\times3m_0$ matrix
\[
\begin{aligned}
\d_{zz} \phi(z,x_1)
& = \frac{1}{m_0m_1\sigma_s^4} P^{\mathsf{T}}_1 D^{\mathsf{T}} S^{\mathsf{T}} \operatorname{diag}\left(\exp \left(-\frac{1}{2\sigma_s^2} S\left(d(z,x_1) \odot d(z, x_1)\right)  \right) \right)S D P_1 \\
& \quad - \frac{1}{m_0m_1\sigma^2_s} P^{\mathsf{T}}_1 \operatorname{diag}\left(S^{\mathsf{T}}  \exp \left(-\frac{1}{2\sigma_s^2} S (d(z,x_1) \odot d(z, x_1))  \right)\right) P_1,
\end{aligned}
\]

The derivatives for $\phi(z,z)$ are computed similarly with the exception that the all-to-all distance vector $d$ (see \secref{s:discdist} for the definition), becomes $d(z,z) = P_1 z - P_2 z = (P_1 - P_2)z$. Thus, the first and second derivative for $\phi(z,z)$ are given by the vector
\[
\d_z \phi(z,z)
= -\frac{1}{m^2_0\sigma^2_s} (P_1 - P_2)^{\mathsf{T}} D^{\mathsf{T}} S^{\mathsf{T}}
\exp \left(-\frac{1}{2\sigma^2_s} S (d(z,z) \odot d(z, z))\right) \in \mathbb{R}^{3m_0}
\]

\noindent and the $3m_0 \times 3m_0$ matrix
\[
\begin{aligned}
\d_{zz} \phi(z,z)
& = \frac{1}{m_0^2\sigma^4_s} (P_1 - P_2)^{\mathsf{T}} D^{\mathsf{T}} S^{\mathsf{T}}
\operatorname{diag}
\left(\exp \left(-\frac{1}{2\sigma^2_s} S (d(z,z) \odot d(z,z))  \right) \right)SD\,(P_1 - P_2) \\
& \quad - \frac{1}{m_0^2\sigma^2_s}  (P_1 - P_2)^{\mathsf{T}}
\operatorname{diag}\left(S^{\mathsf{T}}
\exp \left(-\frac{1}{2\sigma^2_s} S (d(z,z) \odot d(z,z)\right)
\right) (P_1 - P_2)
\end{aligned}
\]

\noindent with $D \defeq \operatorname{diag}(d(z,z))$, respectively.

\begin{remark}
It is inefficient to store the Hessian pieces $\d_{zz} \phi(z,z)$ and $\d_{zz} \phi(z,x_1)$ explicitly. Instead, they are passed as function handles, evaluating their product with a given vector $z$.
\end{remark}

The Newton step $s^{(k)} \in \mathbb{R}^{3m_0}$ at outer iteration $k$ is given by
\begin{equation}
\label{e:newtonstepdist}
H^{(k)} s^{(k)} = -g^{(k)},
\quad
\tilde{x}^{(k+1)}(t^{n+1}) = \tilde{x}^{(k)}(t^{n+1}) + t^{(k)} s^{(k)},
\end{equation}

\noindent where $H \defeq \d_{zz} \phi(z,z) + \d_{zz} \phi(z,x_1) \in \mathbb{R}^{3m_0,3m_0}$, $g \defeq \d_z \phi(z,z) + \d_z \phi(z,x_1) \in \mathbb{R}^{3m_0}$ are the Hessian and gradient of the objective function with $\d_{zz} \phi(z,z)$, $\d_{zz} \phi(z,x_1)$, $\d_z \phi(z,z)$ and $\d_z \phi(z,x_1)$ as defined above, and $t^{(k)} > 0$ is a line search parameter. We use a backtracking line search subject to the Armijo condition~\cite{Nocedal:2006a,Boyd:2004a}. The solution to this linear system at iteration $k$ is the search direction $s^{(k)}$ for $\tilde{x}$ at time $t^{n+1} = 1$. We solve the linear system $H^{(k)} s^{(k)} = -g^{(k)}$ in~\eqref{e:newtonstepdist} iteratively using a matrix-free PCG method. We use a \emph{quadratic forcing sequence} to compute the tolerance $\epsilon_{\text{PCG}}$ for the PCG method. The term \emph{forcing sequence} refers to an adaptive strategy for selecting the tolerance of iterative methods used to solve for the search direction $s^{(k)}$ in \eqref{e:newtonstepdist} that is commonly used in large scale optimization problems. Instead of solving~\eqref{e:newtonstepdist} at every iteration $k$ to high accuracy (i.e., to machine precision, as we would do when considering direct methods), we choose a tolerance that is proportional to the norm of the gradient $g^{(k)}$. This way, we gradually increase the accuracy for the computation of the search direction $s^{(k)}$. Since the overall rate of convergence may not benefit from an accurate estimate of $s^{(k)}$ far away from an optimal point, this allows us to significantly reduce the computational burden. Assuming (optimal) quadratic convergence of our Newton--Krylov solver, we set $\epsilon_{\text{PCG}}$ to $\epsilon_{\text{PCG}} = \min\left(\|g^{(k)}\|_2/\|g^{(0)}\|_2, \tau_{\text{PCG}} \right)$, with $\tau_{\text{PCG}} = \nicefrac{1}{4}$. If we were to assume superlinear convergence we can slightly relax this tolerance by replacing $\|g^{(k)}\|_2/\|g^{(0)}\|_2$ with its square root. We note that other options for selecting $\tau_{\text{PCG}}$ exist. We refer, e.g., to~\cite{Nocedal:2006a,Eisenstat:1996a} for additional details. Notice that the normalization with respect to the initial gradient $\|g^{(0)}\|_2$ at iteration $k = 0$ is not the standard textbook approach. Based on our past experience with solving large scale inverse problems, we prefer this choice since it allows us to avoid issues with the scaling of the gradient $g$.

\subsubsection{Stopping Criteria}

We now discuss the stopping criteria for the operator-splitting iteration introduced in~\eqref{e:opsplit}. We separate this discussion into two parts: a discussion of the algorithm convergence of the operator-splitting iteration and a discussion of the convergence of the mapped template shape to the reference shape.

We monitor the convergence of the splitting algorithm using the primal and dual residuals $r_{\text{prim}}$ and $r_{\text{dual}}$, which are given by
\begin{equation}\label{e:residuals}
r_{\text{prim}}^{(k)} =
  \begin{pmatrix} a^{(k)} \\ x^{(k)} \end{pmatrix}
- \begin{pmatrix} \tilde{a}^{(k)} \\ \tilde{x}^{(k)}  \end{pmatrix}
\quad
\text{and}
\quad
r_{\text{dual}}^{(k)} =
\rho
\left(
  \begin{pmatrix} \tilde{a}^{(k)} \\ \tilde{x}^{(k)} \end{pmatrix}
- \begin{pmatrix} \tilde{a}^{(k-1)} \\ \tilde{x}^{(k-1)} \end{pmatrix}
\right),
\end{equation}

\noindent respectively. Here, the primal residual monitors the convergence of variables to a consensus solution, while the dual residual monitors the convergence to a fixed solution across iterations. Both residuals can be shown to converge to zero~\cite{Parikh:2013a}. Thus, a suitable criteria for convergence of the splitting method in \eqref{e:opsplit} is to stop when the norms of $r_{\text{prim}}$ and $r_{\text{dual}}$ falls below some tolerance $\epsilon_l$, $l \in \{\text{prim}, \text{dual}\}$, respectively. Guidelines for selecting $\epsilon_{\text{prim}}$ and $\epsilon_{\text{dual}}$ are provided in~\cite{Boyd:2011a,Parikh:2013a}.

While these criteria assess the convergence of the optimization algorithm, we are ultimately interested in matching the template shape $x_0$ to the reference shape $x_1$. Consequently, we add an additional convergence criteria that evaluates the discrepancy between the deformed template shape $x$ at $t^{n+1}$ and the reference shape $x_1$. We consider the Hausdorff distance~\cite{Azencott:2010a,Zhang:2021a}. For two generic shapes $x = (x^1, \ldots, x^{m_0}) \in \mathbb{R}^{3m_0}$ and $z = (z^1, \ldots, z^{m_1}) \in \mathbb{R}^{3m_1}$, the Hausdorff distance is defined as
\begin{equation}\label{e:hausdist}
\dist_{H}(x,z) = \max \left\{ d_H(x,z),d_H(z,x)\right\},
\end{equation}

\noindent where $d_H(x,z) = \max_{i \in \{1,\ldots,m_0\}} \left\{ \min_{j \in \{1,\ldots,m_1\}} \|x^i - z^j \|_2 \right\}$. In practice, we use the 95th percentile of the Hausdorff distance to guard against outliers. From a mathematical point of view, this censoring is not required. However, we consider this confidence interval to account for the fact that in practical applications we cannot expect the data to be represented accurately, due to measurement errors, noise, and a limited resolution as well the modelling errors in their representation. Additional sources of errors stem from our problem formulation along with the numerical discretization. Consequently, we do not expect the matching to be perfect everywhere. Adding a 95\% confidence interval allows us to account for these inaccuracies. We terminate the algorithm if this censored Hausdorff distance falls below a user defined threshold $\epsilon_{\text{haus}} > 0$. Moreover, we monitor the update of the Hausdorff distance across five consecutive iterations. Lastly, to safeguard against excessive runtimes, we terminate the optimization if the number of iterations exceeds a user defined upper bound $n_{\text{iter}}$.

In summary, the stopping conditions are
\begin{subequations}
\label{e:stopcond}
\begin{align}
    (C1)\qquad & \dist_{H}(Qx^{(k)},x_1) < \epsilon_{\text{haus}}, \\
    (C2)\qquad & \textstyle\sum_{j=0}^4 |\delta_{\text{haus}}(j)|  < \epsilon_{\text{haus}} / 10^3,\\
    (C3)\qquad & \| r_{\text{prim}}^{(k)} \|_2 < \epsilon_{\text{prim}},\\
    (C4)\qquad & \| r_{\text{dual}}^{(k)} \|_2 < \epsilon_{\text{dual}},\\
    (C5)\qquad & k > n_{\text{iter}},
\end{align}
\end{subequations}

\noindent with $\delta_{\text{haus}}(j) \defeq \dist_{H}(Qx^{(k-j)},x_1) - \dist_{H}(Qx^{(k-(j+1))},x_1)$. We terminate the algorithm if
\begin{equation}\label{e:stopcond-bool}
    (C1) \lor (C2) \lor (C3) \lor (C4)  \lor (C5).
\end{equation}

\noindent We overview the entire algorithm in \ialgref{a:splitting}.

\begin{algorithm}
\caption{Pseudo algorithm for operator splitting approach. We summarize the values and heuristic/empirical policies to select algorithmic and problem specific parameters in~\tabref{t:parameters}.\label{a:splitting}}
\begin{algorithmic}[1]
\State {\bf Input:} $x_0, x_1, a^{(0)} = 0$
\State {\bf Parameters:} $n_{\text{iter}}$, $\epsilon_{l} > 0$, $l \in\{\text{prim}, \text{dual}, \text{haus}\}$, $\rho > 0$, $\alpha > 0$, $\sigma_v, \sigma_s > 0$, $n \in \mathbb{N}$
\State $k \gets 0$
\While{$\neg$ stop}
\State $\begin{pmatrix}x^{(k+1)}, a^{(k+1)}\end{pmatrix} \gets$ given $\begin{pmatrix}u^{(k)}, w^{(k)}\end{pmatrix}$ and $\begin{pmatrix}\tilde{x}^{(k)}, \tilde{a}^{(k)}\end{pmatrix}$ solve kinetic energy subproblem \eqref{e:kineticsubproblem}
\State $\begin{pmatrix}\tilde{x}^{(k+1)}, \tilde{a}^{(k+1)}\end{pmatrix} \gets$ given $\begin{pmatrix}x^{(k+1)}, a^{(k+1)}\end{pmatrix}$, $\begin{pmatrix}u^{(k)}, w^{(k)}\end{pmatrix}$ solve distance subproblem \eqref{e:distsubproblem}
\State $\begin{pmatrix} u^{(k+1)} \\ w^{(k+1)} \end{pmatrix}
\gets \begin{pmatrix} u^{(k)} \\ w^{(k)} \end{pmatrix}
+  \begin{pmatrix} \tilde{x}^{(k+1)} \\ \tilde{a}^{(k+1)} \end{pmatrix}
-  \begin{pmatrix} x^{(k+1)} \\ a^{(k+1)} \end{pmatrix}$
\State $\begin{pmatrix}r_{\text{prim}}^{(k+1)}, r_{\text{dual}}^{(k+1)}\end{pmatrix} \gets$ evaluate \eqref{e:residuals}
\State $k \gets k + 1$
\State stop $\gets$ evaluate \eqref{e:stopcond-bool}
\EndWhile
\State {\bf Output:} $a^{(k)}$
\end{algorithmic}
\end{algorithm}

\begin{remark}
We can use various measures to assess the proximity between geometric structures. In the present case, we are--in essence--interested in matching point clouds. To measure the proximity between points in $\mathbb{R}^3$, the Hausdorff  distance is a natural choice since it has an intuitive geometric interpretation. This is not the case for the kernel distance we consider for optimization. Due to its geometric interpretation, we prefer using the Hausdorff distance as a stopping criterion. Conversely, while the use of the Hausdorff distance to measure the proximity of the deformed shape $y \app s_0$ and the target shape $s_1$ may seem natural, it poses significant challenges in the context of variational optimization (from a theoretical and a practical point of view); it is not smooth with respect to small perturbations of the compared shapes. The kernel distance is. As such, we prefer it for optimization. A more detailed discussion can, e.g., be found in~\cite{Azencott:2010a}.
\end{remark}

\subsubsection{Parameter Choices}\label{s:parameterchoices}

We have to select several parameters. We summarize these parameters in \tabref{t:parameters}. This section includes a brief description of each of those parameters and details how they were selected for our numerical experiments.

\paragraph{Number of Time Steps} The number of cells $n \in \mathbb{N}$ used to discretize the time interval $[0,1]$ into $n+1$ time steps has to be selected as a trade-off between accuracy requirements and computational complexity. In our implementation, it affects the discretization of the kinetic energy and the ODE constraint in~\eqref{e:varoptprob}. As stated in \secref{s:disckinv}, we compute the integral with respect to time using a left Riemann summation~\cite{Quarteroni:2010a}. As stated in \secref{s:timeintegration}, we discretize the ODE using the forward Euler method. It follows that the error associated with the time integration is of order $h = 1 / (n+1)$. Consequently, a finer discretization in time increases the numerical accuracy in a linear way. Increasing the number of time steps increases the computational cost of the problem. Thus, although the kinetic energy and the discretization of the ODE both benefit from a finer discretization, this is tempered by the need to keep a computationally tractable size for $n$. If we register two arbitrary shapes (e.g., shapes representing organs of two distinct individuals), we have to select an adequate number of time steps. For the experiments in this paper, we set $n = 5$. We found empirically that $n=5$ yields a good trade off between computational complexity and accuracy requirements in the considered application. Another scenario we consider in the present work is the diffeomorphic matching of a time series of smooth shapes $s_i$, $i = 0,\ldots,n_f-1$. In this case, we, e.g., monitor (for an individual patient) the motion of an organ of interest in time. (In our case, the motion of the leaflets of mitral valves during the cardiac cycle.) Consequently, the data consists of a fixed number of frames $n_f$. The task is to reconstruct the motion between the individual frames (or, e.g., from frame at pseudo-time $t = 0$ (first frame) to frame at pseudo-time $t=1$ (last frame)). In this case, we also have to select an adequate number of time steps for the discretization. We select to set $n$ equal to the number of frames $n_f$.

\paragraph{Weight for Distance Measure} The problem also requires the selection of the parameter $\alpha > 0$. This parameter balances the two terms of the objective functional in~\eqref{e:varoptprob}: the kinetic energy enforces smoothness of the velocity and the distance term which enforces the matching between the deformed template shape and the reference shape. Selecting this parameter is a much researched topic in the field of inverse problems, and common selection methods include generalized cross-validation, L-curve, Morozov's discrepancy principle, or the unbiased predictive risk estimator method; see e.g.~\cite{Hansen:2010a,Hansen:1998a,Hansen:1992a,Vogel:2002a}. In~\cite{Mang:2015a,Azencott:2010a} we have developed parameter-continuation schemes to select appropriate values for the hyperparameters that balance the kinetic energy (regularization) and the distance measure. The basic idea is to solve a family of optimization problems parameterized by the weight $\alpha$, starting with small values for $\alpha$. Consequently, in a first step we essentially optimize the kinetic energy (i.e., we basically solve a convex, quadratic optimization problem). We can expect quick convergence. Subsequently, we increase $\alpha$ to improve the matching. We use the solution from the former step as an initialization (``warm start''). For simplicity of our numerical study presented here, we do not consider this scheme.

\begin{table}
\caption{Algorithmic parameters. We select these parameters empirically. We describe some of the policies we have developed to select them in \secref{s:parameterchoices}.\label{t:parameters}}
\centering
\begin{small}
\begin{tabular}{llll}\toprule
\bf Symbol                      & \bf Meaning                                                                           & \bf Value \\\midrule
$n$                             & number of time steps (time-discretization)                                            & 5 \\
$\alpha$                        & parameter controlling the contribution of the distance~\eqref{e:discretedist} & 1 \\
$\sigma_v$                      & bandwidth of Gaussian kernel for parameterization of $v$ (see \eqref{e:velrkhs}       & see \eqref{e:bandwidthvel} \\
$\sigma_s$                      & bandwidth of Gaussian kernel for kernel distance \eqref{e:discretedist}               & see \eqref{e:bandwidthdis} \\
$\rho$                          & operator-splitting parameter (see \eqref{e:opsplit})                                  & 1 \\
\midrule
\multicolumn{3}{l}{stopping conditions}\\\midrule
$\epsilon_{\text{prim}}$        & stopping tolerance for primal residual                                                & \num{1e-3} \\
$\epsilon_{\text{dual}}$        & stopping tolerance for dual residual                                                  & \num{1e-3} \\
$\epsilon_{\text{haus}}$        & tolerance for censored Hausdorff distance in \eqref{e:hausdist}                       & see \eqref{e:tolhaus} \\
$n_{\text{iter}}$               & maximum number of iterations                                                          & 100 \\
\bottomrule
\end{tabular}
\end{small}
\end{table}

\paragraph{Standard Deviation for Kernel Matrices} The parameterization of the velocity as elements of an RKHS in~\eqref{e:velrkhs} as well as the kernel distance in~\eqref{e:discretedist} require the selection of the scale parameters/bandwidths $\sigma_v > 0$ and $\sigma_s > 0$, respectively. If the scale parameter $\sigma_v$ increases, the smoothness of the optimal diffeomorphic deformation increases. While smooth deformations are desirable, too large values for $\sigma_v$ may lead to a poor matching quality. Conversely, $\sigma_s$ controls how many points in a local neighborhood $\mathcal{N}(x^i_0) \defeq \{ z \in \mathbb{R}^3: \|x_0^i - y\|_2^2 \leq r, r > 0, x_0^i \in x_0\}$, $i \in \{1,\ldots,m_0\}$, of radius $r > 0$ of a point $x_0^i \in x_0$ influence the movement of point $x^i_0$ of the template shape $x_0$. If $\sigma_s$ is chosen too small, the Gaussian kernels are too narrow; few points contribute locally to the distance functional in~\eqref{e:discretedist} and, thus, a point $x_0^i$ may not move towards points in the target shape $x_1$. In contrast, choosing $\sigma_s$ too large means the Gaussian kernels are too wide. This causes individual points $x_0^i$ to be affected by too many points. This increases the difficulty of optimization, and may result in a poor matching quality, as well. Consequently, policies to control the selection of these parameters are important; we select policies that take into account the geometry of the problem. In particular, they consider the average distance between points in the point cloud (or mesh size) and strike a balance between local and global matching of shapes.

We consider criteria similar to the ones proposed in~\cite{Azencott:2010a}. Recall that the shapes $s_i$, $i=0,1$, are parameterized as a set of points/landmarks $x_i \in \mathbb{R}^{3m_i}$. For these sets of points $x_i$, we can compute a triangulation to approximate the surface of $s_i$. We select $\sigma_v$ to be proportional to the mean edge length $\bar{h}_0$ of the Delaunay triangulation of the template shape $x_0$. That is, let $\mathcal{E}(x^i_0)$ denote a set of all points $x^j_0$, $j \not= i$, of $x_0$ that share an edge with $x^i$. Then,
\[
\bar{h}_0 \defeq
\operatorname{mean}\left(
\bigcup_{i=1}^{m_0} \left\{\| x^j_0 - x^i_0\|_2^2 : x^j_0 \in \mathcal{E}(x^i_0), j = 1,\ldots,m_0\right\}
\right).
\]

We have
\begin{equation}\label{e:bandwidthvel}
\sigma_v = \tau_v\, 2^{-\nicefrac{1}{2}}\,\bar{h}_0,
\end{equation}

\noindent where $\tau_v > 0$ is a user selected parameter.

The bandwidth of the distance is selected to be proportional to the Hausdorff distance between the shapes $x_0$ and $x_1$ according to
\begin{equation}\label{e:bandwidthdis}
\sigma_s = \max\left(\bar{h}_1,\tau_s\,\dist_H\left(x_0,x_1\right)/2\right),
\end{equation}

\noindent where $1\leq \tau_s \leq 2$ is a user defined parameter, $\dist_H$ is the Hausdorff distance~\eqref{e:hausdist} and $\bar{h}_1$ is the average mesh size of the target shape $x_1$. We note that we observed numerical issues if $\sigma_s$ is at the order of or below the average mesh size of the target shape (the computed search direction is no longer a descent direction and the convergence of the solver deteriorates). This is the reason for the lower bound in \eqref{e:bandwidthdis}. We conduct numerical experiments to illustrate the behavior of our methodology for different choices of $\tau_v$ and $\tau_s$, respectively.

\begin{remark}
One possibility to make the selection of the bandwidth completely automatic is to use a parameter continuation in $\sigma_v$ and $\sigma_s$ (or $\tau_v$ and $\tau_s$, respectively; i.e., a multiscale approach). The question on how to combine this multiscale strategy with a parameter continuation in $\alpha$ remains subject to future work.
\end{remark}

\paragraph{Stopping Tolerances} As can be seen in condition (C1) of~\eqref{e:stopcond}, the stopping tolerance for the censored Hausdorff distance is selected proportional to the average mesh size $\bar{h}_1$ of the discretized target shape $x_1$. We use
\begin{equation}\label{e:tolhaus}
\epsilon_{\text{haus}} = \tau_{\text{haus}}\bar{h}_1, \quad \tau_{\text{haus}} > 0.
\end{equation}

\noindent We found $\tau_{\text{haus}} = \nicefrac{1}{2}$ yields a good trade-off between accuracy and runtime requirements. The condition (C2) monitors the change of the Hausdorff distance across five iterations. That is, if we do not make significant progress in five consecutive iterations we terminate the solver. To avoid additional parameters, we select the distance updates as a fraction of the tolerance $\epsilon_{\text{haus}}$ used for the Hausdorff distance.

Since we are mostly interested in driving the mismatch between the data to zero, the remaining tolerances are treated as safeguards to avoid excessive computations. We select them as stated in \tabref{t:parameters}.

\subsection{Quantitative Measures of Shape Variability}\label{s:shapevariability}

After computing the optimal diffeomorphic map $f^{\text{opt}}=(f^{\text{opt}}_t)$ (state variable) that approximately minimizes~\eqref{e:varoptprob} and by that matching two surfaces $s_0, s_1 \in \mathcal{S}$, the terminal $\mathbb{R}^3$-diffeomorphism $y \defeq f^{\text{opt}}_1 = f^{\text{opt}}(t = 1,\cdot)$ is a smooth invertible map from $s_0$ onto a smooth surface $y \app s_0 = y(s_0) \approx s_1$. In a Riemannian setting, we can---according to~\eqref{e:shapemetric}---compute the distance between two arbitrary shapes $s_i$, $i=0,1$, by evaluating the kinetic energy~\eqref{e:kinv} for a minimizer $v^{\text{opt}} = (v^{\text{opt}}_t)$ of~\eqref{e:kinv} that parameterizes the diffeomorphic flow $f_t$ that maps $s_0$ to $s_1$. We refer to~\cite{Trouve:1995a,Miller:2002a,Younes:2020a,Bauer:2014a} for a rigorous mathematical treatment. Intuitively, this corresponds to the minimal energy required to map $s_0$ to $s_1$. Clearly, even in exact arithmetic we can only hope to compute an approximation to this distance, since the considered control problem is high-dimensional after discretization, nonconvex, and nonlinear. Moreover, the computation of the distance is also affected by the implementation choices we made. Aside from numerical tolerances or numerical accuracy of the carried out computations, any changes of the parameterization of the curved shapes, any reparametrization of these shapes, or any changes in the RKHS structure for the parameterization of the velocity yields different results. Moreover, computing geodesic distances between two shapes---even if we assume that the shapes $s_0, s_1$, are in the orbit associated with the group action of the diffeomorphisms generated here---assumes an exact matching, i.e., $y \app s_0 = s_1$. Given the issues described above as well as perturbations introduced by discretization and measurement errors, this assumption does not hold in practice. We can only hope for an approximate matching $y \app s_0 \approx s_1$. Lastly, the minimal kinetic energy involves averages of squared velocities, and hence only provides a global dissimilarity between $s_0$ and $s_1$. Strain analysis of $f$, which we now outline as in~\cite{Zhang:2021a,Dabirian:2022a}, generates the spatial distribution of local distortions between $s_0$ and $s_1$.

Let $x \in s_0$ denote any fixed point and let $\tilde{x} = y(x)$, where $y = f_1$ again denotes the terminal diffeomorphism at $t=1$. (For simplicity, we limit this exposition to the terminal diffeomorphism; the computation naturally translates to all time points $t$ of the flow $f_t$). Moreover, let $\mathcal{T}_x$ and $\mathcal{T}_{\tilde{x}}$ denote the tangent spaces of the initial shape $s_0$ and the deformed shape $\tilde{s}_0$ at $x$ and $\tilde{x}$, respectively, endowed with local surface metrics. Since $y \in \mathcal{G}$ is a smooth bijection, the differential of $y$ is an invertible linear map $J_x : \mathcal{T}_x \to \mathcal{T}_{\tilde{x}}$, $J_x \in\mathbb{R}^{2,2}$. Let $u \in \mathcal{T}_x$ denote any tangent vector with $\|u\| = 1$. With this, we can define the directional strain at $x$ associated with the map $y$ in the direction of $u$ as the dilation/contraction factor $q_{\text{dir}}(x,u) = \|J_x u\|$. Let $J_x^\ast : \mathcal{T}_{\tilde{x}} \to \mathcal{T}_x$ denote the adjoint operator of $J_x$. Moreover, let $\lambda_x^i > 0$, $i = 1, 2$, with $\lambda^1_x \leq \lambda_x^2$ denote the eigenvalues of $J_x^* J_x \succ 0$, with $\sqrt{\lambda_x^i}$, $i=1,2$, equal to the minimal and maximal directional strains at $x$, respectively. To avoid anisotropies associated with the definition of $q_{\text{dir}}$, we consider the isotropic strain $q_{\text{iso}}$ defined as
\[
\textstyle q_{\text{iso}}(x)
= \left(\prod_{i=1}^2\lambda_x^i\right)^{\nicefrac{1}{2}}
= \sqrt{| \det(J_x)|},
\]

\noindent instead. It corresponds to the local isotropic change of length for the tangent vector. We note that $q_{\text{iso}}$ has a simple geometric interpretation. For fixed $x \in s_0$, and any open patch $p_x \subset s_0$ around $x$, define the ratio of surface areas
\[
\text{area}_{\text{rel}}(p_x)= \frac{\text{area}(y(p_x))}{\text{area}(p_x)}.
\]

\noindent Then, $\text{area}_{\text{rel}}$ tends to the square of $q_{\text{iso}}$ at $x$ when the diameter of $p_x$ tends to 0. This observation allows us to compute strain directly from the triangulation after discretization of $s_0$. Since $q_{\text{iso}}$ is a dimensionless quantity of average dilation/contraction around $x$, we convert it into an ``isotropic strain intensity'' defined as
\[
p_{\text{iso}}(x) = |q_{\text{iso}}(x) - 1|.
\]

In~\cite{Zhang:2021a} we have used measures derived from the kinetic energy and strain intensity to study differences between datasets of normal and diseased patients. In~\cite{Dabirian:2022a}, we have extended this idea by exploiting these quantitative measures of shape variability as features for the automatic classification of shapes and shape deformations based on random forests. Here, we merely consider isotropic strain as a visual tool to evaluate complex diffeomorphic deformations.

\section{Numerical Experiments}\label{s:experiments}

In the following, we report results for the numerical software developed by our group for solving the diffeomorphic shape matching problem formulated in~\eqref{e:varoptprob}. The application is cardiac imaging, and, in particular, surface representations of mitral valves of individual patients extracted from image sequences.

\subsection{Data}

\begin{figure}
\centering
\includegraphics[width=0.4\textwidth]{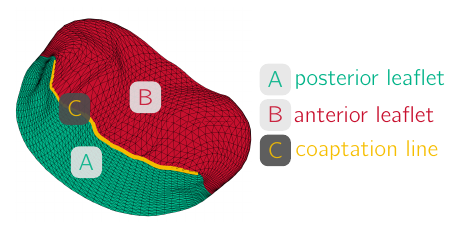}
\caption{Main anatomical structures for a representative mitral valve dataset. We show a Delaunay triangulation of the mitral valve surface. We highlight the main anatomical regions: the posterior leaflet (teal color), the anterior leaflet (red color), and the coaptation line (gold). Each leaflet is discretized by a grid of 800 points, resulting in a total of 1,600 points per 3D surface. We show the Delaunay triangulation of this mesh. (Figure modified from \cite{Dabirian:2022a}.)\label{f:mitral-valve-regions}}
\end{figure}

For a total of 150 patients, transesophageal echocardiography provides dynamic 3D views of their mitral valves at rates of roughly 25 frames per heart cycle with potential mitral valve complications. Computerized segmentation of these 3D images extracts twenty-five 3D snapshots of the mitral valve surface. For time frame $t$, the extracted surface displays the anterior and posterior leaflets of the mitral valve as smooth surfaces discretized by a grid of 800 points per leaflet. Consequently, each 3D surface is discretized by 1\,600 points per frame. These data were prepared and annotated by our longtime collaborators at the DeBakey Heart and Vascular Center of the Houston Methodist Hospital. During each heart cycle, the leaflets close the mitral valve at mid-systole ($t = t_{\text{ms}}$) and open the mitral valve at end-systole ($t = t_{\text{es}}$). The two leaflets define a deformable connected 3D-surface bounded by a flexible ring (the ``annulus''). At $t = t_{\text{ms}}$, the leaflets are in full contact along the \emph{coaptation line} in order to tightly close the mitral valve. The leaflets open the mitral valve progressively until end-systole, then remain fully open during diastole, and start closing again at beginning of systole. When the mitral valve is open, the coaptation line is split into two curved boundary segments sharing the same endpoints. We show the main anatomical regions of these data in \figref{f:mitral-valve-regions}.

In previous studies~\cite{Azencott:2010a,Zhang:2021a,ElTallawi:2021a,ElTallawi:2021b,ElTallawi:2019a}, diffeomorphic shape matching of the mitral valve leaflets between $t = t_{\text{ms}}$ and $t = t_{\text{es}}$ (a problem referred to as ``multi-frame registration'' below) was developed and systematically implemented in order to first reconstruct the smooth deformation $t \to f_t$ of the mitral valve leaflets between mid systole and end systole. We then computed the intensities and spatial distribution of the tissue strain induced by mitral valve deformation between mid systole and end systole (see \secref{s:shapevariability}). This research project aimed to provide cardiologists with patient specific displays of mitral valve leaflet strain intensities as a potential aid to evaluate/compare mitral valve clinical cases~\cite{ElTallawi:2021a,ElTallawi:2021b,ElTallawi:2019a}. In a recent followup paper, we develop diffeomorphic deformation techniques for automatic classification of soft smooth shapes~\cite{Dabirian:2022a}. We have used the described database of 3D mitral valve snapshots as a benchmark to implement and test methodology  for automatic classification of smooth surfaces; more precisely, to discriminate patients diagnosed with \emph{mitral valve regurgitation} from \emph{normal} cases. In the present work, we consider these 3D mitral valve snapshots to study the numerical performance of our improved  numerical solver for the operator splitting problem described above.

\subsection{Bandwidth Selection}\label{s:bandwithselection}

We perform experiments to study the effect of the parameters $\tau_s$ and $\tau_v$ for the bandwidth $\sigma_s$ of the kernel distance in~\eqref{e:discretedist} and $\sigma_v$ for the parameterization of the velocity in~\eqref{e:velrkhs}, respectively. We observe that these parameters have a significant effect on the performance of our method. To obtain comparable results, we execute the solver for a fixed number of $n_{\text{iter}}=100$ iterations; i.e., we neglect all convergence criteria in~\eqref{e:stopcond} with the exception of the upper bound on the number of iterations. This also allows us to explore adequate values for the tolerances for the stopping conditions listed in \eqref{e:stopcond} for the considered problem (registration of mitral valve data). Notice that these experiments are different from our previously reported results~\cite{Zhang:2021a}. In \cite{Zhang:2021a} we performed the registration in a single patient in time (frame-to-frame). Here, we fix the frame to the time point corresponding to mid-systole and register two data sets from different patients. We expect this to be a more challenging problem due to an expected increase in shape variability across patients.

\paragraph{Setup} We limit the exposition in the main text to a single patient. We provide additional results in \secref{s:additional-results-bandwidth}. The overall performance for these additional patients is consistent with the results reported here. We vary the scaling $\tau_s$ in~\eqref{e:bandwidthvel} for the RKHS representation of the bandwidth $\sigma_v$ of the Gaussian kernel in~\eqref{e:velrkhs} between 3 and 8. More precisely, we select values for $\tau_v$ in $\{3,4,6,8\}$. For values smaller than 3 and values larger than 8 we observed a deterioration of the performance of our solver. For each choice of $\tau_v$ we change the scaling $\tau_s$ in~\eqref{e:bandwidthdis} for the bandwidth $\sigma_s$ of the kernel distance in~\eqref{e:discretedist}. More precisely, we select values for $\tau_s$ in $\{\nicefrac{3}{4},1,2,4,6\}$. Likewise to $\tau_v$, we observed a deterioration of the performance of our methodology if we selected parameters outside the range $[\nicefrac{3}{4},6]$ for the considered data sets. Moreover, for $\tau_s = \nicefrac{1}{2}$ the estimate for the bandwidth $\sigma_s$ is below $\bar{h}_1 = \fnum{1.312801e+00}$.

\begin{figure}
\centering
\includegraphics[width=0.98\textwidth]{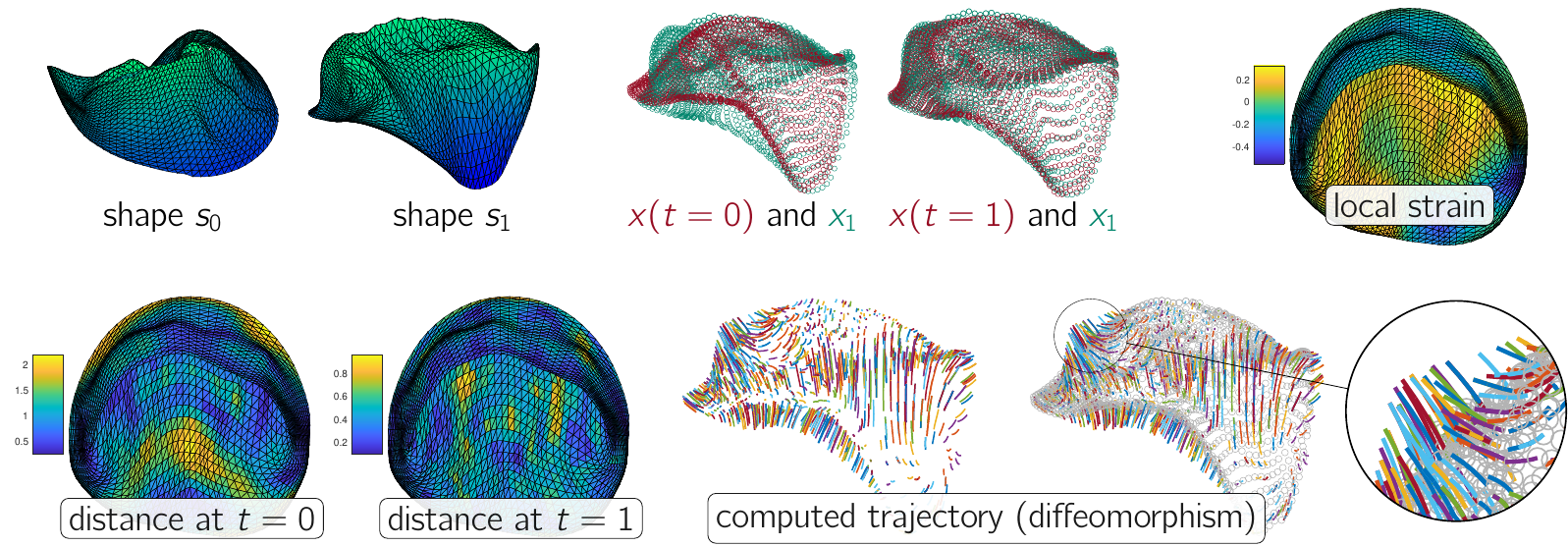}
\caption{Results for the matching of two mitral valve surfaces acquired from two different patients from our database. The results shown here correspond to the run with $\tau_s = 1$ and $\tau_v = 6$ reported in \tabref{t:xpat-BlJa-EgCl-probe-sigma-tauS-tauV} (row 12). The top row (from left to right) shows the original template shape $s_0$ and the target shape $s_1$ (Delaunay triangulations), the associated meshes in overlay before registration (template shape $x(t=0) = x_0$ in red and target shape $x_1$ in teal) and after registration (deformed template shape $x(t=1)$ in red and target shape $x_1$ in teal), and the computed strain intensity values overlaid on the target shape. The first two views in the top row show the entire mitral valve in their initial configuration (both leaflets). The first two views in the bottom row (point-wise distances from the template to the target shape) show the mitral valve from the top. For these two views one can identify the coaptation line that separates both leaflets; the coaptation line appears as a ridge between the anterior and posterior leaflet (see \figref{f:mitral-valve-regions} for an illustration). The bottom row shows (on the left) the point-wise distance from the template to the target shape before (at pseudo-time $t=0$) and after (at pseudo-time $t=1$) registration. The values are normalized with respect to the mean edge length of the target shape (i.e., values of 1.0 are equal to the average edge length). We also show the computed trajectory for 75\% of the points (randomly chosen; left: vector field only; right: vector field overlaid on the template shape $x(t=0) = x_0$). The circle shows a closeup of the computed vector field (zoomed view).\label{f:xpat-BlJa-EgCl-probe-sigma-tauS-1-tauV-6}}
\end{figure}

\paragraph{Results} We visualize results for the considered dataset in \figref{f:xpat-BlJa-EgCl-probe-sigma-tauS-1-tauV-6}. These results correspond to the parameter choices $\tau_v = 6$ and $\tau_s = 1$. In particular, we show the original template shape $s_0$ and the target shape $s_1$ in the top left corner (the surfaces are approximated from the meshes $x_0$ and $x_1$ through a Delaunay triangulation). In the middle block of the first row in \figref{f:xpat-BlJa-EgCl-probe-sigma-tauS-1-tauV-6} we show the template shape $x(t=0) = x_0$ and the reference shape before and after diffeomorphic shape matching. We can observe that the distance between the mesh points is initially large and becomes significantly smaller after diffeomorphic registration. The rightmost plot of the top row in \figref{f:xpat-BlJa-EgCl-probe-sigma-tauS-1-tauV-6} shows the associated isotropic strain intensity (see \secref{s:shapevariability}). We can observe for this particular patient pair that the strain is on average larger on the anterior leaflet than the posterior leaflet. The bottom row shows the point-wise distance of each mesh point from the template to the target shape before (left) and after (right)  diffeomorphic shape matching. Notice that the color map captures a different range of values before and after matching. On the right of the bottom row we illustrate the computed diffeomorphism (trajectory of the points $x_0$, i.e., the state variable $x(t)$) for 75\% of the points on the template shape selected at random. We provide two illustrations, one (left) that only shows the vector field and (right) the vector field together with the template shape $x_0 = x(t=0)$ in faint gray. We also show a closeup of the computed trajectory on the right.

Overall, the most important observations of this experiment are:
\begin{enumerate*}[label=(\roman*)]
\item Although the shapes are drastically different, we obtain a good matching quality, with pointwise distances below $\nicefrac{1}{2}$ of the target mesh resolution.
\item The computed diffeomorphism (trajectory) is quite complicated.
\end{enumerate*}

In \tabref{t:xpat-BlJa-EgCl-probe-sigma-tauS-tauV} we summarize the overall performance of our solver for the considered values for $\tau_v$ and $\tau_s$, respectively. We report the final (censored) Hausdorff distance (after diffeomorphic matching) in the third column from the left. We also report (in brackets) the percentage the obtained distance is with respect to the original Hausdorff distance of \fnum{2.981669e+00}. We also report the absolute and relative values of the norm of the primal and dual residuals in columns four and five. The last column shows the runtime for the solver in seconds (for 100 iterations).

The most important observations are:
\begin{enumerate*}[label=(\roman*)]
\item We can reduce the distance by up to 70\% for the choices $\tau_v = 4$ and $\tau_v = 6$ for $\tau_s \in \{\nicefrac{3}{4},1\}$, respectively.
\item The performance of our methodology in terms of registration quality (Hausdorff distance) as well as runtime is much more sensitive to the choice of $\tau_s$ compared to $\tau_v$.
\item The final norms of the residuals are only slightly affected by the choices for $\tau_v$  and $\tau_s$.
\end{enumerate*}

\begin{table}
\caption{Performance for different choices of the scaling parameters $\tau_v$ and $\tau_s$ that control the bandwidth $\sigma_v$ and $\sigma_s$ of the RKHS parameterization and the kernel distance, respectively. The results are for two representative patients from our database. The average mesh size $\bar{h}_1$ of the target shape is \fnum{1.312801e+00}. We select $\tau_v$ in $\{3,4,6,8\}$ and $\tau_s$ in $\{\nicefrac{3}{4},1,2,4,6\}$. The initial (censored) Hausdorff distance for this dataset is \fnum{2.981669e+00}. We report (from left to right) the final Hausdorff distance (after matching; absolute value and percentage of initial value in brackets), the (relative) norm of the primal residual, the (relative) norm of the dual residual, and the runtime of the solver (in seconds).\label{t:xpat-BlJa-EgCl-probe-sigma-tauS-tauV}}
\centering\small
\begin{tabular}{lllllr}\toprule
$\tau_v$ & $\tau_s$   & {\bf final distance} (\%) & \multicolumn{2}{l}{\bf residuals}   & {\bf runtime} \\
         &            &                     & primal (relative) & dual (relative) \\\midrule
3 & 0.75 & \fnum{1.161874e+00} (38.97) & \snum{4.492032e-03} (\snum{3.596596e-04})  & \snum{5.392893e-01} (\snum{4.317881e-02}) & \fnum{100.387234} \\
  & 1    & \fnum{1.173727e+00} (39.36) & \snum{2.906727e-03} (\snum{2.120993e-04})  & \snum{4.780634e-01} (\snum{3.488352e-02}) & \fnum{83.386913} \\
  & 2    & \fnum{1.450821e+00} (48.66) & \snum{2.089185e-03} (\snum{1.495835e-04})  & \snum{3.810458e-01} (\snum{2.728250e-02}) & \fnum{68.938882} \\
  & 4    & \fnum{2.107427e+00} (70.68) & \snum{2.060969e-03} (\snum{2.088280e-04})  & \snum{2.977642e-01} (\snum{3.017100e-02}) & \fnum{63.820017} \\
  & 6    & \fnum{2.534500e+00} (85.00) & \snum{1.868288e-03} (\snum{2.839526e-04})  & \snum{2.767971e-01} (\snum{4.206914e-02}) & \fnum{62.391137} \\
\midrule
4 & 0.75 & \fnum{9.058977e-01} (30.38) & \snum{5.894792e-03} (\snum{4.719732e-04})  & \snum{5.885302e-01} (\snum{4.712133e-02}) & \fnum{91.347397} \\
  & 1    & \fnum{9.356657e-01} (31.38) & \snum{3.987376e-03} (\snum{2.909525e-04})  & \snum{5.234655e-01} (\snum{3.819644e-02}) & \fnum{82.905749} \\
  & 2    & \fnum{1.211718e+00} (40.64) & \snum{2.781515e-03} (\snum{1.991537e-04})  & \snum{4.206576e-01} (\snum{3.011866e-02}) & \fnum{69.267063} \\
  & 4    & \fnum{1.879761e+00} (63.04) & \snum{2.841239e-03} (\snum{2.878889e-04})  & \snum{3.506940e-01} (\snum{3.553412e-02}) & \fnum{61.823046} \\
  & 6    & \fnum{2.339422e+00} (78.46) & \snum{2.483905e-03} (\snum{3.775175e-04})  & \snum{3.437944e-01} (\snum{5.225176e-02}) & \fnum{60.726389} \\
\midrule
6 & 0.75 & \fnum{9.073484e-01} (30.43) & \snum{7.560293e-03} (\snum{6.053234e-04})  & \snum{6.785546e-01} (\snum{5.432924e-02}) & \fnum{92.390780} \\
  & 1    & \fnum{9.242514e-01} (31.00) & \snum{4.982128e-03} (\snum{3.635379e-04})  & \snum{6.145162e-01} (\snum{4.484027e-02}) & \fnum{78.333367} \\
  & 2    & \fnum{1.143209e+00} (38.34) & \snum{3.259263e-03} (\snum{2.333600e-04})  & \snum{4.899762e-01} (\snum{3.508181e-02}) & \fnum{78.800609} \\
  & 4    & \fnum{1.697406e+00} (56.93) & \snum{3.184282e-03} (\snum{3.226478e-04})  & \snum{4.063452e-01} (\snum{4.117298e-02}) & \fnum{65.628645} \\
  & 6    & \fnum{2.155416e+00} (72.29) & \snum{2.725908e-03} (\snum{4.142984e-04})  & \snum{4.092239e-01} (\snum{6.219609e-02}) & \fnum{64.217187} \\
\midrule
8 & 0.75 & \fnum{9.896912e-01} (33.19) & \snum{8.629760e-03} (\snum{6.909515e-04})  & \snum{7.472190e-01} (\snum{5.982693e-02}) & \fnum{94.089245} \\
  & 1    & \fnum{1.009466e+00} (33.86) & \snum{6.412356e-03} (\snum{4.678994e-04})  & \snum{6.709854e-01} (\snum{4.896074e-02}) & \fnum{86.053542} \\
  & 2    & \fnum{1.237214e+00} (41.49) & \snum{3.535220e-03} (\snum{2.531182e-04})  & \snum{5.214616e-01} (\snum{3.733613e-02}) & \fnum{75.360986} \\
  & 4    & \fnum{1.691207e+00} (56.72) & \snum{2.975507e-03} (\snum{3.014936e-04})  & \snum{4.036373e-01} (\snum{4.089860e-02}) & \fnum{64.831449} \\
  & 6    & \fnum{2.100938e+00} (70.46) & \snum{2.612026e-03} (\snum{3.969900e-04})  & \snum{4.032289e-01} (\snum{6.128493e-02}) & \fnum{59.485455} \\
\bottomrule
\end{tabular}
\end{table}

\begin{figure}
\centering
\includegraphics[width=0.98\textwidth]{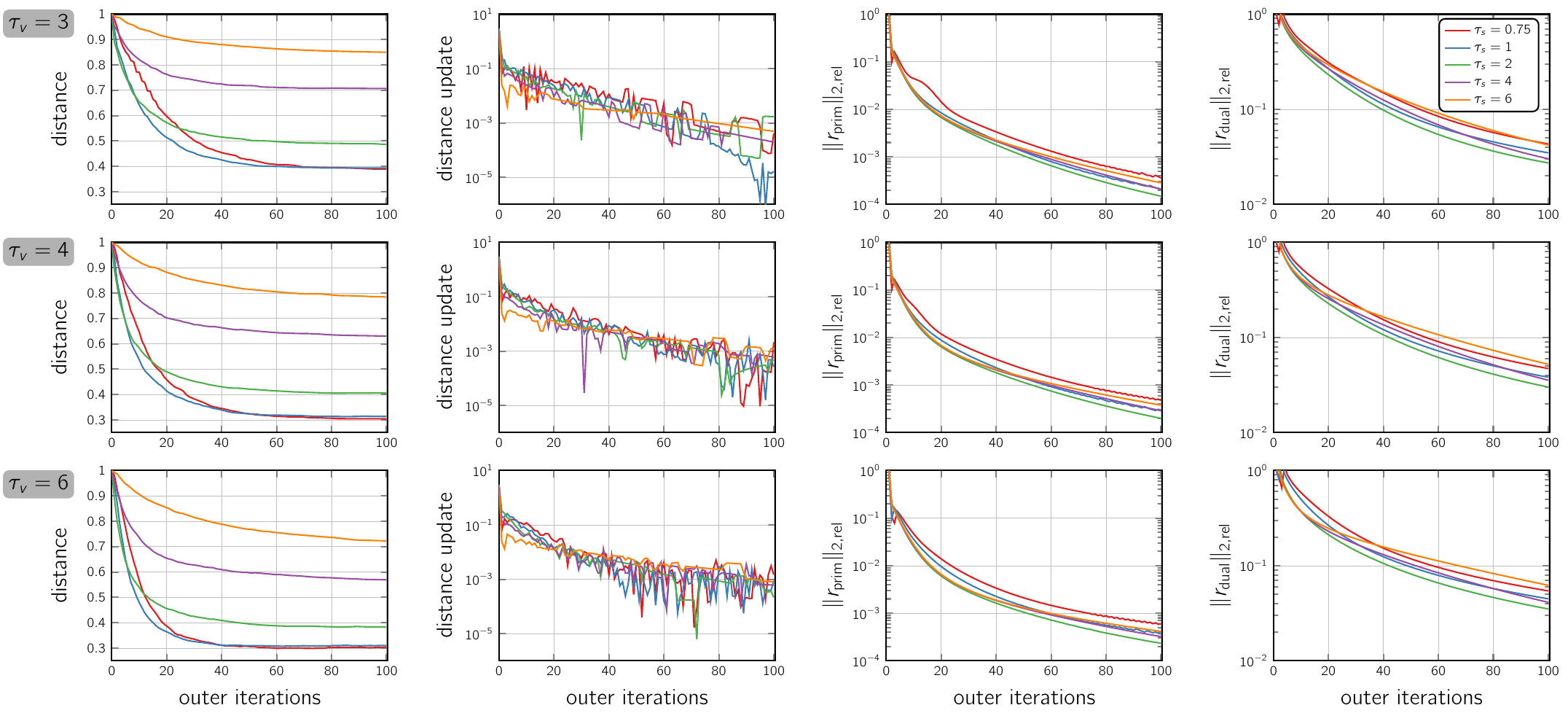}
\caption{Performance for different choices of the scaling parameters $\tau_v$ and $\tau_s$ that control the bandwidth $\sigma_v$ and $\sigma_s$ of the RKHS parameterization and the kernel distance, respectively. The results are for two representative patients from our database. The average mesh size $\bar{h}_1$ of the target shape is \fnum{1.312801e+00}. The results correspond to the first, second and third block in \tabref{t:xpat-BlJa-EgCl-probe-sigma-tauS-tauV}. From top to bottom (row 1 through 3) we report results for $\tau_v=3$, $\tau_v=4$ and $\tau_v=6$. Each plot shows results for $\tau_s \in \{\nicefrac{3}{4},1,2,4,6\}$ (see legend). The first column shows the trend of the Hausdorff distance with respect to the number of iterations (normalized with respect to the initial value). Notice that we do not minimize the Hausdorff distance in our optimization but the kernel distance in \eqref{e:discretedist}. The second column shows the update of the Hausdorff distance per iteration. The third and fourth column shows the norm of the primal and dual residual versus the iteration count. These residuals are normalized with respect to the residual obtained at iteration 2.\label{f:xpat-BlJa-EgCl-probe-sigma-tauS-tauV}}
\end{figure}

In \figref{f:xpat-BlJa-EgCl-probe-sigma-tauS-tauV} we illustrate the convergence behavior of our solver. The horizontal axis is the iteration count of the splitting method (outer iteration index $k$) for all plots. We show (from left to right)
\begin{enumerate*}[label=(\roman*)]
\item the trend of the Hausdorff distance (normalized by the initial distance at iteration 1),
\item the trend of the updates of the Hausdorff distance,
\item the trend of the relative primal residual, and
\item the trend of the relative dual residual.
\end{enumerate*}

The most important observations are the following:
\begin{enumerate*}[label=(\roman*)]
\item The smaller $\tau_s$ the more accurate the matching. The results for $\tau_s = \nicefrac{3}{4}$ and $\tau_s = 1$ are quite similar.
\item Increasing $\tau_v$ improves the registration accuracy; the optimal accuracy is obtained for $\tau_v = 6$. (In \tabref{t:xpat-BlJa-EgCl-probe-sigma-tauS-tauV} we saw that for $\tau_v = 8$, the accuracy decreases again.)
\item The updates for the Hausdorff distance decrease with increasing iterates.
\item The norms of the residuals decrease (almost) monotonically across all choices of parameters $\tau_s$ and $\tau_v$.
\item The trend of the residual norms is only slightly affected by the choices of the parameters $\tau_s$ and $\tau_v$, respectively.
\end{enumerate*}

\subsection{Multi-Frame Registration}

We report results for the registration across time of successive mitral valve surface snapshots acquired between mid-systole and end-systole for one individual patient. That is, we are given data at different points in time and try to compute the diffeomorphic flow from the first frame at mid systole through all data to the last frame at end systole. All available frames enter the kernel distance (i.e., we sum over all frames). This setup is equivalent to the one considered in \cite{Azencott:2010a,Zhang:2021a}. The real time duration of this typical 25 frames sequence is less than half a second.

\begin{remark}
We note that the mitral valve changes topology as we go from mid-systole to end-systole. That is, the mitral valve opens. This constitutes a significant problem when trying to compute a diffeomorphism that maps corresponding points to one another and preserves the integrity of the considered anatomy. Based on our parameterization of the curved surface, we compute a diffeomorphism that takes us from one set of points to the other. However, it is not guaranteed that we do not match points situated on one leaflet to those on another. One simple fix to this problem is to register the leaflets of the mitral valves individually. Our data has the anatomical annotations to do so (see \figref{f:mitral-valve-regions}). In fact, our software supports this, and we have tested this strategy to classify patient data~\cite{Dabirian:2022a}. In the present work, we are solely interested in studying the performance of the proposed solver. As such, we consider the entire mitral valve data and do not register the leaflets individually. Since the boundaries of individual leaflets are well identified even when the mitral valve is closed, one can mitigate the initial topology change by using high weights for leaflet boundary points in the surface matching terms (see \secref{s:discdist} for a discussion).

On a more general note, using diffeomorphic registration in problems where the topology changes and/or structures are not present in both datasets to be registered is a delicate problem, and an active area of research. We refer to~\cite{Scheufele:2019a,Scheufele:2020a,Zacharaki:2009a,Gooya:2012a,Hogea:2008a,Zacharaki:2008a,Li:2012a} for examples for medical images and~\cite{Hsieh:2022b,Franccois:2022a,Antonsanti:2021a,Sukurdeep:2022a} for work on shape representations that attempt to deal with this issue.
\end{remark}

\paragraph{Setup} We use the parameter values for our surface matching solver as suggested by our results in \tabref{t:parameters}.

More precisely,  we set the weight for the kernel distance is set to $\alpha = 1$. The scaling parameters $\tau_s$ and $\tau_v$ are set to 1 and 6, respectively (based on the observations from our prior numerical experiments reported in \secref{s:bandwithselection}). The parameter $\rho$ is set to 1. The tolerances for the primal and dual residual are $\epsilon_{\text{prim}} = \num{1e-3}$ and $\epsilon_{\text{dual}} = \num{1e-3}$, respectively. The scaling for the tolerance of the censored Hausdorff distance is $\tau_{\text{haus}} =\nicefrac{1}{2}$. We limit the maximum number of iterations to 100 (not actually reached).

The main difference compared to the experiments reported in the former section is that we now have several 3D-image data available to control our diffeomorphic registration results: multiple frames that monitor the heart cycle of an individual patients. We compute the diffeomorphism that flows the first shape at time $t=0$ (mid-systole) through all these data until we reach the final shape at $t=1$ (end-systole). See \cite{Zhang:2021a} for additional details.

\begin{figure}
\centering
\includegraphics[width=0.98\textwidth]{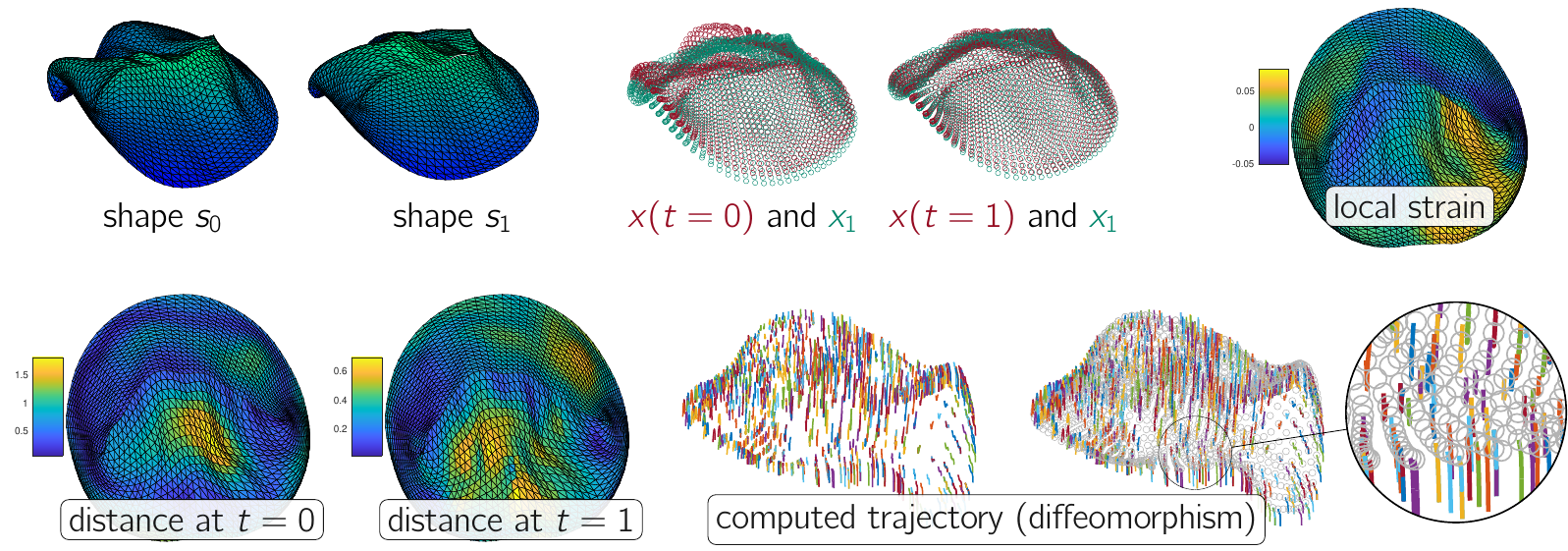}
\caption{Results for  diffeomorphic matching between the  two mitral valve surface snapshots acquired from a single patient at mid systole and end systole (first patient in \tabref{t:results-MF-tauS-1-tauV-6}). We set $\tau_s = 1$ and $\tau_v = 6$. The top row (from left to right) shows the template shape $s_0$ (mitral valve at mid-systole; Delaunay triangulation) and the target shape $s_1$ (mitral valve at end-systole; Delaunay triangulation), the associated meshes in overlay before registration (template shape $x(t=0) = x_0$ in red and target shape $x_1$ in teal) and after registration (deformed template shape $x(t=1)$ in red and target shape $x_1$ in teal), and the computed strain intensity values overlaid on the target shape. The first two views in the top row show the entire mitral valve in their initial configuration (both leaflets). The first two views in the bottom row (point-wise distances from the template to the target shape) show the mitral valve from the top. For these two views  one can identify the coaptation line that separates both leaflets; the coaptation line appears as a ridge between the anterior and posterior leaflet (see \figref{f:mitral-valve-regions} for an illustration). The bottom row shows (on the left) the point-wise distance from the template to the target shape before (at pseudo-time $t=0$) and after (at pseudo-time $t=1$) registration. The distance values are normalized with respect to the mean edge length of the target shape (i.e., values of 1.0 are equal to the average edge length). We also show the computed trajectory for 75\% of the points (randomly chosen; left: vector field only; right: vector field overlaid on the template shape $x(t=0) = x_0$). The circle shows a closeup of the computed vector field (zoomed view).\label{f:results-viz-AdVe-MF-tauS-1-tauV-6}}
\end{figure}

\paragraph{Results} We show a representative registration result in \figref{f:results-viz-AdVe-MF-tauS-1-tauV-6}. The setup of this figure is exactly the same as for the results reported in \figref{f:xpat-BlJa-EgCl-probe-sigma-tauS-tauV}. We refer to \secref{s:bandwithselection} for an explanation of the layout of this figure. The observations are similar to those we have made in \secref{s:bandwithselection}.

In addition, we report performance results for 20 patients from our database in \tabref{t:results-MF-tauS-1-tauV-6}. For each patient, we report the following quantities of interest (from left to right; all times are reported in seconds):
\begin{enumerate*}[label=(\roman*)]
\item The number of outer iterations of our splitting algorithm.
\item The Hausdorff distance after diffeomorphic shape matching and the associated percentage of the original Hausdorff distance (in brackets).
\item The time-to-solution (TTS; runtime of the solver).
\item The overall time spent on solving the distance subproblem and the percentage this amounts to with respect to the time-to-solution (in brackets).
\item The overall time spent on computing the Newton steps during the solution of the distance subproblem. We also provide (in brackets) the percentage this computation amounts to with respect to the time needed to solve the distance subproblem.
\item The time spent to setup the KKT system for the kinetic energy subproblem.
\item The overall time spent on solving the kinetic energy system and the percentage this amounts to with respect to the time-to-solution (in brackets).
\item The overall time spent on solving the KKT system and the percentage this amounts to with respect to the solution of the distance subproblem (in brackets).
\end{enumerate*}

\begin{table}
\caption{Results for  diffeomorphic matching between mid systole and end systole mitral valve surface snapshots for 20 representative patients in our database. All run times are reported in seconds. We report (first four columns from left to right): ($i$) the run id (different patients in each row), ($ii$) the number of outer iterations of our splitting algorithm, ($iii$) the final censored Hausdorff distance after diffeomorphic matching (\% of original distance in brackets), and ($iv$) the overall time to solution (TTS). In the fifth and the sixth column we report run times for the distance subproblem. In particular, we report ($v$) the runtime spent on solving the distance subproblem (percentage of TTS in brackets) and ($vi$) the time spent on the solving for the search direction (Newton step; percentage of time for overall solution of distance subproblem in brackets). The seventh, eighths and ninths column report run times for the kinetic energy subproblem. In particular, we report ($vii$) time spent on the setup of the KKT system, ($viii$) time spent on solving the kinetic energy subproblem (percentage of TTS in brackets), and ($ix$) time spent on solving the KKT system (percentage of time for overall solution of kinetic energy subproblem in brackets).
\label{t:results-MF-tauS-1-tauV-6}
}
\small\centering\resetrid
\begin{tabular}{rrrrrrrrrr}\toprule
\multicolumn{4}{c}{} & \multicolumn{2}{l}{\bf distance subproblem} & \multicolumn{3}{l}{\bf kinetic energy subproblem}\\
\bf run & \bf iter & \bf distance & \bf TTS & \bf solve (\%) & \bf newton  (\%) & \bf setup & \bf solve (\%) & \bf solve KKT  (\%) \\\midrule
\runid & 25 & \fnum{6.5549e-01}~(\inum{42.7209}) & \fnum{ 87.6215} & \fnum{ 50.9406}~(\inum{58.1371}) & \fnum{ 50.9216}~(\inum{99.9626}) & \fnum{1.8002} & \fnum{ 36.5291}~(\inum{41.6897}) & \fnum{ 34.8083}~(\inum{95.2890}) \\
\runid & 12 & \fnum{6.3409e-01}~(\inum{44.0736}) & \fnum{ 51.9937} & \fnum{ 32.5011}~(\inum{62.5098}) & \fnum{ 32.4921}~(\inum{99.9721}) & \fnum{1.6761} & \fnum{ 19.4133}~(\inum{37.3378}) & \fnum{ 18.7097}~(\inum{96.3759}) \\
\runid & 70 & \fnum{7.0181e-01}~(\inum{59.7196}) & \fnum{153.6740} & \fnum{ 92.3071}~(\inum{60.0668}) & \fnum{ 92.2621}~(\inum{99.9512}) & \fnum{1.5030} & \fnum{ 60.9607}~(\inum{39.6688}) & \fnum{ 57.4154}~(\inum{94.1843}) \\
\runid & 49 & \fnum{6.8984e-01}~(\inum{30.8863}) & \fnum{158.6951} & \fnum{109.8692}~(\inum{69.2329}) & \fnum{109.8377}~(\inum{99.9713}) & \fnum{1.4983} & \fnum{ 48.5415}~(\inum{30.5879}) & \fnum{ 46.0214}~(\inum{94.8085}) \\
\runid & 46 & \fnum{6.6010e-01}~(\inum{34.4157}) & \fnum{136.7009} & \fnum{ 81.6284}~(\inum{59.7132}) & \fnum{ 81.5950}~(\inum{99.9590}) & \fnum{1.7067} & \fnum{ 54.7981}~(\inum{40.0861}) & \fnum{ 51.9936}~(\inum{94.8821}) \\
\runid &  5 & \fnum{4.9879e-01}~(\inum{60.0386}) & \fnum{ 19.6431} & \fnum{ 12.4055}~(\inum{63.1548}) & \fnum{ 12.4024}~(\inum{99.9746}) & \fnum{1.4048} & \fnum{  7.1981}~(\inum{36.6443}) & \fnum{  6.9123}~(\inum{96.0302}) \\
\runid &  5 & \fnum{5.4241e-01}~(\inum{43.4489}) & \fnum{ 24.9475} & \fnum{ 14.8915}~(\inum{59.6914}) & \fnum{ 14.8878}~(\inum{99.9751}) & \fnum{1.7091} & \fnum{ 10.0163}~(\inum{40.1494}) & \fnum{  9.6671}~(\inum{96.5142}) \\
\runid & 23 & \fnum{4.4807e-01}~(\inum{39.7045}) & \fnum{114.3949} & \fnum{ 86.5237}~(\inum{75.6360}) & \fnum{ 86.5089}~(\inum{99.9829}) & \fnum{1.4983} & \fnum{ 27.7267}~(\inum{24.2377}) & \fnum{ 26.0083}~(\inum{93.8023}) \\
\runid & 13 & \fnum{5.5921e-01}~(\inum{42.4099}) & \fnum{ 47.8315} & \fnum{ 31.3558}~(\inum{65.5547}) & \fnum{ 31.3475}~(\inum{99.9736}) & \fnum{1.4411} & \fnum{ 16.3905}~(\inum{34.2672}) & \fnum{ 15.7146}~(\inum{95.8764}) \\
\runid & 45 & \fnum{6.2089e-01}~(\inum{48.2060}) & \fnum{ 81.1076} & \fnum{ 51.1301}~(\inum{63.0398}) & \fnum{ 51.1049}~(\inum{99.9508}) & \fnum{1.1906} & \fnum{ 29.7204}~(\inum{36.6432}) & \fnum{ 27.7612}~(\inum{93.4079}) \\
\runid &  6 & \fnum{5.6235e-01}~(\inum{46.2197}) & \fnum{ 26.5705} & \fnum{ 18.0384}~(\inum{67.8887}) & \fnum{ 18.0345}~(\inum{99.9786}) & \fnum{1.4159} & \fnum{  8.4869}~(\inum{31.9409}) & \fnum{  8.1683}~(\inum{96.2463}) \\
\runid & 67 & \fnum{6.6431e-01}~(\inum{31.3638}) & \fnum{263.9978} & \fnum{157.4552}~(\inum{59.6426}) & \fnum{157.4002}~(\inum{99.9650}) & \fnum{1.9665} & \fnum{106.1536}~(\inum{40.2100}) & \fnum{101.3590}~(\inum{95.4833}) \\
\runid & 11 & \fnum{5.2352e-01}~(\inum{46.5362}) & \fnum{ 50.6636} & \fnum{ 31.6562}~(\inum{62.4832}) & \fnum{ 31.6483}~(\inum{99.9749}) & \fnum{1.7606} & \fnum{ 18.9345}~(\inum{37.3729}) & \fnum{ 18.2714}~(\inum{96.4984}) \\
\runid & 14 & \fnum{5.3117e-01}~(\inum{41.1869}) & \fnum{ 67.5811} & \fnum{ 41.0940}~(\inum{60.8069}) & \fnum{ 41.0835}~(\inum{99.9745}) & \fnum{1.7943} & \fnum{ 26.3958}~(\inum{39.0580}) & \fnum{ 25.4460}~(\inum{96.4016}) \\
\runid & 33 & \fnum{6.7477e-01}~(\inum{44.0021}) & \fnum{174.0991} & \fnum{103.0771}~(\inum{59.2060}) & \fnum{103.0470}~(\inum{99.9708}) & \fnum{2.2225} & \fnum{ 70.8218}~(\inum{40.6790}) & \fnum{ 67.6815}~(\inum{95.5658}) \\
\runid & 51 & \fnum{7.2125e-01}~(\inum{39.8295}) & \fnum{230.7465} & \fnum{129.2283}~(\inum{56.0045}) & \fnum{129.1826}~(\inum{99.9646}) & \fnum{2.4128} & \fnum{101.2076}~(\inum{43.8609}) & \fnum{ 96.8654}~(\inum{95.7097}) \\
\runid &  9 & \fnum{4.9534e-01}~(\inum{38.7374}) & \fnum{ 58.6898} & \fnum{ 41.4656}~(\inum{70.6521}) & \fnum{ 41.4552}~(\inum{99.9750}) & \fnum{1.7270} & \fnum{ 17.1581}~(\inum{29.2351}) & \fnum{ 16.6072}~(\inum{96.7897}) \\
\runid & 58 & \fnum{5.5777e-01}~(\inum{56.1411}) & \fnum{204.4451} & \fnum{115.2271}~(\inum{56.3609}) & \fnum{115.1813}~(\inum{99.9602}) & \fnum{1.9478} & \fnum{ 88.8809}~(\inum{43.4742}) & \fnum{ 84.6596}~(\inum{95.2506}) \\
\runid & 41 & \fnum{6.7397e-01}~(\inum{35.8763}) & \fnum{176.5596} & \fnum{101.9649}~(\inum{57.7510}) & \fnum{101.9327}~(\inum{99.9685}) & \fnum{2.1086} & \fnum{ 74.3524}~(\inum{42.1118}) & \fnum{ 71.3570}~(\inum{95.9714}) \\
\runid & 24 & \fnum{5.2741e-01}~(\inum{58.8882}) & \fnum{ 75.1826} & \fnum{ 47.5213}~(\inum{63.2079}) & \fnum{ 47.5068}~(\inum{99.9694}) & \fnum{1.4970} & \fnum{ 27.5179}~(\inum{36.6015}) & \fnum{ 26.2066}~(\inum{95.2348}) \\
\bottomrule
\end{tabular}
\end{table}

The most important observations are:
\begin{enumerate*}[label=(\roman*)]
\item We can see that we can reduce the distance by up to 60\% compared to the original distance between the shape at $t=0$ and $t=1$. This gain is smaller than for the multi-patient case, but we also (generally) start with a smaller initial distance between the surfaces to be matched. This is also reflected by the smaller Hausdorff distance after diffeomorphic matching compared to the best results reported in \tabref{t:xpat-BlJa-EgCl-probe-sigma-tauS-tauV}.
\item The time until convergence ranges from about 25 seconds up to about 250 seconds (depending on the dataset).
\item On average, we spend about 60\% of the total runtime on solving the distance subproblem and about 40\% of the runtime on solving the kinetic energy subproblem.
\item The runtime in each subproblem is dominated by solving the optimality conditions/computing the search direction; the percentages are on average at the order of almost 100\% for the distance subproblem and about 95\% for the kinetic energy subproblem.
\end{enumerate*}

\section{Conclusions}\label{s:conclusions}

We have presented and studied an improved matrix-free implementation of an operator splitting optimization algorithm for the solution of variational optimal control formulations for diffeomorphic shape matching. As opposed to many other works in the field, we consider a discretize-then-optimize approach to solve this problem. The associated infinite-dimensional variational problem is naturally recast as a finite-dimensional optimization problem in a RKHS. Our algorithmic improvements yield a speedup at the order of 2$\times$ compared to our prior work in~\cite{Zhang:2021a}. Moreover, we have reduced the memory pressure since we no longer explicitly form and store matrix operators as they appear in the optimality conditions of the individual subproblems. Aside from introducing these numerical developments, we have reported new results for the registration of smooth 3D surfaces acquired from different subjects (multi-subject registration) and for time series of smooth 3D surfaces acquired from individual patients (multi-frame registration). In addition, we have surveyed our prior work in this area, which spans a period of more than one decade, with work touching on the mathematical problem formulation~\cite{Azencott:2008a,Azencott:2010a}, algorithmic developments~\cite{Azencott:2010a,Freeman:2014a,Jajoo:2011a,Qin:2013a,Zhang:2019a,Zhang:2021a}, clinical applications~\cite{Zekry:2018a,ElTallawi:2021a,ElTallawi:2021b,ElTallawi:2019a}, as well as machine learning approaches for classification in shape space~\cite{Dabirian:2022a}.

Our current and future work is targeted towards developing computational approaches for classification in shape space as well as improving computational throughput for algorithms for diffeomorphic shape matching and image registration. For the work described in the present paper we have identified several future research directions. An open problem is the tuning of hyper-parameters. The search for an adequate weight $\alpha$ for the shape distance can be performed based on a parameter continuation scheme~\cite{Azencott:2010a,Mang:2015a,Mang:2019a}. However, an open question is how to most effectively combine this search with a multi-scale policy for selecting the parameters $\sigma_v$ and $\sigma_s$, respectively. From a numerical perspective, several components of our solver can be improved. The most critical component is the implementation of the kernel matrices $K$. One possibility for improvements is to explore other, more efficient kernel representations~\cite{Jain:2013a}. An efficient implementation for the associated kernel operations is readily available in Python, and described in~\cite{Charlier:2021a}. Moreover, although we have significantly reduced the memory footprint of our solver in the implementation described in this manuscript, we still store the kernel matrices and evaluate them using tensor algebra. Implementing more efficient kernel evaluations that are entirely matrix-free remains subject to future work. Lastly, our implementation is not deployed on any modern computing architectures. Porting our software to graphics-processing units may significantly improve the computational throughput~\cite{Brunn:2021a,Brunn:2020a}.

\textbf{Acknowledgements:} This work was partly supported by the National Science Foundation ({\bf NSF}) under the awards DMS-1854853, DMS-2012825, and DMS-2145845. Any opinions, findings, and conclusions or recommendations expressed herein are those of the authors and do not necessarily reflect the views of NSF. This work was completed in part with resources provided by the Research Computing Data Core at the University of Houston. We thank our collaborators Dr. William Zoghbi (MD) and  Dr. Carlos El Tallawi (MD) , from the DeBakey Heart and Vascular Center of the Houston Methodist Hospital for sharing and curating the annotated database of 3D echocardiology image sequences of the mitral valve acquired from a cohort of mitral valve patients according to local clinical protocols.

Our thanks also go out to our colleague, Dr. James Herring, and former PhD students, Drs. Jeff Friedman, Saurhab Jain, Aarti Jajoo, Yue Qin, and Peng Zhang, listed in alphabetical order. Their significant contributions over the course of this research project have been invaluable. We would like to convey our special gratitude to our esteemed late colleagues, Dr. Roland Glowinski and Dr. Ronald Hoppe. Their profound insights and advice, stemming from numerous discussions over the past two decades, have left an indelible mark on our work.

\begin{appendix}
\section{Hardware}\label{s:hardware}

Standard desktop computer: Mac Studio (Mac13,1) equipped with an Apple M1 Max 10-core CPU and 32 GB LPDDR5 memory. The code is implemented in Matlab (R2022b, Release 9.13.0.2105380).

\section{Additional Numerical Results}\label{s:additional-results-bandwidth}

\begin{figure}
\centering
\includegraphics[width=0.98\textwidth]{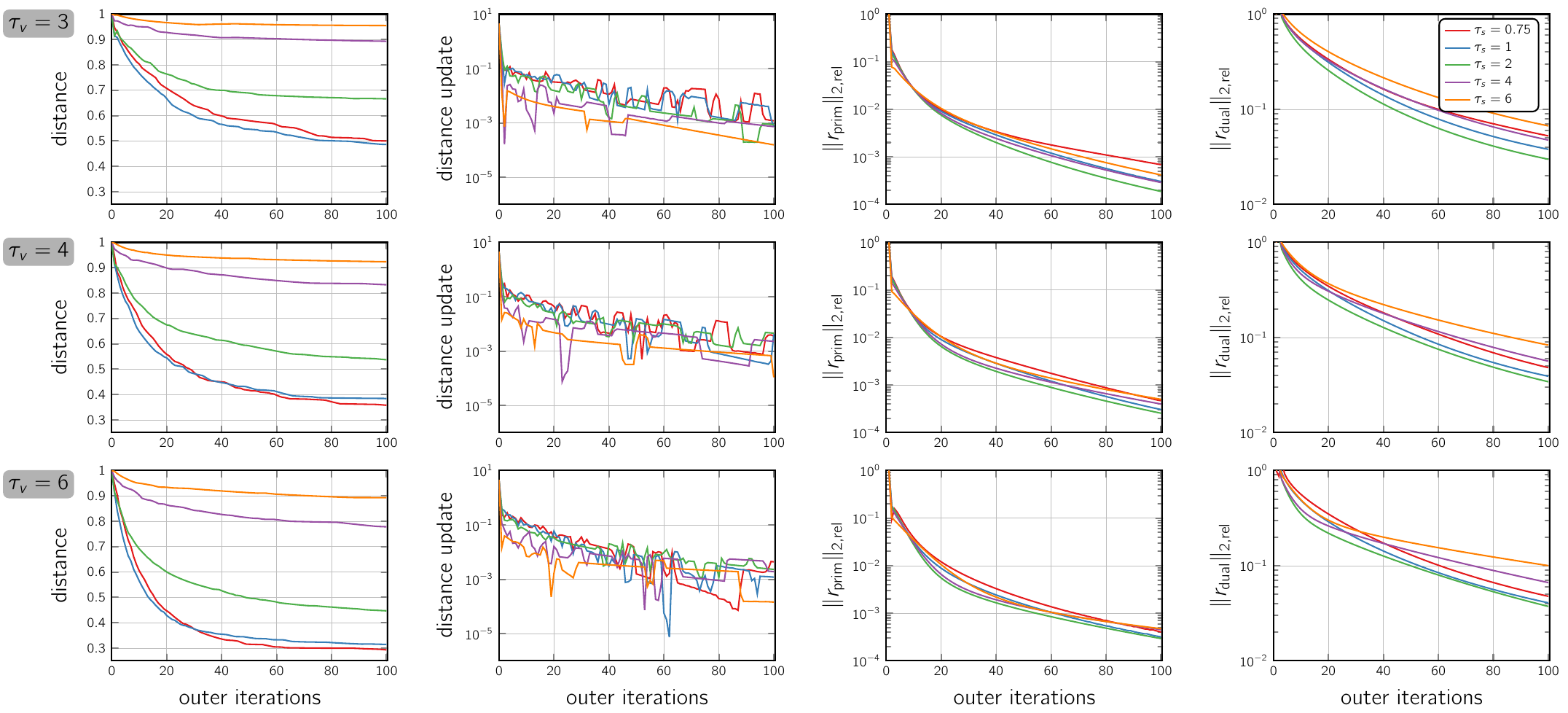}
\caption{Performance for different choices of the scaling parameters $\tau_v$ and $\tau_s$ that control the bandwidth $\sigma_v$ and $\sigma_s$ of the RKHS parameterization and the kernel distance, respectively. The results are for two representative patients from our database. The average mesh size $\bar{h}_1$ of the target shape is \fnum{1.239477e+00}. The results correspond to the first, second and third block in \tabref{t:xpat-ShCr-ShCh-probe-sigma-tauS-tauV}. From top to bottom (row 1 through 3) we report results for $\tau_v=3$, $\tau_v=4$ and $\tau_v=6$. Each plot shows results for $\tau_s \in \{\nicefrac{3}{4},1,2,4,6\}$ (see legend). The first column shows the trend of the Hausdorff distance with respect to the number of iterations (normalized with respect to the initial value). Notice that we do not minimize the Hausdorff distance in our optimization but the kernel distance in \eqref{e:discretedist}. The second column shows the update of the Hausdorff distance per iteration. The third and fourth column shows the norm of the primal and dual residual versus the iteration count. These residuals are normalized with respect to the residual obtained at iteration 2.\label{f:xpat-ShCr-ShCh-probe-sigma-tauS-tauV}}
\end{figure}

\begin{figure}
\centering
\includegraphics[width=0.98\textwidth]{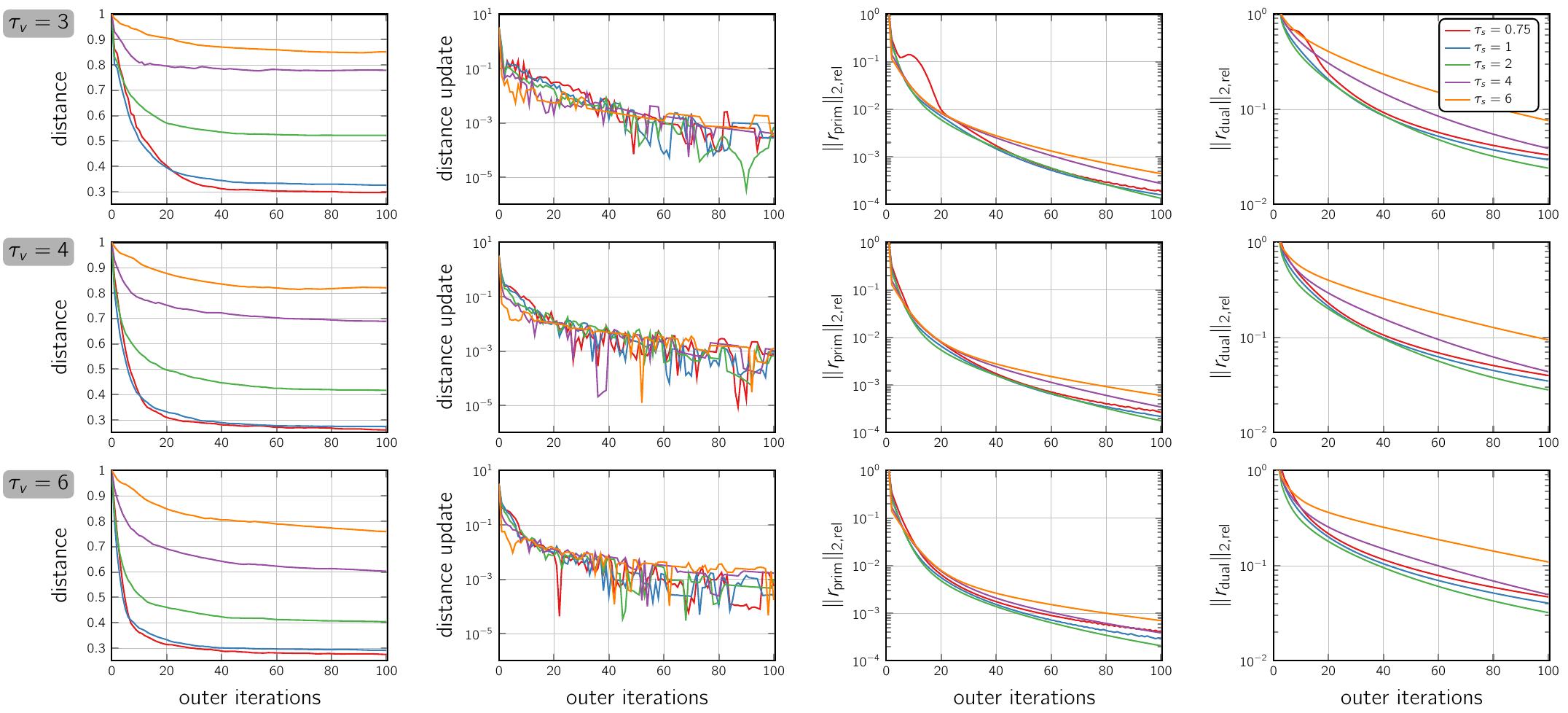}
\caption{Performance for different choices of the scaling parameters $\tau_v$ and $\tau_s$ that control the bandwidth $\sigma_v$ and $\sigma_s$ of the RKHS parameterization and the kernel distance, respectively. The results are for two representative patients from our database. The average mesh size $\bar{h}_1$ of the target shape is \fnum{1.076166e+00}. The results correspond to the first, second and third block in \tabref{t:xpat-WaSh-YaSa-probe-sigma-tauS-tauV}. From top to bottom (row 1 through 3) we report results for $\tau_v=3$, $\tau_v=4$ and $\tau_v=6$. Each plot shows results for $\tau_s \in \{\nicefrac{1}{4},1,2,4,6\}$ (see legend). The first column shows the trend of the Hausdorff distance with respect to the number of iterations (normalized with respect to the initial value). Notice that we do not minimize the Hausdorff distance in our optimization but the kernel distance in \eqref{e:discretedist}. The second column shows the update of the Hausdorff distance per iteration. The third and fourth column shows the norm of the primal and dual residual versus the iteration count. These residuals are normalized with respect to the residual obtained at iteration 2.\label{f:xpat-WaSh-YaSa-probe-sigma-tauS-tauV}}
\end{figure}

\begin{figure}
\centering
\includegraphics[width=0.98\textwidth]{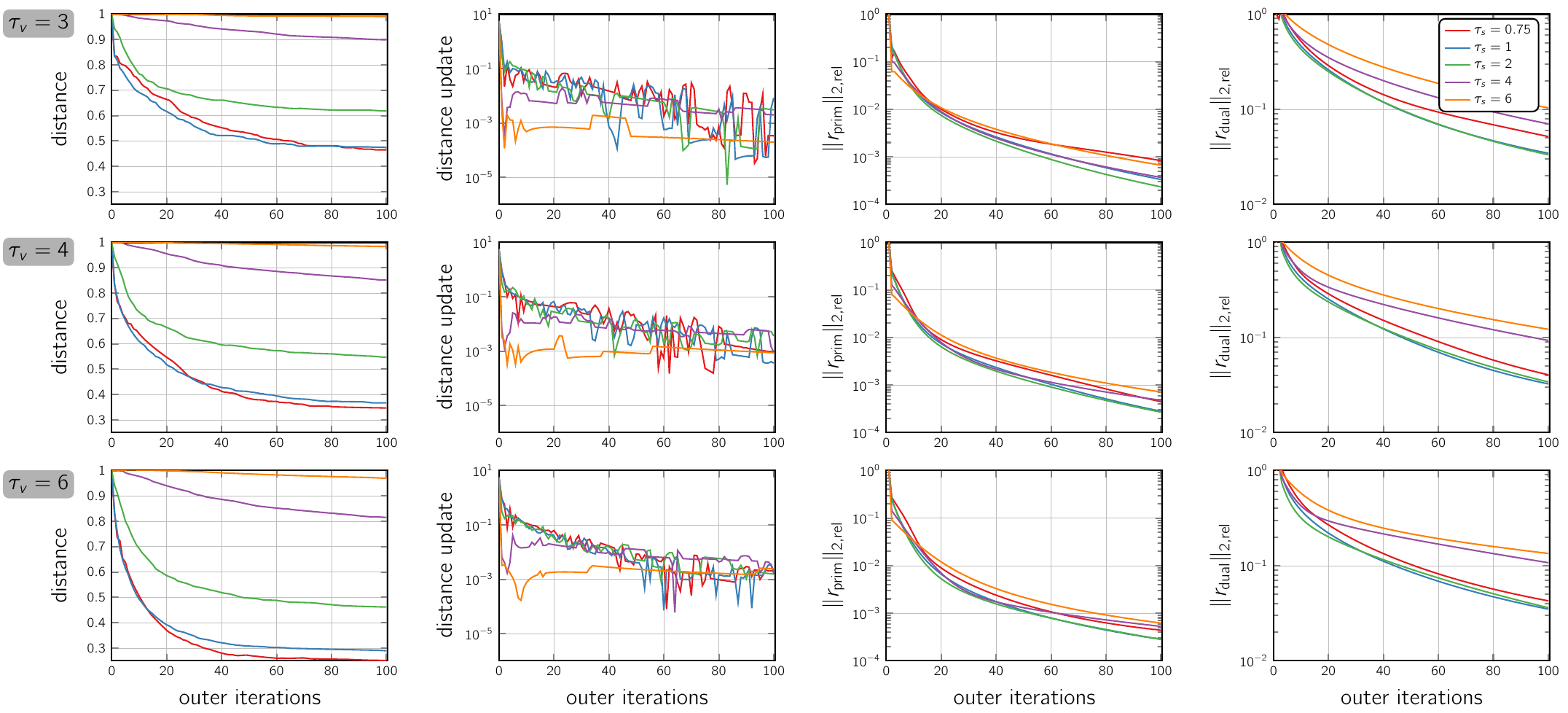}
\caption{Performance for different choices of the scaling parameters $\tau_v$ and $\tau_s$ that control the bandwidth $\sigma_v$ and $\sigma_s$ of the RKHS parameterization and the kernel distance, respectively. The results are for two representative patients from our database. The average mesh size $\bar{h}_1$ of the target shape is \fnum{1.160142e+00}. The results correspond to the first, second and third block in \tabref{t:xpat-NaAd-OnGe-probe-sigma-tauS-tauV}. From top to bottom (row 1 through 3) we report results for $\tau_v=3$, $\tau_v=4$ and $\tau_v=6$. Each plot shows results for $\tau_s \in \{0.75,1,2,4,6\}$ (see legend). The first column shows the trend of the Hausdorff distance with respect to the number of iterations (normalized with respect to the initial value). Notice that we do not minimize the Hausdorff distance in our optimization but the kernel distance in \eqref{e:discretedist}. The second column shows the update of the Hausdorff distance per iteration. The third and fourth column shows the norm of the primal and dual residual versus the iteration count. These residuals are normalized with respect to the residual obtained at iteration 2.\label{f:t:xpat-NaAd-OnGe-probe-sigma-tauS-tauV}}
\end{figure}

We report additional results for six additional patients selected from our database for the experiments performed in \secref{s:bandwithselection}. We refer to \secref{s:bandwithselection} for more details about the setup of these numerical experiments. We report the final Hausdorff distance, the norm (absolute and relative) of the primal and dual residual and the runtime for different choices of the scaling parameters $\tau_s$ and $\tau_v$ for the bandwidth of the Gaussian kernels used in our kernel distance and the parameterization in a RKHS in \tabref{t:xpat-ShCr-ShCh-probe-sigma-tauS-tauV}, \tabref{t:xpat-WaSh-YaSa-probe-sigma-tauS-tauV}, and \tabref{t:xpat-NaAd-OnGe-probe-sigma-tauS-tauV}, respectively. Each table contains results for a different pair of patients. In addition, we illustrate the convergence behavior of our solver in \figref{f:xpat-ShCr-ShCh-probe-sigma-tauS-tauV}, \figref{f:xpat-WaSh-YaSa-probe-sigma-tauS-tauV} and \figref{t:xpat-NaAd-OnGe-probe-sigma-tauS-tauV} associated with a subset of these experiments. Again, each figure corresponds to a different patient. The layout for the tables as well as the figures is the same as in \secref{s:bandwithselection}.

\begin{table}
\caption{Performance for different choices of the scaling parameters $\tau_v$ and $\tau_s$ that control the bandwidth $\sigma_v$ and $\sigma_s$ of the RKHS parameterization and the kernel distance, respectively. The results concern the diffeomorphic matching of two mitral valve surfaces acquired at midsystole from two distinct patients (representative of our patients data base). The average mesh size $\bar{h}_1$ of the target shape is \fnum{1.239477e+00}. We select $\tau_v$ in $\{3,4,6,8\}$. We select $\tau_s$ in $\{\nicefrac{3}{4},1,2,4,6\}$. The initial (censored) Hausdorff distance for this dataset is \fnum{4.565145e+00}. We report (from left to right) the final Hausdorff distance (after matching; absolute value and percentage of initial value in brackets), the (relative) norm of the primal residual, the (relative) norm of the dual residual, and the runtime of the solver (in seconds).\label{t:xpat-ShCr-ShCh-probe-sigma-tauS-tauV}}
\centering\small
\begin{tabular}{lllllr}\toprule
$\tau_v$ & $\tau_s$   & {\bf final distance} (\%) & \multicolumn{2}{l}{\bf residuals}   & {\bf runtime} \\
         &            &                     & primal (relative) & dual (relative) \\\midrule
3 & 0.75 & \fnum{2.282962e+00} (\inum{50.01}) & \snum{1.318694e-02} (\snum{6.860209e-04}) & \snum{1.010133e+00} (\snum{5.254987e-02}) & \inum{107.682910} \\
  & 1    & \fnum{2.221254e+00} (\inum{48.66}) & \snum{6.406709e-03} (\snum{3.049397e-04}) & \snum{7.974050e-01} (\snum{3.795403e-02}) & \inum{89.391413} \\
  & 2    & \fnum{3.041994e+00} (\inum{66.64}) & \snum{3.928970e-03} (\snum{1.875003e-04}) & \snum{6.259779e-01} (\snum{2.987323e-02}) & \inum{67.086934} \\
  & 4    & \fnum{4.078304e+00} (\inum{89.34}) & \snum{3.038588e-03} (\snum{2.890389e-04}) & \snum{4.990895e-01} (\snum{4.747478e-02}) & \inum{63.531854} \\
  & 6    & \fnum{4.360767e+00} (\inum{95.52}) & \snum{2.012641e-03} (\snum{4.170102e-04}) & \snum{3.232280e-01} (\snum{6.697140e-02}) & \inum{60.876852} \\
\midrule
4 & 0.75 & \fnum{1.634233e+00} (\inum{35.80}) & \snum{9.008702e-03} (\snum{4.686574e-04}) & \snum{9.275919e-01} (\snum{4.825588e-02}) & \inum{101.220236} \\
  & 1    & \fnum{1.752065e+00} (\inum{38.38}) & \snum{6.415600e-03} (\snum{3.053628e-04}) & \snum{8.217829e-01} (\snum{3.911434e-02}) & \inum{81.714471} \\
  & 2    & \fnum{2.452084e+00} (\inum{53.71}) & \snum{5.359213e-03} (\snum{2.557551e-04}) & \snum{7.179372e-01} (\snum{3.426176e-02}) & \inum{67.096352} \\
  & 4    & \fnum{3.802236e+00} (\inum{83.29}) & \snum{4.197686e-03} (\snum{3.992956e-04}) & \snum{5.928163e-01} (\snum{5.639034e-02}) & \inum{62.158182} \\
  & 6    & \fnum{4.217449e+00} (\inum{92.38}) & \snum{2.417339e-03} (\snum{5.008620e-04}) & \snum{4.017875e-01} (\snum{8.324859e-02}) & \inum{60.273473} \\
\midrule
6 & 0.75 & \fnum{1.337887e+00} (\inum{29.31}) & \snum{7.736206e-03} (\snum{4.024587e-04}) & \snum{9.157639e-01} (\snum{4.764056e-02}) & \inum{91.521943} \\
  & 1    & \fnum{1.435812e+00} (\inum{31.45}) & \snum{6.670047e-03} (\snum{3.174737e-04}) & \snum{8.546688e-01} (\snum{4.067961e-02}) & \inum{79.450615} \\
  & 2    & \fnum{2.040929e+00} (\inum{44.71}) & \snum{6.152043e-03} (\snum{2.935909e-04}) & \snum{7.817732e-01} (\snum{3.730817e-02}) & \inum{67.817065} \\
  & 4    & \fnum{3.553041e+00} (\inum{77.83}) & \snum{4.750529e-03} (\snum{4.518836e-04}) & \snum{6.984289e-01} (\snum{6.643650e-02}) & \inum{57.793488} \\
  & 6    & \fnum{4.077730e+00} (\inum{89.32}) & \snum{2.281732e-03} (\snum{4.727647e-04}) & \snum{4.831260e-01} (\snum{1.001016e-01}) & \inum{57.364824} \\
\midrule
8 & 0.75 & \fnum{1.347715e+00} (\inum{29.52}) & \snum{8.659638e-03} (\snum{4.504982e-04}) & \snum{9.938346e-01} (\snum{5.170201e-02}) & \inum{88.565961} \\
  & 1    & \fnum{1.423136e+00} (\inum{31.17}) & \snum{7.629370e-03} (\snum{3.631345e-04}) & \snum{9.510106e-01} (\snum{4.526518e-02}) & \inum{79.702667} \\
  & 2    & \fnum{2.058415e+00} (\inum{45.09}) & \snum{5.886876e-03} (\snum{2.809364e-04}) & \snum{7.940777e-01} (\snum{3.789537e-02}) & \inum{66.541629} \\
  & 4    & \fnum{3.512109e+00} (\inum{76.93}) & \snum{4.642316e-03} (\snum{4.415901e-04}) & \snum{7.048431e-01} (\snum{6.704663e-02}) & \inum{58.084464} \\
  & 6    & \fnum{4.076401e+00} (\inum{89.29}) & \snum{1.998382e-03} (\snum{4.140559e-04}) & \snum{4.896961e-01} (\snum{1.014628e-01}) & \inum{54.157883} \\
\bottomrule
\end{tabular}
\end{table}

\begin{table}
\caption{Performance for different choices of the scaling parameters $\tau_v$ and $\tau_s$ that control the bandwidth $\sigma_v$ and $\sigma_s$ of the RKHS parameterization and the kernel distance, respectively. The results are for two representative patients from our database. The average mesh size $\bar{h}_1$ of the target shape is \fnum{1.076166e+00}. We select $\tau_v$ in $\{3,4,6,8\}$. We select $\tau_s$ in $\{\nicefrac{3}{4},1,2,4,6\}$. The initial (censored) Hausdorff distance for this dataset is \fnum{3.235954e+00}. We report (from left to right) the final Hausdorff distance (after matching; absolute value and percentage of initial value in brackets), the (relative) norm of the primal residual, the (relative) norm of the dual residual, and the runtime of the solver (in seconds).\label{t:xpat-WaSh-YaSa-probe-sigma-tauS-tauV}}
\centering\small
\begin{tabular}{lllllr}\toprule
$\tau_v$ & $\tau_s$   & {\bf final distance} (\%) & \multicolumn{2}{l}{\bf residuals}   & {\bf runtime} \\
         &            &                     & primal (relative) & dual (relative) \\\midrule
3 & 0.75 & \fnum{9.616103e-01} (\inum{29.72}) & \snum{2.994087e-03} (\snum{1.861310e-04}) & \snum{5.328101e-01} (\snum{3.312277e-02}) & \inum{94.079729} \\
  & 1    & \fnum{1.054698e+00} (\inum{32.59}) & \snum{2.739112e-03} (\snum{1.591885e-04}) & \snum{5.076694e-01} (\snum{2.950413e-02}) & \inum{81.451607} \\
  & 2    & \fnum{1.689346e+00} (\inum{52.21}) & \snum{2.184564e-03} (\snum{1.344331e-04}) & \snum{3.912323e-01} (\snum{2.407554e-02}) & \inum{73.702778} \\
  & 4    & \fnum{2.520968e+00} (\inum{77.90}) & \snum{2.453746e-03} (\snum{2.754564e-04}) & \snum{3.458618e-01} (\snum{3.882628e-02}) & \inum{65.766635} \\
  & 6    & \fnum{2.755975e+00} (\inum{85.17}) & \snum{1.902508e-03} (\snum{4.442333e-04}) & \snum{3.254364e-01} (\snum{7.598901e-02}) & \inum{63.211539} \\
\midrule
4 & 0.75 & \fnum{8.430208e-01} (\inum{26.05}) & \snum{4.284224e-03} (\snum{2.663339e-04}) & \snum{6.414984e-01} (\snum{3.987951e-02}) & \inum{88.064567} \\
  & 1    & \fnum{8.862896e-01} (\inum{27.39}) & \snum{3.768049e-03} (\snum{2.189870e-04}) & \snum{5.992679e-01} (\snum{3.482754e-02}) & \inum{78.592223} \\
  & 2    & \fnum{1.347934e+00} (\inum{41.65}) & \snum{2.878026e-03} (\snum{1.771071e-04}) & \snum{4.572368e-01} (\snum{2.813731e-02}) & \inum{75.175425} \\
  & 4    & \fnum{2.229514e+00} (\inum{68.90}) & \snum{3.076610e-03} (\snum{3.453788e-04}) & \snum{3.893044e-01} (\snum{4.370313e-02}) & \inum{64.571813} \\
  & 6    & \fnum{2.656403e+00} (\inum{82.09}) & \snum{2.585690e-03} (\snum{6.037554e-04}) & \snum{4.043116e-01} (\snum{9.440627e-02}) & \inum{60.919303} \\
\midrule
6 & 0.75 & \fnum{8.904025e-01} (\inum{27.52}) & \snum{6.827124e-03} (\snum{4.244163e-04}) & \snum{7.550348e-01} (\snum{4.693764e-02}) & \inum{91.073585} \\
  & 1    & \fnum{9.427237e-01} (\inum{29.13}) & \snum{4.935048e-03} (\snum{2.868093e-04}) & \snum{6.937871e-01} (\snum{4.032070e-02}) & \inum{81.469029} \\
  & 2    & \fnum{1.307830e+00} (\inum{40.42}) & \snum{3.370887e-03} (\snum{2.074366e-04}) & \snum{5.221418e-01} (\snum{3.213141e-02}) & \inum{69.458165} \\
  & 4    & \fnum{1.951229e+00} (\inum{60.30}) & \snum{3.435081e-03} (\snum{3.856206e-04}) & \snum{4.425240e-01} (\snum{4.967754e-02}) & \inum{61.391864} \\
  & 6    & \fnum{2.459062e+00} (\inum{75.99}) & \snum{2.996709e-03} (\snum{6.997279e-04}) & \snum{4.685647e-01} (\snum{1.094093e-01}) & \inum{57.993948} \\
\midrule
8 & 0.75 & \fnum{1.009600e+00} (\inum{31.20}) & \snum{8.358934e-03} (\snum{5.196431e-04}) & \snum{7.992024e-01} (\snum{4.968337e-02}) & \inum{97.204811} \\
  & 1    & \fnum{1.053831e+00} (\inum{32.57}) & \snum{6.433122e-03} (\snum{3.738726e-04}) & \snum{7.402340e-01} (\snum{4.302004e-02}) & \inum{83.965953} \\
  & 2    & \fnum{1.359764e+00} (\inum{42.02}) & \snum{3.321021e-03} (\snum{2.043680e-04}) & \snum{5.351628e-01} (\snum{3.293269e-02}) & \inum{68.959240} \\
  & 4    & \fnum{1.886035e+00} (\inum{58.28}) & \snum{3.310927e-03} (\snum{3.716831e-04}) & \snum{4.451063e-01} (\snum{4.996742e-02}) & \inum{60.693723} \\
  & 6    & \fnum{2.353753e+00} (\inum{72.74}) & \snum{2.983355e-03} (\snum{6.966099e-04}) & \snum{4.723249e-01} (\snum{1.102873e-01}) & \inum{58.287894} \\
\bottomrule
\end{tabular}
\end{table}

\begin{table}
\caption{Performance for different choices of the scaling parameters $\tau_v$ and $\tau_s$ that control the bandwidth $\sigma_v$ and $\sigma_s$ of the RKHS parameterization and the kernel distance, respectively. The results are for two representative patients from our database. The average mesh size $\bar{h}_1$ of the target shape is \fnum{1.160142e+00}. We select $\tau_v$ in $\{3,4,6,8\}$. We select $\tau_s$ in $\{\nicefrac{3}{4},1,2,4,6\}$. The initial (censored) Hausdorff distance for this dataset is \fnum{5.607648e+00}. We report (from left to right) the final Hausdorff distance (after matching; absolute value and percentage of initial value in brackets), the (relative) norm of the primal residual, the (relative) norm of the dual residual, and the runtime of the solver (in seconds).\label{t:xpat-NaAd-OnGe-probe-sigma-tauS-tauV}}
\centering\small
\begin{tabular}{lllllr}\toprule
$\tau_v$ & $\tau_s$   & {\bf final distance} (\%) & \multicolumn{2}{l}{\bf residuals}   & {\bf runtime} \\
         &            &                     & primal (relative) & dual (relative) \\\midrule
3 & 0.75 & \fnum{2.609566e+00} (\inum{46.54}) & \snum{2.574608e-02} (\snum{8.352317e-04}) & \snum{1.591545e+00} (\snum{5.163150e-02}) & \inum{108.800464} \\
  & 1    & \fnum{2.657062e+00} (\inum{47.38}) & \snum{1.138272e-02} (\snum{3.343344e-04}) & \snum{1.175546e+00} (\snum{3.452828e-02}) & \inum{92.787476} \\
  & 2    & \fnum{3.465720e+00} (\inum{61.80}) & \snum{7.098221e-03} (\snum{2.338866e-04}) & \snum{1.011519e+00} (\snum{3.332960e-02}) & \inum{71.657277} \\
  & 4    & \fnum{5.039531e+00} (\inum{89.87}) & \snum{3.459248e-03} (\snum{3.696201e-04}) & \snum{6.576166e-01} (\snum{7.026622e-02}) & \inum{65.715832} \\
  & 6    & \fnum{5.558619e+00} (\inum{99.13}) & \snum{2.093782e-03} (\snum{6.691898e-04}) & \snum{3.243178e-01} (\snum{1.036546e-01}) & \inum{61.161741} \\
\midrule
4 & 0.75 & \fnum{1.946453e+00} (\inum{34.71}) & \snum{1.370692e-02} (\snum{4.446679e-04}) & \snum{1.252759e+00} (\snum{4.064091e-02}) & \inum{102.311561} \\
  & 1    & \fnum{2.058129e+00} (\inum{36.70}) & \snum{9.589375e-03} (\snum{2.816602e-04}) & \snum{1.106125e+00} (\snum{3.248923e-02}) & \inum{88.192574} \\
  & 2    & \fnum{3.065749e+00} (\inum{54.67}) & \snum{8.163451e-03} (\snum{2.689860e-04}) & \snum{1.034979e+00} (\snum{3.410260e-02}) & \inum{69.874286} \\
  & 4    & \fnum{4.773156e+00} (\inum{85.12}) & \snum{4.467787e-03} (\snum{4.773823e-04}) & \snum{8.722164e-01} (\snum{9.319618e-02}) & \inum{60.779706} \\
  & 6    & \fnum{5.516478e+00} (\inum{98.37}) & \snum{2.236132e-03} (\snum{7.146861e-04}) & \snum{3.818292e-01} (\snum{1.220358e-01}) & \inum{57.587639} \\
\midrule
6 & 0.75 & \fnum{1.406057e+00} (\inum{25.07}) & \snum{1.350149e-02} (\snum{4.380036e-04}) & \snum{1.314412e+00} (\snum{4.264099e-02}) & \inum{93.420104} \\
  & 1    & \fnum{1.624450e+00} (\inum{28.97}) & \snum{9.597825e-03} (\snum{2.819084e-04}) & \snum{1.189118e+00} (\snum{3.492690e-02}) & \inum{85.307116} \\
  & 2    & \fnum{2.589942e+00} (\inum{46.19}) & \snum{8.539388e-03} (\snum{2.813731e-04}) & \snum{1.097518e+00} (\snum{3.616327e-02}) & \inum{69.632610} \\
  & 4    & \fnum{4.572970e+00} (\inum{81.55}) & \snum{4.945769e-03} (\snum{5.284546e-04}) & \snum{1.006337e+00} (\snum{1.075269e-01}) & \inum{57.449638} \\
  & 6    & \fnum{5.437616e+00} (\inum{96.97}) & \snum{1.941372e-03} (\snum{6.204785e-04}) & \snum{4.204996e-01} (\snum{1.343951e-01}) & \inum{53.617893} \\
\midrule
8 & 0.75 & \fnum{1.346575e+00} (\inum{24.01}) & \snum{1.697474e-02} (\snum{5.506796e-04}) & \snum{1.343328e+00} (\snum{4.357907e-02}) & \inum{93.490618} \\
  & 1    & \fnum{1.585189e+00} (\inum{28.27}) & \snum{1.108558e-02} (\snum{3.256070e-04}) & \snum{1.218573e+00} (\snum{3.579208e-02}) & \inum{85.495666} \\
  & 2    & \fnum{2.442877e+00} (\inum{43.56}) & \snum{8.431196e-03} (\snum{2.778082e-04}) & \snum{1.095663e+00} (\snum{3.610212e-02}) & \inum{68.242436} \\
  & 4    & \fnum{4.522326e+00} (\inum{80.65}) & \snum{4.902813e-03} (\snum{5.238647e-04}) & \snum{1.000835e+00} (\snum{1.069390e-01}) & \inum{56.992847} \\
  & 6    & \fnum{5.406784e+00} (\inum{96.42}) & \snum{1.617577e-03} (\snum{5.169908e-04}) & \snum{4.133701e-01} (\snum{1.321165e-01}) & \inum{48.545881} \\
\bottomrule
\end{tabular}
\end{table}

\end{appendix}

\end{document}